\documentclass[a4paper,10pt]{amsart}
\usepackage{tikz}
\usetikzlibrary{matrix,arrows,decorations.pathmorphing}
\usepackage{tikz-cd}
\usepackage{amsfonts, amssymb, amsthm, amsmath}  
\usepackage{hyperref}
\usepackage{graphicx}
\newtheorem{thm}{Theorem}  

\newtheorem{cor}[thm]{Corollary}  
\newtheorem{lemma}[thm]{Lemma}  
\newtheorem{remark}[thm]{Remark}  
\newtheorem{defn}[thm]{Definition}  
\newtheorem{prop}[thm]{Proposition}  
\newtheorem{claim}[thm]{Claim}  
  
\numberwithin{thm}{section}  
\def\pf{\noindent\emph{Proof: }}  
\def\stop{\hfill$\square$}  
\def\s{\mathfrak s}

\newcommand{\Charts}{{\mathcal Charts}}
\DeclareMathOperator{\Poly}{Poly}

\providecommand{\poly}[4]{\Poly^{#2}_{#1}(#3,#4)}
\providecommand{\ov}[1]{\hspace{-.1cm}\downarrow_{#1}}

\providecommand{\Y}{\mathcal Y}

\providecommand{\totl}[1]{\ensuremath{\lceil #1\rceil }}
\providecommand{\totb}[1]{\ensuremath{\underline{ #1}}}
\DeclareMathOperator{\End}{End}

\newcommand{\exte}{\subset_{e}}
\newcommand{\rexte}{\ {}_{e}\!\!\supset}

\newcommand{\ex}{\bold}
\providecommand {\e}[1]{\mathfrak t^{#1}}
\newcommand{\co}{\mathfrak{o}}

\providecommand{\C}[2]{\ensuremath {C^{#1,\underline{#2}}}}

\newcommand{\Msw}{\mathcal M^{st}}

\DeclareMathOperator{\id}{id}

\newcommand{\dbar}{\bar{\partial}}

\providecommand{\et}[2]{\ensuremath{\bold T^{#1}_{#2}}}
\providecommand{\lrb}[1]{\ensuremath{\left(#1\right)}}
\providecommand{\abs}[1]{\left\lvert #1\right\rvert}

\author{Brett Parker   }
  
\title{Integral counts of pseudo-holomorphic curves}

\begin{document}

\thanks{This work was supported by ARC grant DP1093094.}

\begin{abstract} In \cite{FOinteger}, Fukaya and Ono outlined a way of counting pseudo-holomorphic curves in a general compact symplectic manifold to obtain integer valued invariants. This paper contains the details of Fukaya and Ono's suggested construction for any compact symplectic manifold and a large class of exploded manifolds. 
\end{abstract}

\maketitle

\tableofcontents

\newpage

\section{Introduction}

Because a curve $f$ with a group $\abs G$ of automorphisms must be counted with weight $1/\abs{G}$, Gromov-Witten invariants of a symplectic manifold $X$ take values in the rational numbers. In \cite{FOinteger}, Fukaya and Ono outlined a way to construct an integer part of Gromov-Witten invariants. Fukaya and Ono's integer counts of holomorphic curves  remove curves with nontrivial automorphism group from the count as described below: 

\

 The section $\dbar$, whose intersection with $0$ defines holomorphic curves, is modified to a section $s$ which is complex in directions normal to curves with nontrivial automorphism group. For a generic such normally complex section $s$, the resulting perturbed moduli space, $\mathcal N:= s^{-1}(0)$  need not be smooth, but is compact and breaks up into smooth strata depending on the number of nodes a curve has and the automorphism group of the curve. In particular, the curves of genus $g$ with no automorphisms, nodes, or marked points, representing a homotopy class $\beta\in [\Sigma_{g},X]/MCG(\Sigma_{g})$ form a smooth, oriented moduli space $\mathcal N_{g,\beta}\subset \mathcal N$. This automorphism free part of the moduli space $\mathcal N_{g,\beta}$  has the expected dimension $d$ of the moduli space of holomorphic curves and  $\overline{\mathcal N_{g,\beta}}\setminus \mathcal N_{g,\beta}$ is a finite union of smooth manifolds with dimension at most $d-2$. Integer invariants can be defined from the homology class represented by $\overline{\mathcal N_{g,\beta}}$ as follows: 
 
 If the expected dimension $d=0$, then the signed number of curves
 $n_{g,\beta}$
 in $\mathcal N_{g,\beta}$ is a symplectic invariant which does not depend on the choices used to define the normally complex perturbation $s$. In the more general case that the expected dimension $d\geq 0$,  we can count curves in $\mathcal N_{g,\beta}$ passing through generic oriented submanifolds $\alpha_{1},\dotsc, \alpha_{k}\subset X$ with codimensions summing to $d+2k$. The signed number of curves in $\mathcal N_{g,\beta}$ with points $x_{i}$ sent to $\alpha_{i}$ is an integer  invariant
 \[n_{g,\beta}(\alpha_{1},\dotsc,\alpha_{k})\]
which depends only on   $g$, $\beta$,  the cohomology classes Poincare dual to $\alpha_{i}$ in $H^{*}(X,\mathbb Z)$, and the deformation class of the symplectic form on $X$. 

\

As outlined in section \ref{future section}, there are several generalizations of the invariants $n_{g,\beta}$ which use some version of `normally complex sections'. In particular a version of gravitational descendants may easily be incorporated into the invariants $n_{g,\beta}$ by using generic normally complex sections of tautological line bundles.  
With a little more work, Floer homology can be defined using integer instead of rational coefficients, as outlined by Fukaya and Ono in \cite{FOinteger}. There should also be a version of `normally complex sections' which may be used to define an invariant which throws away contributions from  curves which factor through other holomorphic curves.

\

The proofs of this paper are carried out in the setting of curves in a family $\hat{\ex B}$ of exploded manifolds. There are invariants analogous to $n_{g,\beta}(\alpha_{1},\dotsc,\alpha_{k})$ which are defined for exploded manifolds, however because the invariance depends on a `codimension $2$ singularities' argument, these invariants may change in families of exploded manifolds. The reason that exploded manifolds must be used in this paper is that the analysis of families of curves in which nodes form is easier and more natural in the setting of exploded manifolds. Using exploded manifolds in \cite{evc} allowed the construction of a particularly concrete Kuranishi structure with a type of almost complex structure. This almost complex structure was the key extra ingredient required to construct Fukaya and Ono's integer invariants. 

\subsection{Outline of the construction of normally complex sections}

\

Let $\Msw(X)$ indicate the moduli stack of stable $\C\infty1$ curves in a symplectic manifold $X$, as defined in \cite{evc}. This is a stack which may be thought of as describing families of (not-necessarily holomorphic) stable curves
\[\begin{tikzcd}\ex C\dar \rar{\hat f}& X
\\ \ex F \end{tikzcd}\]
where  $\ex C\longrightarrow \ex F$ is a $\C\infty1$ family of curves in the category of exploded manifolds, and $\hat f$ is a $\C\infty1$ map.  Without using exploded manifolds, the family $\ex C\longrightarrow \ex F$ may loosely be thought of as a family of Riemann surfaces in which nodes are allowed to form, and the regularity $\C\infty1$ may be thought of as meaning smooth away from nodes, and a kind of regularity at nodes which is as good as smooth for all practical purposes.  The reader unfamiliar with stacks should think of $\Msw(X)$ as a Frechet orbifold or polyfold representing a moduli space of (not-necessarily holomorphic) stable curves in $X$--- so a family of curves corresponds to a map $\ex F\longrightarrow \Msw(X)$.

Choose an  almost complex structure $J$ on $X$ so that the symplectic form is positive on $J$-holomorphic planes. The $\dbar$ equation defines a section 
\[\dbar:\Msw(X)\longrightarrow \Y\]
of some sheaf $\Y$ defined over $\Msw(X)$, so that the moduli stack of $J$-holomorphic curves correspond to the intersection of $\dbar$ with $0$. (If you are thinking of $\Msw(X)$ as a Frechet orbifold, then think of $\Y$ as a complex,  infinite dimensional, Frechet vector bundle over $\Msw(X)$.)

At any curve $f\in \Msw(X)$, there is a natural tangent space $T_{f}\Msw(X)$ defined in \cite{evc}. If $f$ is a nodal curve, then corresponding to variations of infinitesimal gluing data at each node, there is a copy of $\mathbb C\subset T_{f}\Msw$. The quotient of $T_{f}\Msw$ by this gluing data corresponds to the space of variations of $f$ which do not smooth any of the nodes.

 If $f$ is $J$-holomorphic, then the derivative of $\dbar$ at $f$ makes sense and is a Fredholm map:
\[D\dbar_{f}:T_{f}\Msw(X)\longrightarrow \Y(f)\]
At holomorphic curves $f$, the tangent space $T_{f}\Msw$ has a natural complex structure, but $D\dbar_{f}$ is not usually complex unless $J$ is integrable.  If $f$ has an automorphism group $G$, then $G$ acts naturally on $T_{f}\Msw(X)$ and $\Y(f)$ so that $D\dbar_{f}$ is $G$-equivariant. The linear homotopy of $D\dbar_{f}$ to its complex linear part $D\dbar_{f}^{\mathbb C}$ is Fredholm and $G$-equivariant. There therefore exists a finite dimensional, complex, $G$-invariant subspace $V(f)\subset \Y(f)$ so that the linear homotopy of $D\dbar_{f}$ to its complex linear part is transverse to $V(f)$.

As proved in \cite{evc}, $V(f)$ may be extended to a nice vectorbundle $V\subset \Y$ defined on a neighborhood $\mathcal U$ of $f\in\Msw$ so that the moduli stack of curves in $\mathcal U$ with $\dbar$ contained in $V$ is represented by $\hat f/G$, where $\hat f$ is some $\C\infty1$ family of curves with a group $G$ of automorphisms. The data $(\mathcal U,V,\hat f/G)$ is called an embedded Kuranishi chart.  

Let the domain of $\hat f$ be $\ex C\longrightarrow \ex F$. The family of curves $\hat f$ corresponds to a map $\ex F\longrightarrow \mathcal U\subset\Msw$. Use the notation $V$ to indicate the pullback of $V$ under the inclusion $\ex F\longrightarrow \mathcal U$. Note that $\dbar\hat f$ defines a $G$-equivariant section of $V$ over $\ex F$. We therefore have the following finite dimensional local model for the $\dbar$ equation
\begin{equation}\label{lm}\begin{tikzcd}\arrow[loop left]{}{G} V\dar
\\ \arrow[loop left]{}{G}\ex F\uar[bend right][swap]{\dbar}
\end{tikzcd}\end{equation}
 
 In order to define a virtual fundamental class for the moduli stack of holomorphic curves, we need some way of defining the virtual intersection of $\dbar$ with $0$. A generic $G$-equivariant section of $V$ may not be transverse to the zero section, and its intersection with $0$ may have codimension $1$ boundary. Fukaya and Ono overcame this problem in \cite{FO} by perturbing $\dbar$ to a $G$-equivariant weighted branched  section of $V$. The resulting transverse intersection with $0$ is only good for defining a virtual fundamental class over $\mathbb Q$.

We may do better by incorporating information from the complex structure of $T_{f}\Msw$ into the local model from figure (\ref{lm}). The linear homotopy of $D\dbar_{f}$ to $D\dbar_{f}^{\mathbb C}$ defines a homotopy class of $G$-equivariant complex structure on $T_{f}\ex F$. There is a canonical complex structure on $T_{f}\Msw$ only if $f$ is holomorphic, however in \cite{evc}, it is proved that so long as $\mathcal U$ is chosen small enough,  there is a canonical homotopy class of $G$-equivariant almost complex structure on $\ex F$.

For the purposes of exposition, assume for the moment that we have a $G$-equivariant complex structure on $\ex F$ and have chosen a $G$-equivariant holomorphic structure on the complex vector bundle $V$ over $\ex F$. We may then define a normally complex section of $V$ to be a section which is a finite sum of complex, $G$-equivariant sections  times smooth, $G$-invariant $\mathbb C$-valued functions. 

 Note that at a point in $\ex F$ around which $G$ acts freely, normally complex sections are the same as smooth sections of $V$. Similarly, on a strata of $(\ex F,G)$ with isotropy group $H$, normally complex sections only use the complex structure in directions normal to the strata, and do not see the complex structure in the $H$-invariant directions of the strata. The actual construction of normally complex sections uses a complex structure normal to each strata of $(\ex F,G)$. In this paper, we make a global construction of such complex structures normal to strata of $(\ex F,G)$ so that there is a compatible notion of normally complex sections on a collection of Kuranishi charts covering the entire moduli stack of holomorphic curves.

Because of the flexibility introduced by allowing multiplication by smooth (not holomorphic) $\mathbb C$-valued functions, it is easy to achieve a kind of transversality using equivariant normally complex sections.   Standard invariant theory gives that normally complex sections are locally generated by a finite number of complex algebraic sections. As complex algebraic singularities are codimension $2$ and nice, the intersection of a generic $G$-equivariant normally complex section with the zero section will have nice, codimension $2$ singularities. As outlined by Fukaya and Ono in \cite{FOinteger}, we may define integer holomorphic curve invariants by perturbing $\dbar$ to a generic normally complex section, then counting curves with no automorphisms as described in the previous section. 

\subsection{Example: how rigid curves with automorphisms are counted}

\

Suppose that $f$ is a holomorphic curve so that that $D\dbar_{f}$ is an isomorphism.  In this case, $f$  has a well defined contribution to $n_{g,\beta}$ which we shall now compute in terms of $D\dbar_{f}$ by considering a single Kuranishi chart containing $f$.

Let the automorphism group of $f$ be $G$. The only case in which $D\dbar_{f}$ can be an isomorphism is if the domain of $f$ is a smooth (node free) curve on which $G$ acts freely.  As calculated in Remark \ref{multiple cover}, the index of $D\dbar_{f}^{\mathbb C}$ is $0$ in complex, $G$-equivariant $K$-theory in this case if and only if the index of $D\dbar_{f}$ is $0$ in ordinary $K$-theory.  For a generic choice of almost complex structure, $D\dbar_{f_{0}}$ is an isomorphism if $f_{0}$ is   somewhere injective, has a smooth domain, and the index of $D\dbar_{f_{0}}$ is $0$.  It can also be expected that $D\dbar_{f}$ is an isomorphism for any unbranched multiple cover $ f$ of such a curve.

We shall compute the contribution of such unbranched multiple covers to $n_{g,\beta}$. The calculation of the contribution of all the other possible curves which factor through $f_{0}$ is significantly more delicate because such curves come in a moduli stack with dimension greater than $0$.

\

As $D\dbar_{f}$ is $G$-equivariant, it splits into factors corresponding to the irreducible real representations $\rho$ of $G$. \[D\dbar_{f}=\oplus_{\rho}D\dbar_{f,\rho}\]
Although $D\dbar_{f,\rho}$ it not itself complex, the domain and range of $D\dbar_{f,\rho}$ both have complex structure.  The fact that $[D\dbar_{f}^{\mathbb C}]$ is $0$ in complex $G$-equivariant $K$-theory implies that $D\dbar_{f,\rho}$ is homotopic through $G$-equivariant Fredholm maps to a complex linear isomorphism. We shall need to know if $D\dbar_{f,\rho}$ is isotopic through linear $G$-equivariant isomorphisms to a complex isomorphism.

The endomorphism ring, $\End(\rho)$, of each irreducible real representation $\rho$  is either $\mathbb R$, $\mathbb C$, or $\mathbb H$, because Schur's lemma implies that it is a division ring over $\mathbb R$.  If  $\End(\rho)$ is isomorphic to $\mathbb C$ or $\mathbb H$, then $D\dbar_{f,\rho}$ is isotopic through $G$-equivariant isomorphisms to a complex linear isomorphism because the $G$-equivariant Fredholm maps which are not isomorphisms have codimension $ 2$ and $4$ respectively.  

In the remaining case that $\End(\rho)=\mathbb R$,  we shall define $\sigma_{f,\rho}$ below as a kind of generalized sign of the determinant of $D\dbar_{f,\rho}$. This generalized sign will live in the following ring:

\begin{defn} Let $\mathbb Z[I]_{I\leq G}$ indicate the ring over the integers generated by subgroups of $G$ with multiplication corresponding to intersection of subgroups. (Note that the element $G$ acts as the multiplicative identity in this ring, so we shall use $1$ to mean $1G$.) Let $<\cdot,\cdot>$ indicate the inner product on $\mathbb Z[I]_{I\leq G}$ with orthonormal basis given by the set of subgroups $I\leq G$, so 
\[<\sum_{I} a_{I}I,\sum_{I} b_{I}I>=\sum_{I}a_{I}b_{I}\] 
\end{defn}

The multiplication in $Z[I]_{I\leq G}$ is designed to reflect the fact that the product of two points with isotropy groups $I_{1}$ and $I_{2}$ will have isotropy group $I_{1}\cap I_{2}$. 

\begin{defn}
If $D\dbar_{f,\rho}$ is isotopic to a complex isomorphism through $G$-equivariant isomorphisms,  define 
\[\sigma_{f,\rho}:=1\in \mathbb Z[I]_{I\leq G}\]
In the remaining case that $D\dbar_{f,\rho}$ is not isotopic to a complex isomorphism through $G$-equivariant isomorphisms, then define $\sigma_{f,\rho}\in \mathbb Z[I]_{I\leq G}$ as the unique element which satisfies the following equation for all subgroups $H\leq G$:
\[<\sigma_{\rho}H,H>=(-1)^{\dim \rho^{H}}\]
where $\dim \rho^{H}$ indicates the dimension of the subspace of the representation $\rho$ on which $H$ acts trivially. 

Then define 
\[\sigma_{f}:=\prod_{\rho}\sigma_{f,\rho}\]
\end{defn}

In the special case that $\rho$ is the trivial representation, and $D\dbar_{f,\rho}$ is not isotopic to a complex map,  we regard $\rho$ as given by $G$ acting trivially on $\mathbb R$. Then $\dim \rho^{H}=1$ for all $H$, so $\sigma_{f,\rho}=- 1 $. In other words, if $\rho$ is the trivial representation, then $\sigma_{f,\rho}$ is the sign of the determinant of $D\dbar_{f,\rho}$ times the identity, $1\in\mathbb Z[I]_{I\leq G}$.

The coefficient $<\sigma_{f,\rho},I>$ of $I$ in $\sigma_{f,\rho}$ may be computed from knowledge of $\dim \rho^{H}$ using M\"obius inversion, so the above equations uniquely determine $\sigma_{f,\rho}$. We shall see that $<\sigma_{f,\rho},I>$ counts the number of zeros of a normally complex perturbation of $D_{f,\rho}$ which have isotropy group $I$ (so $<\sigma_{f,\rho}H,H>$ counts the number of zeros with isotropy group containing $H$.)  From this interpretation of $\sigma_{f,\rho}$ as given by a count of zeros of a normally complex perturbation, we shall be able to prove that $\sigma_{f,\rho}^{2}=1$, so $\sigma_{f}$ may be thought of as a generalized sign.

\begin{lemma}\label{f contribution}The contribution of $f$ to $n_{g,\beta}$ is the integer
\[\frac 1{\abs G}<\sigma_{f},\{1\}>\]
\end{lemma}

\pf

Consider a Kuranishi chart containing $f$ which gives the local model

\[\begin{tikzcd}\arrow[loop left]{}{G} V\dar
\\ \arrow[loop left]{}{G}\ex F\uar[bend right][swap]{\dbar}
\end{tikzcd}\]
We must modify $\dbar$ to a normally complex section $s$.
As specified in Definition \ref{c def} part \ref{c1}, our normally complex section $s$ must approximate $\dbar$ outside of a compactly contained neighborhood of $f\in \ex F$ in the sense that outside of this neighborhood, $\abs{\dbar - s}<\frac 12\abs \dbar$. In particular, as both $s$ and $\dbar$ are transverse to $0$, this implies that the algebraic number of zeros of $s$ is the same as the algebraic number of zeros of $\dbar:\ex F\longrightarrow V$, which is the sign of the determinant of $D\dbar_{f}$. (Recall that $D\dbar_{f}$ is equal to the sum of the derivative of $\dbar:\ex F\longrightarrow V$ with a $G$-equivariant isomorphism which is isotopic through $G$-equivariant isomorphisms to a complex linear map. We can therefore identify the sign of the determinant of $D\dbar_{f}$ with the sign of the determinant of the derivative of $\dbar$ at $f$. )

Given any subgroup $H\leq G$,  let $\ex F^{H}$ indicate the closed submanifold of $\ex F$ consisting of points which are fixed by $H$, and $V^{H}\longrightarrow \ex F^{H}$ be the sub vector bundle of $V$ restricted to $\ex F^{H}$ on which $H$ acts trivially. Again, both $s$ and $\dbar$ are transverse to $0$ when restricted to sections  $\ex F^{H}\longrightarrow V^{H}$, so they have the same algebraic number of zeros. The normally complex section $s$ is  complex in directions normal to $\ex F^{H}$, therefore the algebraic  number of zeros of $s$ which happen to be contained in $\ex F^{H}$ is the same as the algebraic number of zeros of $\dbar$ restricted to be a section of $\ex F^{H}\longrightarrow  V^{H}$, which is the sign of the determinant of $D\dbar_{f}^{H}$.   This observation suffices to compute the number of zeros of $s$ with any given isotropy group. 

In particular define 
\[\sigma:=\sum_{I\leq G}n_{I}I\]
where $n_{I}$ is the algebraic number of zeros of $s$ with isotropy group $I$. We shall show that $\sigma=\sigma_{f}$. For now, we have shown that
\[<\sigma H,H>=\text{ sign of the determinant of }D\dbar_{f}^{H}\]
Recalling that the product of $I_{1}$ with $I_{2}$ is $I_{1}\cap I_{2}$, it follows that for any two elements $\sigma_{i}\in \mathbb Z[I]_{I\leq G}$, 
\[<\sigma_{1}\sigma_{2}H,H>=<\sigma_{1}H,H><\sigma_{2}H,H>\]
The sign of the determinant  of $D\dbar_{f}^{H}$ is the product of the signs of the determinants of $D\dbar_{f,\rho}^{H}$. Therefore, to prove that $\sigma_{f}:=\prod_{\rho}\sigma_{f,\rho}=\sigma$, it suffices to prove that  $<\sigma_{f,\rho}H,H>$ is the sign of the determinant of $D\dbar_{f,\rho}^{H}$ for all $\rho$. 

If $D\dbar_{f,\rho}$ is isotopic to a complex map through $G$-equivariant isomorphisms, then the sign of the determinant of  $D\dbar_{f,\rho}^{H}$ is $1$. In this case $\sigma_{f,\rho}$ is defined to be $1$, so $<\sigma_{f,\rho}H,H>=1$ for all $H$ as required.

The remaining case to check is when $D\dbar_{f,\rho}$ is not isotopic through $G$-equivariant isomorphisms to a complex linear isomorphism. Let $V_{\rho}$ be the vector space corresponding to the irreducible real representation $\rho$. The only case in which $\sigma_{f,\rho}$ is not $1$ is the case in which the endomorphism ring of $V_{\rho}$ is equal to $\mathbb R$, and $V_{\rho}\oplus iV_{\rho}$ is an irreducible complex representation of $G$. Then the direct sum of $D\dbar_{f,\rho}$ with the complex conjugation map $V_{\rho}\oplus iV_{\rho}\longrightarrow V_{\rho}\oplus i V_{\rho}$ is isotopic through $G$-equivariant isomorphisms to a complex linear map, therefore the sign of the determinant of $D\dbar_{f,\rho}^{H}$ is equal to the sign of the determinant of the $H$-invariant part of the complex conjugation map $V_{\rho}\oplus iV_{\rho}\longrightarrow V_{\rho}\oplus i V_{\rho}$, which is $(-1)^{\dim V_{\rho}^{H}}$.
By definition,  $<\sigma_{\rho,f}H,H>=(-1)^{\dim V_{\rho}^{H}}$ in this case, therefore $<\sigma_{\rho,f}H,H>$ is equal to the sign of the determinant of $D\dbar_{f,\rho}^{H}$, as required. 

We have now proven that
\[<\sigma_{f}H,H>=<\sigma H,H>\]
for all $H\leq G$, therefore $\sigma_{f}=\sigma$. The contribution of $f$ to $n_{g,\beta}$ is 
\[\frac {n_{\{1\}}}{\abs G}=\frac {<\sigma \{1\},\{1\}>}{\abs G}=\frac{<\sigma_{f}\{1\},\{1\}>}{\abs G}\]
as required. Note that this really is an integer because $n _{\{1\}}$ is the number of zeros of $s$ on which $G$ acts freely, therefore $\abs G$ divides $n _{\{1\}}$.

\stop

\

We shall now consider how to calculate the contribution of all unbranched multiple covers  $f$ of a given curve $f_{0}$ under the assumption that $D\dbar_{f}$ is an isomorphism. (This assumption should hold if the almost complex structure is generic, and $D\dbar_{ f_{0}}$ is an isomorphism.) 

\

Let the domain of $f_{0}$ be $\Sigma$. The finite, connected multiple covers of $\Sigma$ correspond to finite index normal subgroups $K\leq \pi_{1}(\Sigma)$. Given such a finite index normal subgroup $K$, let $f_{K}$ indicate the corresponding cover of $f_{0}$. The automorphism group $G$ of a cover used above is equal to $\pi_{1}(\Sigma)/K$, so each irreducible $G$-representation $\rho$ corresponds to an irreducible $\pi_{1}(\Sigma)$ representation which we shall also refer to as $\rho$. The kernel $K_{\rho}$ of $\rho$ is a finite index normal subgroup of $\pi_{1}(\Sigma)$, and therefore corresponds to a cover $f_{K_{\rho}}$ of $f_{0}$. 

Denote by 
\[\sigma_{\rho}\in\mathbb Z[I]_{I\leq \pi_{1}(\Sigma)}\] the image of $\sigma_{f_{K_{\rho}},\rho}$ under the obvious homomorphism 
\[\mathbb Z[I]_{I\leq \pi_{1}(\Sigma)/{K_{\rho}}}\longrightarrow \mathbb Z[I]_{I\leq\pi_{1}(\Sigma)}\]

The fact that relates the contribution of all multiple covers of $f_{0}$ is the following: $D\dbar_{f_{K},\rho}$ is equal to the pullback of $D\dbar_{f_{K_{\rho}},\rho}$ under the covering map corresponding to the inclusion $K\leq K_{\rho}$. In particular, this implies that 
\[\sigma_{f_{K},\rho}=\pi_{K}\sigma_{\rho}\]
where 
\[\pi_{K}:\mathbb Z[I]_{I\leq \pi_{1}(\Sigma)}\longrightarrow\mathbb Z[I]_{I\leq\pi_{1}(\Sigma)/K}\]
is the  homomorphism which sends $I$ to $I/K$ if $K\leq I$, and which sends $I$ to $0$ otherwise. 
 Therefore,
\[\sigma_{f_{K}}=\pi_{K}\prod_{\rho\  \vert\  K_{\rho}\geq K } \sigma_{\rho}\]
The contribution of $f_{K}$ to $n_{g,\beta}$ is the coefficient of the trivial subgroup in $\sigma_{f_{K}}$ divided by $\abs{\pi_{1}(\Sigma)/K}$. The contribution of $f_{K}$ to $n_{g,\beta}$ is therefore also equal to the coefficient of $K$ in $\prod_{\rho\  \vert\  K_{\rho}\geq K } \sigma_{\rho}$ divided by the index of $K$.

\

For example, if $\Sigma$ is a torus, $\pi_{1}(\Sigma)=\mathbb Z^{2}$. The only effective irreducible representations of finite abelian groups which have endomorphism ring $\mathbb R$ are the trivial representation and the sign representation of $\mathbb Z_{2}$ on $\mathbb R$. Therefore $\sigma_{f_{K_{\rho}},\rho}$ can differ from the identity only in the following four cases:

\begin{itemize}
\item The trivial representation $\rho_{0}$. If $D\dbar_{f_{0}}$ is not isotopic to a complex isomorphism, then $\sigma_{\rho_{0}}=-1$, and otherwise $\sigma_{\rho_{0}}=1$.
\item The representations $\rho_{1}, \rho_{2},\rho_{3}$ which factor through the sign representation via  homomorphisms with kernels 
\[K_{1}=<2a,b>\leq \mathbb Z^{2}\]
\[K_{2}=<a,2b>\leq \mathbb Z^{2}\]
\[K_{3}=<a+b,a-b>\leq\mathbb Z^{2}\]
If $D\dbar_{f_{K_{\rho_{i}}},\rho_{i}}$ is not isotopic to a complex linear map through equivariant isomorphisms, then 
\[\sigma_{\rho_{i}}=1-2K_{i}\]
and otherwise $\sigma_{\rho_{i}}=1$.
\end{itemize}

The contribution of $f_{K}$ to $n_{g,\beta}$ is the coefficient of $K$ in
\[\sigma_{\rho_{0}}\times \prod_{K_{i}\geq K}\sigma_{\rho_{i}}\]
divided by the index of $K$.
The only cases in which this coefficient can be nonzero are as follows:
\begin{itemize}\item
$K=\pi_{1}(\Sigma)$. Then the contribution of $f_{K}=f_{0}$ to $n_{g,\beta}$ is $\sigma_{\rho_{0}}$, as is usual in Gromov-Witten theory.
\item $K=K_{i}$. Then the contribution of the double cover $f_{K}$ to $n_{g,\beta}$ is $0$ if $\sigma_{\rho_{i}}=1$, and $-\sigma_{\rho_{0}}$ otherwise.
\item $K=K_{1}\cap K_{2}=K_{2}\cap K_{3}=K_{3}\cap K_{1}$. Then the contribution of the $4$-fold cover $f_{K}$ to $n_{g,\beta}$ is
\begin{itemize}
\item $0$  if  at least two of $\sigma_{\rho_{1}},\sigma_{\rho_{2}},\sigma_{\rho_{3}}$ are $1$,
\item  and $\sigma_{\rho_{0}}$ if  less than two of $\sigma_{\rho_{1}},\sigma_{\rho_{2}},\sigma_{\rho_{3}}$ are $1$.

\end{itemize}
\end{itemize}

\

\subsection{Comparison to Taubes' Gromov invariant}

\

In \cite{taubessw2}, Taubes defines an integer invariant $Gr$ which, for a generic choice of compatible almost complex structure, counts embedded holomorphic submanifolds of a symplectic $4$-manifold with $b_{2}^{+}\geq 1$. Taubes' Gromov invariant involves a count of collections of embedded holomorphic curves with mutual intersection numbers $0$. (In particular, such collections of embedded holomorphic curves allowed multiple copies of a single embedded  holomorphic torus with trivial normal bundle, but only a single copy of all other types of curves.) In section 5e of \cite{taubessw2}, Taubes shows that a collection of $m$ identical copies of a given torus must count with a weight determined by the  coefficient of $z^{m}$ in a formal power series in the form
\[P(z)^{a}P(z^{2})^{b}P(z^{4})^{c}\]
In our setup, 
\begin{itemize}
\item $a$ is the contribution of the given embedded torus to $n_{(1,\beta)}$,
\item $b$ is the contribution of two-fold covers of the given torus to $n_{(1,\beta)}$,
\item and $c$ is the contribution of $4$-fold covers of the given torus to $n_{(1,\beta)}$.
\end{itemize}

As noted in section 5e of \cite{taubessw2}, any choice of $P(z)$ with first two terms $1+z$ may be used to define an invariant count of holomorphic submanifolds.  In order to get equality of $Gr$ with the Seiberg-Witten invariants in \cite{taubessw3}, Taubes chose $P(z)=(1-z)^{-1}$. As all other components of curves involved in Taubes' Gromov invariant are honestly embedded and counted with the usual signs, the above observation allows us to reconstruct $Gr$ from our integer invariants. 

We shall  state the correspondence between $Gr(A)$ and our invariants in the case that $Gr(A)$ counts rigid embedded holomorphic submanifolds representing homology class $A$. (This rigidity is equivalent to the identity  $A\cdot A+\int_{A}c_{1}=0$.) Consider  the following formal power series in the commuting variables $(g,\beta)$ 
\[\prod_{(g,\beta)}(1-(g,\beta))^{-n_{g,\beta}}\]
Then $Gr(A)$ is the sum of the coefficients  of the terms in the form $(g_{1},\beta_{1})(g_{2},\beta_{2})\dotsb (g_{n},\beta_{n})$ so that
\begin{itemize}
\item $A=\sum_{i}[\beta_{i}]$, where $[\beta_{i}]$ is the homology class represented by $\beta_{i}$,
\item The adjunction inequality predicts that holomorphic curves in the class $(g_{i},\beta_{i})$ will be embedded i.e. $2(g_{i}-1)=[\beta_{i}]\cdot [\beta_{i}]-\int_{\beta_{i}}c_{1}$.
\item $[\beta_{i}]\cdot[\beta_{j}]=0$ for all $i\neq j$.
\end{itemize}
Note that because we do not consider terms with $(g,\beta)^{2}$ unless $[\beta]\cdot[\beta]=0$, we may replace $(1-(g,\beta))^{-n_{g,\beta}}$ with $1+n_{g,\beta}(g,\beta)$ in every case except  the case of tori with trivial normal bundle.

More generally, $Gr$ counts embedded holomorphic submanifolds passing through some generic points and $1$ dimensional submanifolds.  In this situation, the multiply covered tori already considered are still the only non-embedded curves which must be considered, so $Gr$ may also be recovered easily from our integer invariants.

\

The construction of normally complex perturbations of $\dbar$ applies equally well to counting holomorphic curves with disconnected domains. In a symplectic $4$-manifold, one could define an invariant $Gr'(A)$ which uses a normally complex perturbation to count automorphism-free  curves passing through some generic points and one-manifolds  with the extra conditions that each component  satisfies the adjunction formula for embedded curves, the intersection number of different components is $0$, and the total homology class of the curve is $A$. Such an invariant $Gr'$ is identical to Taubes' Gromov invariant, except the weights associated to multiples of rigid tori use $P(z)=1+z$ instead of $P(z)=(1-z)^{-1}$.

  Suppose that a normally complex perturbation of $\dbar$ has already been chosen for curves with connected domains. A Kuranishi chart for a curve $f$ with a disconnected domain may be obtained by simply multiplying together Kuranishi charts corresponding to each connected component. We can similarly pull back our normally complex perturbation of $\dbar$ to this product Kuranishi chart, however the resulting product perturbation will not be normally complex in general, because our product Kuranishi chart will have extra automorphisms acting which correspond to permuting identical components  of $f$. 

In the case of counting rigid curves, the effect of modifying a pulled back perturbation to be normally complex can be calculated, giving a formula for how multiple copies of a given rigid curve should be counted. In particular, $n$ identical copies of a given curve count as $(-1)^{n}$ if the given curve counts $-1$ towards $n_{g,\beta}$, and more than $2$ identical copies of a given curve count as $0$ if the given curve  counts $+1$ towards $n_{g,\beta}$. Equivalently, the number of curves with $n$ components, all in the class $(g,\beta)$, is given by the coefficient of $z^{n}$ in the power series 
\[(1+z)^{n_{g,\beta}}\]
 In general, the above implies that counts of rigid curves with disconnected domains may be encoded in the following generating function
\[\prod_{(g,\beta)}(1+(g,\beta))^{n_{g,\beta}}\]
which should be considered as a formal power series in the commuting  variables $(g,\beta)$. The coefficient of $(g_{1},\beta_{1})\dotsc (g_{n},\beta_{n})$ in the above generating function is the number of automorphism free, rigid curves with $n$ components of genus $g_{1},\dotsc,g_{n}$, and homotopy class $\beta_{1},\dotsc,\beta_{n}$ respectively.  To compute $Gr'(A)$, in the case that  no points or $1$-manifolds are needed to rigidify curves,  we just need to sum the coefficients of the terms $(g_{1},\beta_{1})\dotsc (g_{n},\beta_{n})$ so that $2(g-1)=[\beta_{i}]\cdot [\beta_{i}]-\int_{\beta}c_{1}$, and   $[\beta_{i}]\cdot[\beta_{j}]=0$, and $\sum_{i}[\beta_{i}]=A$.

\

Although counting automorphism-free curves using a normally complex perturbation does not give Taubes' Gromov invariant, we can still consider Taubes' invariant to be a count of curves using a normally complex perturbation as follows: Choose a normally complex perturbation of $\dbar$ on the space of curves with possibly disconnected domain. Then $Gr(A)$ counts (isomorphism classes of) such curves whose automorphisms only permute the different components, and whose components satisfy the usual homological conditions--- namely, each component satisfies the adjunction formula for embedded curves, the intersection number of different components is $0$, and the total homology class of the curve is $A$. (In the case that the expected dimension $A\cdot A+\int_{A}c_{1}>0$, then these curves are counted passing through an appropriate number of generic points and  one-manifolds.)

In the generic situation, the only such curves which will actually have automorphisms will be the curves with several copies of the same rigid torus. If such a rigid torus counts $1$ towards $n_{(1,\beta)}$, then it is easily calculated that $n$ copies of this torus will count $1$ towards our curve count. This has the effect of changing the polynomial $P(z)$ used to calculate weights from $(1+z)$ to $(1-z)^{-1}$. It is an involved combinatorial exercise to calculate directly that after perturbing, all the contributions corresponding more than 1 copy of a negative rigid torus will count $0$ towards $Gr$.

\

From now on, following this paper will require some knowledge of exploded manifolds from \cite{iec}, however the reader unfamiliar with exploded manifolds may still be able to understand sections \ref{integer section} and \ref{future section}. This is a preliminary version of this paper which has not yet been sent to a journal, so I would appreciate any comments or suggestions.

\section{Finite group actions on exploded manifolds}

\ 

We shall now be working with exploded manifolds. As the regularity of embedded Kuranishi charts constructed in \cite{evc} is $\C\infty1$, we shall be working in the category of $\C\infty1$ exploded manifolds. The reader need not have a full understanding of this regularity in order to understand this paper, because we shall show later on that our embedded Kuranish charts may be `rigidified' to have a structure which is even stronger than smooth. 

This section contains local considerations of exploded manifolds $\ex F$ with the action of a finite group $G$. Because the topology on an exploded manifold $\ex F$ is not Hausdorff, we shall need to take note of when $G$ sends a point in $\ex F$ to a topologically indistinguishable point. 

\begin{defn}\label{weakly fixed}Suppose that a group  $G$ acts on a $\C\infty1$ exploded manifold $\ex F$. Say that a point $p\in \ex F$ is weakly fixed by $g\in G$ if the image of $p$ in $\totl{\ex F}$ is fixed by the induced action of $g$ on $\totl{\ex F}$. Let the weak isotropy group $I_{p}\leq G$ of $p$ be the subgroup of  elements that weakly fix $p$.  \end{defn}

When defining how we `rigidify' Kuranishi charts later on, we shall need the following notion of a rigid map between coordinate charts:
 
\begin{defn}\label{rigid} A rigid map between coordinate charts $\mathbb R^{n}\times \et mP\longrightarrow \mathbb R^{n'}\times \et {m'}Q$ is map of the form
\[(x,\tilde z_{1},\dotsc,\tilde z_{m})\mapsto (Ax,e^{a_{1}(x)}\tilde z^{\alpha_{1}},\dotsc ,e^{a_{m'}(x)}\tilde z^{\alpha_{m'}})\]  
where $A$ is affine, the $a_{i}$ are complex valued affine maps, and the $\alpha_{i}$ are  multi-indices in $\mathbb Z^{m}$.

\

A group action on $\mathbb R^{n}\times \et mP$ is rigid if it acts by rigid maps.
\end{defn}

A second way of thinking of a rigid map is as follows: There is a standard connection on any coordinate chart which preserves canonical basis of the tangent space, given by $\frac \partial{\partial x_{i}} $ and the real and imaginary parts of $\tilde z_{j}\frac \partial{\partial \tilde z_{j}}$. A map between coordinate charts is rigid if and only if it preserves the standard connection. 

We shall now work towards proving that we can choose coordinate charts so that a finite group action on an exploded manifold is rigid in the above sense. The following lemma allows us to modify maps to be $G$-equivariant.

\begin{lemma}\label{modify to G map}Let $\ex F$ be a connected $\C\infty1$ exploded manifold with an action of a finite group $G$, and let 
\[f:\ex F\longrightarrow \mathbb R^{n}\times \et mP\]
be any $\C\infty1$ map. Given any rigid action of $G$ on $\mathbb R^{n}\times \et mP$ so that the tropical part $\totb{f}$ of $f$ is $G$-equivariant, and so that there is there is a point  $p \in\ex F$ weakly fixed by $G$ so that $f$ is $G$-equivariant at $p$, there exists an isotopy $f_{t}$ of $f$ to a $G$-equivariant map 
\[f_{1}:\ex F\longrightarrow \mathbb R^{n}\times \et mP\]
so that
\begin{enumerate}\item  $f_{t}(p)$ is independent of $t$, 
\item the set of points $x$ at which $f$ is $G$-equivariant and $f_{t}(x)$ is independent of $t$ is separated from the set of points $x'$ at which $f$ is $G$-equivariant but $f_{t}(x)$ is not independent of $t$.
\item whenever $f_{t}(x)$ is independent of $t$,   $T_{x}f_{t}$ is independent of $t$ if $Tf$ is also  $G$-equivariant at $T_{x}\ex F$.
\item Given any $G$-equivariant, rigid map 
\[\phi:\mathbb R^{n}\times \et mP\longrightarrow \mathbb R^{n'}\times \et {m'}{P'}\] $\phi\circ f_{t}$ is independent of $t$ on the connected component containing $p$ of the set of points $x$ at which  $\phi\circ f$ is $G$-equivariant.
\end{enumerate}\end{lemma}

\pf

The standard basis for $T(\mathbb R^{n}\times \et mP)$ consists of vectors in the form $\frac\partial {\partial x_{i}}$ and the real and imaginary parts of $\tilde z_{i}\frac\partial{\partial \tilde z_{i}}$. An action of $G$ on $\mathbb R^{n}\times \et mP$ is rigid if and only if the connection which preserves this basis is $G$ invariant. Let $\psi_{v}$ denote exponentiation of a vector(field) $v$ using this connection.
The construction of $f_{1}$ is essentially to average $g^{-1}\circ f\circ g$ over all $g\in G$ using this connection. 

For  $g\in G$, denote the various $g$ actions simply by $g$. Because $\totb f\circ g=g\circ \totb f$, for each $g\in  G$,  there exists a unique $\C\infty1$ section $v_{g}$ of $f^{*}T(\mathbb R^{n}\times \et mP)$ which vanishes at $p$ so that at any point $x\in \ex F$
\[\lrb{g\circ \psi_{v_{g}(x)}}( f(x)) =f (g(x))\]
Note that if $f$ is $G$ equivariant at $x$, then  $\psi_{v_{g}(x)}$ is the identity, so $v_{g}(x)$ is some sum of integer multiples of $4\pi$ times the imaginary parts of $ \tilde z_{i}\frac\partial{\partial \tilde z_{i}}$. In particular, this implies that the set of $G$-equivariant points where $v_{g}$ is $0$ is separated from the set of $G$-equivariant points where $v_{g}\neq 0$.  Similar considerations for the composition $\phi\circ f$ imply that $d\phi(v_{g})=0$ at $p$ and on the connected component containing $p$ of the set of points $x$ at which $\phi\circ f$ is equivariant. Note also that  if $v_{g}(x)=0$ and $Tf$ is equivariant at $x$, then the derivative of $v_{g}$ also vanishes at $x$.

Consider the following calculation of the relationship between $v_{gh}$, $v_{g}$ and $v_{h}$.
\[\begin{split}(gh\circ \psi_{v_{gh}(x)})( f(x))&=(f\circ  g)(h(x))
\\& =(g\circ \psi_{v_{g}(h(x))}) f(h(x))
\\ &=g\circ \psi_{v_{g}(h(x))}\lrb{\lrb{ h\circ \psi_{v_{h}(x)}}f(x)}\end{split}\]
so in particular, at a point $x$ in $\ex F$,
\[\psi_{v_{gh}(x)}=h^{-1}\circ \psi_{v_{g}(h(x))}\circ  h\circ \psi_{v_{h}(x)}\]
As $G$ preserves the connection used to define $\psi_{v}$, we may exchange the above equation with 
\begin{equation}\label{psivgh}\psi_{v_{gh}(x)}=\psi_{h^{-1}v_{g}(h(x))}\circ  \psi_{v_{h}(x)}\end{equation}
where $h^{-1}v$ indicates the action of $h^{-1}$ on $T(\mathbb R^{n}\times \et mP)$ composed with  $v$ thought of as a map $v:F\longrightarrow T(\mathbb R^{n}\times \et mP)$. Noting that our connection is flat, and using the standard basis for $T(\mathbb R^{n}\times \et mP)$ to identify vectors at different points, we may reformulate equation (\ref{psivgh}) to read
\[v_{gh}(x)=h^{-1}\lrb{v_{g}(h(x))}+v_{h}(x)\]
Strictly speaking, equation (\ref{psivgh}) only implies that the above equation holds mod $4\pi$ times the imaginary parts of $\tilde z_{i}\frac\partial {\partial \tilde z_{i}}$, however, it holds exactly at $x=p$ and both sides are continuous, therefore it must hold exactly at all points.

Now, define
\[v:=\frac 1{\abs G}\sum_{g\in G}v_{g}\]
We have that at any point $x\in\ex F$
\[v(x)=h^{-1}\lrb{v(h(x))}+v_{h}(x)\]

Now we may define $f_{t}:=\psi_{tv}\circ f$. The map $f_{1}$ is $G$-equivariant, as shown by the following calculation:
\[\begin{split}(f_{1}\circ h)(x)&=\psi_{v(h(x))}\lrb{ f( h(x))}
\\& = \psi_{v(h(x))}\lrb{ h\lrb{ \psi_{v_{h}(x)} \lrb{f (x)}}}
\\&=h\circ (\psi_{h^{-1}v(h(x))}\circ \psi_{v_{h}(x)})\lrb{f(x)} 
\\ &= h\circ \psi_{h^{-1}\lrb{v(h(x))}+v_{h}(x)}(f(x))
\\& =(h\circ f_{1})(x)\end{split}\]

Note that the set of $x$ where $f$ is equivariant and $v(x)=0$ is separated from the set of $x'$ where $f$ is equivariant but $v(x')\neq 0$. At all $x$ in the first set, $f_{t}(x)$ is independent of $t$, and if $Tf$ is also equivariant at such $x$, the derivative of $v$ vanishes, so $T_{x}f_{t}$  is independent of $t$. Similarly $\phi\circ f_{t}$ is independent of $t$ at $p$ and on the connected component containing $p$ of the set of points where $\phi\circ f$ is equivariant.  

\stop

\begin{lemma}\label{G action}Any finite group action on a $\C\infty1$ exploded manifold $\ex F$ is locally rigid in the sense that around every point with weak isotropy group $G$, there exists some $G$-equivariant neighborhood $U$ with a $G$-equivariant isomorphism onto an open subset of $\mathbb R^{n}\times \et mP$ with a rigid $G$ action in sense of  Definition \ref{rigid}.
\end{lemma}

\pf 

 Suppose that a point in  $\ex F$ has weak isotropy group $G$. The set of points topologically equivalent to $p$ in $\ex F$ is isomorphic to $\et m{P^{\circ}}$ where $P^{\circ}$ is an open polytope. The action of each element of  $G$ on $\et m{P^{\circ}}$ is in the form
\[(\tilde z_{1},\dotsc,\tilde z_{m})\mapsto (c_{1}\tilde z^{\alpha_{1}},\dotsc ,c_{m}\tilde z^{\alpha_{m}})\]  
where the $c_{i}$ are constant nonzero complex numbers, and the $\alpha_{i}$ are  multi-indices in $\mathbb Z^{m}$.

The tangent bundle $T\ex F$ restricted to $\et m{P^{\circ}}$ is isomorphic to $\mathbb R^{n}\times  T\et m{P^{\circ}}$, where the $\mathbb R^{n}$-bundle is the sub bundle of  $T\ex F$ which is orthogonal to $T\et m{P^{\circ}}\subset T\ex F$ using some $G$-invariant metric. The action of $G$ on this vector bundle preserves $T\et m{P^{\circ}}\subset \mathbb R^{n}\times T\et m{P^{\circ}}$, and therefore preserves the orthogonal $\mathbb R^{n}$ bundle. Note that our $G$ action sends $\C\infty1$ sections to $\C\infty1$ sections and any $\C\infty1$ section of a vector bundle over $\et m{P^{\circ}}$ is constant. There is therefore an induced linear action of $G$ on $\mathbb R^{n}$. 

A neighborhood $U$ of $\et m{P^{\circ}}$ in $\ex F$ is isomorphic to $\mathbb R^{n}\times\et mP$ where $P$ is a polytope with $P^{\circ}$ as its interior strata. The isomorphism is via a map
\[\Phi:U\longrightarrow \mathbb R^{n}\times \et mP\longrightarrow\ex F\]
which sends $\et m{P^{\circ}}\subset U\subset \ex F$ identically to  $0\times \et m{P^{\circ}}\subset \mathbb R^{n}\times \et mP$, and with derivative at $\et m{P^{\circ}}$ which preserves the splitting of the tangent space into the tangent space to $\mathbb R^{n}$ and the tangent space to $\et m{P^{\circ}}$. Our linear action of $G$ on $\mathbb R^{n}$ and the action of $G$ on $\et m{P^{\circ}}$ extend uniquely to a rigid action of $G$ on $\mathbb R^{n}\times \et mP$.
 Our map $\Phi$ is not necessarily $G$-equivariant,  but its tropical part $\totb\Phi$ is equivariant, and $T\Phi$ is $G$-equivariant restricted to $\et m{P^{\circ}}$, so Lemma \ref{modify to G map} allows us to modify $\Phi$ to a $G$-equivariant  map $\Phi_{1}$ equal to $\Phi$ at $\et m{P^{\circ}}$ to first order. It follows that $\Phi_{1}$ restricted to some neighborhood of $\et m{P^{\circ}}$ is a $G$-equivariant isomorphism  onto an open subset of $\mathbb R^{n}\times\et mP$.

\stop

\begin{remark}The proof of the above lemma not only shows that any finite group action on $\ex F$ is locally rigid--- it also shows that coordinates may be chosen so that this action is a split action on $\mathbb R^{n}\times \et mP$ consisting of maps in the form 
\[(x,\tilde z)\mapsto (Ax,c_{1}\tilde z^{\alpha_{1}},\dotsc, c_{m}\tilde z^{\alpha_{m}})\]
where $A$ is a linear map and the $c_{i}$ are constants in $\mathbb C^{*}$. \end{remark}

A consequence of Lemma \ref{G action} is that the $G$-fixed set of $\ex F$ is always a closed $\C\infty1$ sub exploded  manifold of $\ex F$. The weakly  $G$-fixed set of $\ex F$ is also a locally finite union of closed $\C\infty1$ sub exploded manifolds of $\ex F$.

\begin{lemma}\label{G map model} Every $G$-equivariant submersion is locally isomorphic to an equivariant rigid map. 

In particular, given a $G$-equivariant submersion \[\phi:\ex F\longrightarrow \ex X\  ,\] for any point $p$ in $\ex F$, there exist $I_{p}$-invariant open neighborhoods $U\subset \ex F$ of $p$ and  $V\subset \ex X$ of $\phi(U)$, and a commutative $I_{p}$-equivariant diagram
\[\begin{tikzcd}\ex F\supset U\dar{\phi} \rar{\psi}&\mathbb R^{n}\times \et mP\dar{\Phi}
\\ \ex X\supset V\rar{\psi_{0}}&\mathbb R^{n'}\times \et {m'}{P'}\end{tikzcd}\]
where $\psi$ and $\psi_{0}$ are $I_{p}$-equivariant isomorphisms on their image, and $\Phi$ is a rigid  $I_{p}$-equivariant map, and the action of $I_{p}$ on the domain and target of $\Phi$ is rigid.
\end{lemma}

\pf The weak stabilizer of $\phi(p)$ must contain $I_{p}$, so Lemma \ref{G action} implies that we can choose an $I_{p}$-invariant open neighborhood $V$ of $\phi(p)$ with an $I_{p}$-equivariant isomorphism $\psi_{0}$ onto an open subset of $\mathbb R^{n'}\times \et {m'}{P'}$ with a rigid action of $I_{p}$.

Similarly, there is an open neighborhood $U'\subset \phi^{-1}(V)$ of $p$ with an $I_{p}$-equivariant isomorphism 
\[\psi':U'\longrightarrow \mathbb R^{n}\times \et mP\]
onto an open subset of $\mathbb R^{n}\times \et mP$ with a rigid action of $I_{p}$.
There is then a unique rigid map 
\[\Phi:\mathbb R^{n}\times \et mP\longrightarrow \mathbb R^{n'}\times \et {m'}{P'} \]
with the property that the diagram 
\[\begin{tikzcd} U'\dar{\phi} \rar{\psi'}&\mathbb R^{n}\times \et mP\dar{\Phi}
\\  V\rar{\psi_{0}}&\mathbb R^{n'}\times \et {m'}{P'}\end{tikzcd}\]
commutes when restricted to points $x\in U'$ topologically indistinguishable from $p$, and the derivative of the maps in the above diagram also commute when restricted to the same points. As all other maps in the diagram are $I_{p}$-equivariant, it follows that for any $g$ in $I_{p}$,  $g\Phi g^{-1}$ satisfies the same commutativity properties. As the action of $I_{p}$ on the domain and target of $\Phi$ is rigid, it follows that $g\Phi g^{-1}$ is also rigid, therefore $g\Phi g^{-1}$ must be equal to $\Phi$, and $\Phi$ is a $I_{p}$-equivariant rigid map.

As the above diagram commutes up to first order at any point topologically indistinguishable from $p$, and $\Phi$ is a submersion, $\psi'$ may be modified on some neighborhood $U''$ of $p$ to a map $\psi''$ which is an isomorphism onto an open subset of $\mathbb R^{n}\times \et mP$ so that the diagram  
\[\begin{tikzcd} U''\dar{\phi} \rar{\psi''}&\mathbb R^{n}\times \et mP\dar{\Phi}
\\  V\rar{\psi_{0}}&\mathbb R^{n'}\times \et {m'}{P'}\end{tikzcd}\]
commutes, and so that $\psi''$ is equal to $\psi'$ to first order at any point topologically indistinguishable to $p$. We have that $\psi''$ is $I_{p}$-equivariant to first order at any point topologically indistinguishable  from $p$, and $\Phi\circ \psi''$ is $I_{p}$-equivariant therefore lemma \ref{modify to G map} allows us to modify $\psi''$ to a  $I_{p}$-equivariant map $\psi$ which is equal to $\psi''$ to first order at any point topologically indistinguishable to $p$ and so that $\Phi\circ \psi=\Phi\circ \psi''$. As $\psi''$ was an isomorphism onto an open set, $\psi$ must be an isomorphism onto an open set when restricted to a sufficiently small neighborhood $U$ of $p$.

 We have now constructed the  required local isomorphism of $\phi$ with an equivariant rigid map $\Phi$.
 \[\begin{tikzcd} U\dar{\phi} \rar{\psi}&\mathbb R^{n}\times \et mP\dar{\Phi}
\\  V\rar{\psi_{0}}&\mathbb R^{n'}\times \et {m'}{P'}\end{tikzcd}\]

\stop

\begin{lemma}\label{G modify}Suppose that $\ex F$ and $\ex Y$ are $\C\infty1$ exploded manifolds with an action of a finite group $G$. Then, given any $\C\infty1$ map 
\[f:\ex F\longrightarrow \ex Y\]
so that $\totb{f}$ is $G$-equivariant, 
there exists an isotopy $f_{t}$ of $f$ to a map which is $G$ equivariant in an open neighborhood of the set where $f$ is $G$-equivariant. Moreover, we can choose $f_{t}$ so that if $x$ is any point at which $f$ is $G$-equivariant, $f_{t}(x)$ is independent of $t$, and if $Tf$ is also $G$-equivariant at $x$, $T_{x}f_{t}$ is also independent of $t$.

\

If moreover, there is a $G$-equivariant submersion $\phi:\ex Y\longrightarrow \ex X$, then $f_{t}$ may be chosen so that $\phi \circ f_{t}$ is independent of $t$ at any point $x$ at which $\phi \circ f$ is equivariant in a neighborhood of the set of points where $f$ is equivariant.  
\end{lemma}

\pf 

We prove the more general case in which there is a $G$-equivariant submersion $\phi:\ex Y\longrightarrow \ex X$.
Lemma \ref{G map model} implies that any point $y$ in $\ex Y$ is contained in some connected $I_{y}$-invariant neighborhood $U$ on which the action of $I_{y}$ is rigid and $\phi$ is isomorphic to a rigid $I_{y}$-equivariant  map. We may also choose $U$ so that the every $G$-translate of $U$ is either disjoint from $U$ or equal to $U$. Let $U'$ be the $G$-orbit of any compactly contained, open subset of $U$. Lemma \ref{modify to G map} implies that there is an isotopy of $f$ to a map which is equivariant on a neighborhood of the points $x$ in $f^{-1}(U')$ where $f$ is equivariant, and that we may choose this isotopy 
\begin{itemize}\item  to be constant at points $x$ where $f$ is already $G$-equivariant, 
\item  to have constant derivative where $Tf$ is $G$-equivariant, 
\item and to have constant composition with $\phi$ for any point $x$ at which $\phi\circ f$ is $G$-equivariant in a neighborhood of the points where $f$ is $G$-equivariant. 
\end{itemize}

Choose a countable cover of $\ex Y$ by such open sets $U'_{i}$ where the above isotopies may be constructed. Suppose that an isotopy with the given properties has been constructed on the union of  $f^{-1}(U'_{i})$ for $i<k$. We may follow this with an isotopy with the required properties on a neighborhood of $f^{-1}(U'_{k})$ which is constant on $f^{-1}(U'_{i})$ for all $i<k$. 

Doing the same for all $i$ gives a countable sequence of isotopies which are eventually  constant on a neighborhood of any point. We may therefore combine them into a single isotopy  with the required properties.

\stop

\section{Normally rigid structure}

In this section, we shall construct a kind of rigid structure normal to certain strata of $\ex F$ with the action of a finite group $G$. If $\ex F$ is a smooth manifold, these strata are the sub manifolds of $\ex F$ consisting of points with isotropy subgroup $I$ so that the action of $I$ on the tangent space is isomorphic to a particular representation.  In this case, the rigid structure is a $G$-equivariant map of the normal bundle of each strata to $\ex F$ so that the maps from different strata are compatible in a way specified by Definition \ref{normally rigid structure} below.

\begin{defn}[Strata of $(\ex F,G)$] \label{strata def} Use the notation $\s$ to denote 
\begin{itemize}\item A $m$-dimensional polytope, $P_{\s}$,
 \item a finite group $I_{\s}$,
\item an action of $I_{\s}$ on $\et m{P_{\s}}$ by rigid transformations,
\item and a linear representation of $I_{\s}$ on $\mathbb R^{n}$ which fixes only $0$.
\end{itemize}
Let $(\ex F,G)^{\s}\subset \ex F$ denote the subset of $\ex F$ consisting of points $p$ with weak isotropy group $I_{p}$ isomorphic to $I_{\s}$ and tropical structure $\mathcal P(p)$ isomorphic to $P_{\s}$ so that the action of $I_{\s}$ in a neighborhood of $p$ is isomorphic to the action of $I_{\s}$ on an open subset of $\mathbb R^{n}\times \et m{P_{\s}}$ times some $\mathbb R^{k}$ with a trivial group action. 

Call $(\ex F,G)^{\s}$ a strata of $(\ex F,G)$. When there is no ambiguity about the group action on $\ex F$, we shall use $\ex F^{\s}$ instead of $(\ex F,G)^{\s}$. Given a point $f\in \ex F$, use the notation $\s(f)$ to indicate the strata containing $f$, so $f\in\ex F^{\s(f)}$.

Say that $\s\geq \s'$ if $(\mathbb R^{n}\times \et m{P_{\s}}, I_{\s})^{\s'}$ is nonempty. Say that $\s=\s'$ if $I_{\s}$ is isomorphic to $I_{\s'}$, $P_{\s}$ is isomorphic to $P_{\s'}$, and the action of $I_{\s}$ and $I_{\s'}$ on $\mathbb R^{n}\times \et m{P_{\s}}$ is isomorphic.
 \end{defn}

\begin{remark}	
\begin{itemize}
\item Lemma \ref{G action} implies that $\ex F^{\s}$ is an exploded submanifold of $\ex F$, and each point in $\ex F$ is in $\ex F^{\s}$ for some $\s$.
\item The closure of $\ex F^{\s'}$ contains $\ex F^{\s}$ if and only if $\s\geq \s'$.
\item The action of $G$ on $\ex F$ preserves $\ex F^{\s}$, so $\ex F^{\s}$ has a given $G$ action.
\item $\geq$ defines a partial order on the set of possible $\s$. 
\item If $\s\geq \s'$ then $I_{\s'}$ is isomorphic to a subgroup of $I_{\s}$ and $P_{\s'}$ is isomorphic to a face of $P_{\s}$. 
\item If  $\s>\s'$, then either $I_{\s'}$ has fewer elements than $I_{\s}$ and $P_{\s'}$ has dimension at most the dimension of $P_{\s}$ or $\abs{I_{\s'}}\leq \abs{I_{\s}}$ and $\dim P_{\s'}<\dim P_{\s}$.
\item So long as $G$ is finite and $\ex F$ is finite dimensional, there is a bound on the length of any chain $\dotsb <\s_{i}<\s_{i+1}<\dotsb$ so that each $\ex F^{\s_{i}}$ is nonempty.

\end{itemize}
\end{remark}

The usual definition of a normal bundle for $\ex F^{\s}$ as the quotient of the pullback of $T\ex F$ by $T \ex F^{\s}$ is not adequate for our purposes, because it will not be isomorphic to any  neighborhood of $\ex F^{\s}$ in $\ex F$ if $P_{\s}$ has non empty interior $P^{\circ}_{\s}$. Instead we shall define a notion of a normal bundle to $\ex F^{\s}$ which is a rigid  $\mathbb R^{n}\times \et m{P_{\s}}$ bundle in the sense defined below. 

\begin{defn}\label{etmp bundle} A rigid $\mathbb R^{n}\times\et mP$ bundle over a manifold $M$ is an exploded manifold $\ex E$ with a surjective submersion $\pi_{M}:\ex E\longrightarrow M$ with fibers isomorphic to $\mathbb R^{n}\times\et mP$ 
\[\begin{tikzcd}\mathbb R^{n}\times\et mP\rar& \ex E\dar{\pi_{M}}
\\ & M\end{tikzcd}\]
and a sheaf, $\Charts$, of sets on $M$ so that each element of $\Charts (U)$ is a map
 \[\psi:\pi^{-1}_{M}(U)\longrightarrow  \mathbb R^{n}\times \et mP\]
so that 
\begin{itemize}\item Each $\psi$ in $\Charts(U)$ is an isomorphism when restricted to $\pi_{M}^{-1}(p)$ for any point $p$ in $U$.
\item  Each point in $M$ has a neighborhood so that $\Charts(U)$ is nonempty.
\item If $\psi$ is in $\Charts(U)$, then a map $\psi'$ is in $\Charts(U)$ if and only if there exists a commutative diagram of smooth maps
\[\begin{tikzcd}\pi^{-1}(U)\dar[swap]{(\pi_{M},\psi)}\ar{dr}{\psi'}
\\ U\times  \mathbb R^{n}\times\et mP \rar{f}&  \mathbb R^{n}\times \et mP\end{tikzcd}\]
so that $f$ is a rigid map when restricted to each $ \mathbb R^{n}\times\et mP$ slice of $U\times\mathbb R^{n}\times  \et mP$.
\end{itemize}

\end{defn}

\begin{remark}Any  map from $ \et m{P^{\circ}}$ to itself extends uniquely to a rigid  map of $ \et mP$ to itself if and only if the tropical part of the map $P^{\circ}\longrightarrow P^{\circ}$ extends to a map $P\longrightarrow P$. The same holds for smooth families of rigid maps. It follows that any  $ \et m{P^{\circ}}$ bundle extends uniquely to a rigid $ \et mP$ bundle if and only if the corresponding $P^{\circ}$ bundle extends to a $P$ bundle. 

Similarly, any rigid $\mathbb R^{n}\times \et m{P^{\circ}}$ bundle extends uniquely to a rigid $\mathbb R^{n}\times \et mP$ bundle if the corresponding $P^{\circ}$ bundle extends to a $P$ bundle.

For example, $\ex F^{\s}$ is a $\et m{P^{\circ}_{\s}}$ bundle over its smooth part $\totl{\ex F^{\s}}$. As any $\et m{P^{\circ}_{\s}}$ bundle is automatically rigid, this bundle extends uniquely to a rigid $\et m{P_{\s}}$ bundle. 

\end{remark}

The tangent space $T\ex F$ of $\ex F$ is the pullback of a vector bundle on $\totl{\ex F}$, so $T_{p}\ex F$ has a canonical identification with $T_{p'}\ex F$ whenever $\totl p=\totl{p'}$. All points in $\totl{\ex F^{\s}}$ are  fixed by $I_{\s}$, so $I_{\s}$ has an action on $T_{p}\ex F$ for all $p$ in $\ex F^{\s}$ which is given by the derivative of the action on $\ex F$ by $h\in I_{\s}$ followed by the identification of $T_{h*p}\ex F$ with $T_{p}\ex F$. 

The usual definition of a normal bundle to $\ex F^{\s}$ is the quotient by $T\ex F^{\s}$ of $T\ex F$ restricted to $\ex F^{\s}$. As $I_{\s}$-invariant vectors are got rid of in this quotient, we can put slightly more structure on our normal bundle by starting with the  subbundle $S\subset T\ex F\rvert_{\ex F^{\s}}$ of $I_{\s}$-anti-invariant vectors. The pre normal bundle to $\ex F^{\s}$ is defined by quotienting $S$ by an equivalence relation which may be thought of as the kernel of the exponential map $S\longrightarrow \ex F$ using any connection. The result of defining the pre normal bundle this way is that it contains a canonical copy of $\ex F^{\s}$ and  has an inclusion into $\ex F$ which is canonical up to first order at $\ex F^{\s}$. 

\begin{defn}[pre normal bundle] Define the pre normal bundle of $\ex F^{\s}$ as follows:
Consider the sub vector-bundle $S$ of $T\ex F$ restricted to $\ex F^{\s}$ consisting of the $I_{\s}$ anti-invariant vectors; so
\[S:=\{(p,v)\subset T\ex F\}\]
where $p\in \ex F^{\s}$ and $v\in T_{p}\ex F$ has the property that
\[\sum_{g\in I_{\s}}g*v=0\]
Let $(p,v_{0})$ be any vector in $T\ex F^{\s}$ which is also in $S$. Then, because $v_{0}$ is $I_{\s}$ anti-invariant, it must be tangent to a $\et m{P^{\circ}}$ fiber of $\ex F^{\s}$, and exponentiation $\psi_{v_{0}}(p)$ of $v_{0}$ is independent of what smooth connection is used. As $\psi_{v_{0}}(p)$ and $p$ will be in the same $\et m{P^{\circ}}$-fiber of $\ex F^{\s}$, $T_{p}\ex F$ and $T_{\psi_{v_{0}}(p)}\ex F$ are canonically isomorphic. 

Define the pre normal bundle $PN_{\ex F^{\s}}$ of $\ex F^{\s}$ to be the quotient of $S$ by the equivalence relation
\begin{equation}\label{er} (p,v+v_{0})=(\psi_{v_{0}}(p),v)\end{equation}
whenever $v_{0}\in (T_{p}\ex F)\cap S$.

Let $n$ be the codimension of $\ex F^{\s}\subset \ex F$. The pre normal bundle $PN_{\ex F^{\s}}$ of $\ex F^{\s}$ has a natural structure of a rigid $\mathbb R^{n}\times \et m{P_{\s}^{\circ}}$ bundle, described as follows:
The above equivalence relation (\ref{er}) preserves fibers of the projection 
\[\begin{array}{ccc}S&\longrightarrow &\totl{\ex F^{\s}}
\\ (p,v)&\mapsto&\totl p\end{array}\]
 so there is a corresponding projection of the pre normal bundle to $\totl{\ex F^{\s}}$. A chart on the pre-normal bundle is equivalent to a (smooth exploded) locally defined map   $S\longrightarrow \mathbb R^{n}\times \et m{P^{\circ}_{\s}}$ which  is rigid and surjective restricted to each fiber of $S\longrightarrow \totl{\ex F^{\s}}$, and so that points in a fiber of $S$ are sent to the same point in $\mathbb R^{n}\times \et m{P^{\circ}_{\s}}$ if and only if they are related by the above equivalence relation (\ref{er}).

\end{defn}

The action of $G$ on $T\ex F$ preserves $S$, and the corresponding action of $G$ on $S$ preserves the equivalence relation (\ref{er}), so there is an action of $G$ on the pre normal bundle of $\ex F^{\s}$. This action of $G$ preserves the rigid structure of this pre normal bundle.

There is a canonical $G$-equivariant inclusion \[\ex F^{\s}\longrightarrow PN_{\ex F^{\s}}\] given by the inclusion $\ex F^{\s}\subset S$ defined by $p\mapsto (p,0)$ followed by the projection $S\longrightarrow PN_{\ex F^{\s}}$.

 There is also a canonical identification of $T(PN_{\ex F^{\s}})\rvert_{\ex F^{\s}}$ with $T\ex F\rvert_{\ex F^{\s}}$ which should be thought of as follows: Given any choice of connection on $\ex F$,  the exponential map defines a map from $S\subset T\ex F$ to $\ex F$. 
 \[(p,v)\mapsto \psi_{v}(p)\]
 The derivative of this exponential map at the zero section $0\times \ex F^{\s}$ is independent of the choice of connection. In particular, at $(p,0)$ it is given by 
 \[(w,v)\mapsto w+v\]
 where $w\in T_{p}\ex F^{\s}$ and $v$ is in $S$ restricted to $p$. The kernel of this map is obviously $(v_{0},-v_{0})$ where $v_{0}\in S\cap T_{p}\ex F^{\s}$. The tangent space to the pre normal bundle at $p$ is equal to the quotient of $T_{(p,0)}S$ by this kernel. Therefore, the above map gives a well defined injective map from the the tangent space of the pre-normal bundle at $\ex F^{\s}$ to $T\ex F$ at $\ex F^{\s}$.  This map may be inverted by the map which sends $v\in T_{p}\ex F$ to $(v^{I_{\s}},v-v^{I_{\s}})\in T_{(p,0)}(S)$, where $v^{I_{\s}}$ indicates the $I_{\s}$-invariant part of $v$.
 
 Said another way: given any connection, the exponential map $S\longrightarrow \ex F$ factorizes through the projection $S\longrightarrow PN_{\ex F^{\s}}$, and the corresponding exponential map $PN_{\ex F^{\s}}\longrightarrow \ex F$ has derivative at $\ex F^{\s}\subset PN_{\ex F^{\s}}$ which is this canonical isomorphism $T(PN_{\ex F^{\s}})\rvert_{\ex F^{\s}}\longrightarrow T\ex F\rvert_{\ex F^{\s}}$.

\begin{defn}[normal bundle] Define the normal bundle $\ex N_{\ex F^{\s}}$ of $\ex F^{\s}\subset \ex F$ to be the unique extension of the pre normal bundle of $\ex F^{\s}\subset \ex F$ to a rigid $\mathbb R^{n}\times \et m {P_{\s}}$ bundle. 
  \end{defn}
  
Note that there is a canonical inclusion \[\iota:\ex F^{\s}\hookrightarrow NP_{\ex F^{\s}}\hookrightarrow \ex N_{\ex F^{\s}}\] and a canonical isomorphism 
\begin{equation}\label{ci}\iota^{*}T \ex N_{\ex F^{\s}}=T\ex F\rvert_{\ex F^{\s} }\end{equation}

  There is also a (non-canonical) map of $\ex N_{\ex F^{\s}}$ into $\ex F$ so that the following diagram commutes
 \[\begin{tikzcd}\ex F^{\s}\ar[bend left,hook]{rr}\rar[hook]{\iota}\ar{dr} &\ex N_{\ex F^{\s}}\dar\rar& \ex F
 \\ & \totl{\ex F^{\s}}\end{tikzcd}\]

   and the induced map
\[\iota^{*}T\ex N_{\ex F^{\s}}\longrightarrow T\ex F\rvert_{\ex F^{\s}}\]
is the canonical isomorphism (\ref{ci}).

\

\

The following definition of a normally rigid structure should be thought of as a consistent identification of a neighborhood of each strata $\ex F^{\s}$ of $\ex F$  with its normal bundle $\ex N_{\ex F^{\s}}$ (or at least a neighborhood of $\ex F^{\s}\subset \ex N_{\ex F^{\s}}$).

\begin{defn}\label{normally rigid structure} A normally rigid structure on $\ex F$ is for each $\s$,
a map
\[r_{\ex F^{\s}}:\ex N_{\ex F^{\s}}\longrightarrow \ex F\] which is  $G$-equivariant map  
so that 
\begin{enumerate}
\item\label{xc1} The inclusion $\ex F^{\s}\hookrightarrow \ex F$ factors as
\[\begin{tikzcd}\ex F^{\s}\ar[bend left,hook]{rr}\rar[hook]{\iota} &\ex N_{\ex F^{\s}}\rar& \ex F\end{tikzcd}\]
\item \label{xc2}The induced map of vector bundles over $\ex F^{\s}$
\[\iota^{*}T\ex N_{\ex F^{\s}}\longrightarrow T\ex F\rvert_{\ex F^{\s}}\]
is the canonical isomorphism (\ref{ci}).
\item \label{xc3}

 For all $\s$, there exists an open neighborhood $U^{\s}$ of $\ex F^{\s}\subset \ex N_{\ex F^{\s}}$ on which $r_{\ex F^{\s}}$ is injective, so that $r_{\ex F^{\s}}(U^{\s})$ intersects $\ex F^{\s'}$ if and only if $\s'\leq \s$,  and so that  for each $\s'<\s$, there is a commutative diagram 
\[\begin{tikzcd} \ex N_{\ex F^{\s'}}\ar{dr}{r_{\ex F^{\s'}}}
 \\ U^{\s'}\cap \Phi^{-1}_{\s'<\s}(U^{\s})\uar[hook] \dar[swap]{\Phi_{\s'<\s}} &\ex F
\\ \ex N_{\ex F^{\s}}\ar{ur}[swap]{r_{\ex F^{\s}}}\end{tikzcd}\]
where $\Phi_{\s'<\s}$ is the unique fiberwise rigid map
\[\ex N_{\ex F^{\s'}}\rvert_{r_{\ex F^{\s}}(U^{\s})}\xrightarrow{\Phi_{\s'<\s}}  \ex N_{\ex F^{\s}}\]
with the property that the differential of $ r_{\ex F^{\s}}\circ \Phi_{\s'<\s}$ restricted to $\ex F^{\s'}$ is the canonical isomorphism (\ref{ci}). The notation $\ex N_{\ex F^{\s'}}\rvert_{r_{\ex F^{\s}}(U^{\s})}$ means the restriction of the bundle $\ex N_{\ex F^{s'}}\longrightarrow \totl{\ex F^{\s'}}$ to the intersection of $\totl{\ex F^{\s'}}$ with the image of $U^{\s}$ under $\totl{r_{\ex F^{s}}}$.
\end{enumerate}

\end{defn}

To reduce notation complexity, the subscript shall be dropped from $r_{\ex F^{\s}}$ where no ambiguity is possible.

Note that after $r_{\ex F^{\s}}$ and $U^{\s}$ are chosen, $\Phi_{\s'<\s}$ is uniquely specified because rigid maps are specified by their derivatives. After $\Phi_{\s'<\s}$ is specified, $r_{\ex F^{\s'}}$ restricted to $U^{\s'}\cap \Phi^{-1}_{\s'<\s}(U^{\s})$ is determined. Normally rigid structures shall be constructed in Lemma \ref{nr construction} below. Roughly speaking, within a neighborhood of $\ex F^{\s'}$, $r_{\ex F^{\s'}}$ is determined by $r_{\ex F^{\s}}$ on $r(U^{\s})$ for all $\s>\s'$.  A choice of map satisfying all conditions on $r_{\ex F^{\s'}}$ apart from $G$-equivariance may easily be constructed, then modified to a $G$-equivariant map using Lemma \ref{G modify}. In this way, normally rigid structures on $\ex F$ may be constructed starting with $r_{\ex F^{\s}}$ for maximal $\s$ and then constructing  the maps $r_{\ex F^{\s'}}$ with smaller $\s$.  Note that once $r_{\ex F^{\s_{i}}}$  have been chosen compatibly for all $\s_{i}>\s$, if $\s<\s_{1}<\s_{2}$, then the diagram
\[\begin{tikzcd}\ex N_{\ex F^{\s}}\rvert_{\bigcap_{i}r(U^{\s_{i}})}\rar{\Phi_{\s<\s_{1}}}\ar{dr}[swap]{\Phi_{\s<\s_{2}}}& \ex N_{\ex F^{\s_{1}}}\rvert_{r(U^{\s_{2}})}\dar{\Phi_{\s_{1}<\s_{2}}}
&\lar[hook]  U^{\s_{1}}\cap \Phi^{-1}_{\s_{1}<\s_{2}}(U^{\s_{2}}) \dar{r_{\ex F^{\s_{1}}}}
\\ & \ex N_{\ex F^{\s_{2}}}\rar{r_{\ex F^{\s_{2}}}}& \ex F\end{tikzcd}\]
commutes, so the conditions imposed on $r_{\ex F^{\s}}$ by each $r_{\ex F^{\s_{i}}}$ are compatible.

\begin{defn}[Normally rigid map]\label{normally rigid map} Suppose that both $(\ex F,G)$ and $(\ex X,H)$ have  a normally rigid structure, and that there is a given surjective homomorphism $G\longrightarrow H$. A $G$-equivariant map \[\psi:\ex F\longrightarrow \ex X\] is normally rigid if for each $\ex F^{\s'}$ and $\ex X^{\s}$,  there exists an  open neighborhood $U\subset \ex N_{\ex F^{\s'}}$  of $\ex F^{\s'}\cap \psi^{-1}\ex X^{\s}$ and a fiberwise rigid map
 \[ U\longrightarrow  \ex N_{\ex X^{\s}}\]
so that the following diagram commutes
\[\begin{tikzcd} \ex N_{\ex F^{\s'}}\dar{r_{\ex F^{\s'}}}&U\rar \lar[hook] & \ex N_{\ex X^{\s}}\dar{r_{\ex X^{\s}}}
\\ \ex F\ar{rr}{\psi}& &\ex X\end{tikzcd}\]
and so that $U$ maps into a neighborhood of $\ex X^{\s}\subset \ex N_{\ex X^{\s}}$ on which $r_{\ex X^{\s}}$ is injective.
\end{defn}

For example, $G$ acts on an exploded manifold with a normally rigid structure by normally rigid maps. Any $G$-equivariant open subset $O$ of an exploded manifold $\ex F$ with a normally rigid structure has a normally rigid structure so that the inclusion $O\longrightarrow \ex F$ is a normally rigid map.

\begin{remark}\label{smooth remark} Even if $\ex F$ only has a $\C\infty 1$ structure, $\ex N_{\ex F^{\s}}$ has a smooth structure, and a choice of normally rigid structure on $\ex F$ also defines a smooth structure on $\ex F$ so that the maps $r_{\ex F^{\s}}$ are smooth.  The condition of being a normally rigid map is stronger than the condition of being a smooth map.
\end{remark}

\begin{lemma}The composition of two normally rigid maps is a normally rigid map.\end{lemma}

\pf

Suppose that $\psi_{1}:\ex F\longrightarrow \ex X$ and $\psi_{2}:\ex X\longrightarrow \ex Y$ are normally rigid maps. The definition of normally rigid maps gives several commutative diagrams
\[\begin{tikzcd}
\\\ex N_{\ex F^{\s'}}\dar{r_{\ex F^{\s'}}}&\lar[hook] U'_{i}\rar& \ex N_{\ex X^{\s_{i}}}\dar{r_{\ex X^{\s_i}}}&\lar[hook]U''_{i}\rar& \ex N_{\ex Y^{\s}}\dar{r_{\ex Y^{\s}}}
\\ \ex F\ar{rr}{\psi_{1}} &&\ex X\ar{rr}{\psi_{2}}&&\ex Y\end{tikzcd}\]
where $U_{i}'\subset \ex N_{\ex F^{\s'}}$ is an open neighborhood of $\ex F^{\s'}\cap \psi_{1}^{-1}(\ex X^{\s_i})$, and $U''_{i}\subset \ex N_{\ex X^{\s_i}}$ is an open neighborhood of $\ex X^{\s_i}\cap \psi_{2}^{-1}(\ex Y^{\s})$. For a fixed strata $\ex Y^{\s}$, we may choose each $U''_{i}$ small enough that it maps into a fixed open neighborhood $O$ of $\ex Y^{\s}\subset \ex N_{\ex Y^{\s}}$ on which $r_{\ex Y^{\s}}$ is injective. We may also choose $U_{i}''$ small enough that $r_{\ex X^{\s_i}}$ is injective restricted to  $U''_{i}$, and $r_{\ex X^{\s_i}}(U''_{i})$ is an open subset of $\ex X$. Then
\[U_{i}:=U_{i}'\cap (\psi_{1}\circ r_{\ex F^{\s'}})^{-1}(r_{\ex X^{\s_i}}(U''_{i}))\]
is an open subset of $\ex N_{\ex F^{\s'}}$ which contains $\ex F^{\s'}\cap\psi_{1}^{-1}(\ex X^{\s_i})\cap (\psi_{2}\circ \psi_{1})^{-1}(\ex Y^{\s})$.  

There is a commutative diagram
\[\begin{tikzcd}\ex N_{\ex F^{\s'}}\dar{r_{\ex F^{\s'}}}&\lar[hook] U_{i}\rar&O\subset\ex N_{\ex Y^{\s}}\dar{r_{\ex Y^{\s}}}
\\ \ex F\ar{rr}{\psi_{2}\circ \psi_{1}}&& \ex Y\end{tikzcd}\]
where the top right arrow is a fiberwise rigid map, uniquely determined by the requirement that the above diagram commutes because $r_{\ex Y^{\s}}$ is injective restricted to $O$.

Therefore, each these maps for different $i$ agree on their overlap, so there is a commutative diagram 
\[\begin{tikzcd}\ex N_{\ex F^{\s'}}\dar{r_{\ex F^{\s'}}}&\lar[hook] U \rar&O\subset\ex N_{\ex Y^{\s}}\dar{r_{\ex Y^{\s}}}
\\ \ex F\ar{rr}{\psi_{2}\circ \psi_{1}}&& \ex Y\end{tikzcd}\]
where the top right arrow is a fiberwise rigid map, and 
\[U:=\bigcup_{i}U_{i}\]
is an open neighborhood in $\ex N_{\ex F^{\s'}}$ of $\ex F^{\s'}\cap (\psi_{2}\circ \psi_{1})^{-1}(\ex Y^{\s})$. 

\stop

\

Note that if we define a category with morphisms given by normally rigid maps, an isomorphism class only remembers a germ of the map $r_{\ex F^{\s}}$ around $\ex F^{\s}\subset \ex N_{\ex F^{\s}} $.

\

The following lemma allows us to construct normally rigid structures. The statement is complicated by our future need to be able to extend normally rigid structures which are already defined on an exploded  submanifold and an open subset.

\begin{lemma}\label{nr construction} Suppose that 
\begin{itemize}\item  $\ex F$ is an exploded manifold with an action of a finite group  $G$.
\item A $G$-invariant  open subset $\ex A\subset\ex F$ has a normally rigid structure \item There is a $G$-invariant sub exploded manifold $\ex B\subset \ex F$ with a normally rigid structure so that $\ex B$ is equal to the transverse intersection of a $G$-equivariant section of some vector bundle  which is defined  over an open subset $O$ of $\ex F$.
\item The inclusion \[\ex B\cap \ex A\longrightarrow \ex A\] is a normally rigid map.
 \item There is an exploded manifold $\ex X$ with  a normally rigid structure  and a $G$-equivariant submersion
 \[\psi:\ex F\longrightarrow \ex X\]
 which is a normally rigid submersion restricted to $\ex A$ and $\ex B$.
  \end{itemize}

Then, given a choice of subset $C_{0}\subset \ex A$, $C_{1}\subset \ex B$  so that the inclusion $C_{i}\longrightarrow \ex F$ is proper,  there exists a normally rigid structure on $\ex F$ so that $\psi:\ex F\longrightarrow \ex X$ is normally rigid, and each inclusion $\ex A\longrightarrow \ex F$ and $\ex B\longrightarrow \ex F$ is normally rigid on an open neighborhood of $C_{i}$.
\end{lemma}

\pf

 The fact that the closure of $\ex F^{\s}$ contains $\ex F^{\s'}$ if and only if $\s<\s'$ implies that $\s<\s'$ only if the dimension of $\totl{\ex F^{\s}}$ is strictly greater than the dimension of $\totl{\ex F^{\s'}}$--- so $\s'$ is maximal if the dimension of $\totl{\ex F^{\s'}}$ is $0$. 

Suppose that the maps
\[r_{\ex F^{\s'}}:\ex N_{\ex F^{\s'}}\longrightarrow \ex F\]
have already been chosen compatibly for all $\s'>\s$. We shall now construct an appropriate map $r_{\ex F^{\s}}:\ex N_{\ex F^{\s}}\longrightarrow \ex F$.

As well as items (\ref{xc1}) and (\ref{xc2}) from Definition \ref{normally rigid structure}, there are $4$ different compatibility conditions which $r_{\ex F^{\s}}$ must satisfy. After enumerating these compatibility conditions below, we shall consider them in more detail.
\begin{enumerate}
\item \label{r1}$r_{\ex F^{\s}}$ must be compatible with the map 
\[r_{\ex A^{\s}}:\ex N_{\ex A^{\s}}\longrightarrow \ex A\]
\item \label{r2} $r_{\ex F^{\s}}$ must be compatible with the maps 
\[r_{\ex B^{\s'}}:\ex N_{\ex B^{\s'}}\longrightarrow \ex B\] for all $\s'$ so that the image of $\ex B^{\s'}$ intersects $\ex F^{\s}$. 
\item \label{R3}$r_{\ex F^{\s}}$ must be compatible with the previously constructed $r_{\ex F^{\s'}}$ for all $\s'>\s$.
\item\label{R4} and $r_{\ex F^{\s}}$ must be compatible with the condition that $\psi:\ex F\longrightarrow \ex X$ is to be a normally rigid map.
\end{enumerate}

Condition $\ref{r1}$ determines $r_{\ex F^{\s}}$ on a neighborhood of $\ex A\cap \ex F^{\s}$. Note that such an $r_{\ex F^{\s}}$ is $G$-equivariant and automatically satisfies items  (\ref{xc1}) and (\ref{xc2}) from Definition \ref{normally rigid structure} because $r_{\ex A^{\s}}$ does. There is, however, no guarantee that a map defined on an open neighborhood will extend, which is why we must reduce the size of the region where condition \ref{r1} holds to an open neighborhood of a  subset $C_{0}$ of $\ex A$ which is a proper subset of $\ex F$. 

\

 We shall now examine condition $\ref{r2}$. Let $\s'$ be so that $\ex B^{\s'}$ has image which intersects $\ex F^{\s}$. Note that the intersection  of the image of  $\ex B^{\s'}$ with $\ex F^{\s}$ will always consist of the image of some number of connected components of $\ex B^{\s'}$. To reduce notational complexity we shall use $\ex B^{\s'}$ to refer to these components which are sent to $\ex F^{\s}$.

Corresponding to the commutative diagram of $G$-equivariant maps,
\[\begin{tikzcd}\ex B\times_{\ex F}\ex A\rar \dar &\ex B \dar
\\ \ex A\rar & \ex F\end{tikzcd}\]
the commutative diagram
\[\begin{tikzcd}\ex B^{\s'}\times_{\ex F}\ex A^{\s}\rar \dar &\ex B^{\s'} \dar
\\ \ex A^{\s}\rar & \ex F^{\s}\end{tikzcd}\]
 extends to a unique commutative diagram of   fiberwise rigid maps 
\[\begin{tikzcd}\ex N_{\ex B^{\s'}}\rvert_{\ex A^{\s}}\rar \dar &\ex N_{\ex B^{\s'}} \dar
\\ \ex N_{\ex A^{\s}}\rar & \ex N_{\ex F^{\s}}\end{tikzcd}
\]
with derivatives determined by the original commutative diagram. (Recall that $T\ex N_{\ex F^{\s}}$ restricted to $\ex F^{\s}$ is canonically isomorphic to $T\ex F$ restricted to $\ex F^{\s}$.) Condition \ref{r2} requires that the diagram 
\[\begin{tikzcd}\ex N_{\ex B^{\s'}}\dar \rar{r_{\ex B^{\s'}}}&\ex B\dar
\\ \ex N_{\ex F^{\s}}\ar{r}{r_{\ex F^{\s}}}&\ex F\end{tikzcd}\]
commutes when $r_{\ex B^{\s'}}$ is restricted to a neighborhood of $C_{1}$ intersected with $\ex B^{\s'}$. Because the map $\ex N_{\ex B^{\s'}}\longrightarrow \ex N_{\ex F^{\s}} $ extends the map $\ex B^{\s'}\longrightarrow \ex F^{\s}$ and has derivative at $\ex B^{\s'}$ which may be identified with the derivative of $\ex B\longrightarrow \ex F$ using the isomorphism (\ref{ci}), there locally exists a map $r_{\ex F^{\s}}$ satisfying condition \ref{r2} and items (\ref{xc1}) and (\ref{xc2}) from Definition \ref{normally rigid structure}. This map may not be $G$-equivariant, but it may be modified to a $G$-equivariant map satisfying the same conditions using Lemma \ref{G modify}. 

 The condition that the normally rigid structures on $\ex A$ and $\ex B$ are compatible implies that the diagram 
\[\begin{tikzcd}\ex N_{\ex A^{\s}}\dar{r_{\ex A^{\s}}}&\lar \ex N_{\ex B^{\s'}}\rvert_{\ex A^{\s}} \dar{r_{\ex B^{\s'}}}
\\ \ex A&\lar\ex B\times_{\ex F} \ex A \end{tikzcd}\]
commutes when restricted to appropriate neighborhoods, therefore on regions where condition \ref{r1} holds, condition \ref{r2} holds automatically, so conditions \ref{r1} and \ref{r2} are compatible.

\

Now consider condition \ref{R3}. 
Suppose that we have chosen neighborhoods $U^{\s'}$ of $\ex F^{\s'}\subset\ex N_{\ex F^{\s'}}$ for all $\s'>\s$, so that the conditions of Definition \ref{normally rigid structure} are satisfied for all $\s'>\s$.

Then there is a unique normally rigid map 
\[\ex N_{\ex F^{\s}}\rvert_{r(U^{\s'})}\xrightarrow{\Phi_{\s<\s'}} \ex N_{\ex F^{\s'}}\]
so that its composition with $r_{{\ex F}^{\s'}}$ is the identity on  $\ex F^{\s}$ and has derivative equal to the identity restricted to $\ex F^{\s}$. Condition \ref{R3} requires that the diagram 
\[\begin{tikzcd}\ex N_{\ex F^{\s}}\rvert_{r(U^{\s'})}\ar{dr}{r_{\ex F^{\s}}}\rar{\Phi_{\s<\s'}}& \ex N_{\ex F^{\s'}}\dar{r_{\ex F^{\s'}}}
\\ & \ex F\end{tikzcd}\]
commutes when the upper left hand corner is restricted to some neighborhood of $\ex F^{\s}\cap r(U^{\s'})$.

We may choose our open subsets $U^{\s'}$ so that the corresponding diagrams
\[\begin{tikzcd}\ex N_{\ex F^{\s}}\rvert_{r(U^{\s'})\cap r(U^{\s''})}\ar{dr}[swap]{\Phi_{\s<\s''}}\rar{\Phi_{\s<\s'}} &\ex N_{\ex F^{\s'}}\rvert_{r(U^{\s''})}\dar[swap]{\Phi_{\s'<\s''}} \ar{dr}{r_{\ex F^{\s'}}}
\\ & \ex N_{\ex F^{\s''}}\rar{r_{\ex F^{\s''}}}&\ex F\end{tikzcd}\]
commute when the domains of $r$  are restricted to appropriate open subsets, therefore the conditions corresponding to different $\s'>\s$ used in condition \ref{R3} are mutually compatible and determine $r_{\ex F^{\s}}$ on a neighborhood of the intersection of $\ex F^{\s}$ with a neighborhood of all these $\ex F^{\s'}$ for $\s'>\s$.

 On the regions where condition \ref{R3} holds, $r_{\ex F^{\s}}$ automatically satisfies Definition \ref{normally rigid structure}, and conditions \ref{r1} and \ref{r2} automatically also hold because of the  following diagrams are commutative and defined when restricted to the appropriate open subsets
 \[\begin{tikzcd}\ex N_{\ex A^{\s}}\rar \dar \ar[bend right=60]{dd}  &\ex N_{\ex F^{\s}}\dar \ar[bend left=60]{dd}
 \\ \ex N_{\ex A^{\s'}}\dar \rar & \ex N_{\ex F^{\s'}}\dar
 \\ \ex A\rar & \ex F\end{tikzcd}\] 
\[\begin{tikzcd}\ex N_{\ex B^{\s''}}\rar \dar \ar[bend right=60]{dd}  &\ex N_{\ex F^{\s}}\dar \ar[bend left=60]{dd}
 \\ \ex N_{\ex B^{\s'''}}\dar \rar & \ex N_{\ex F^{\s'}}\dar
 \\ \ex B\rar & \ex F\end{tikzcd}\]

\
 
 Now consider condition $\ref{R4}$, which determines the composition of $r_{\ex F^{\s}}$ with the submersion $\psi:\ex F\longrightarrow \ex X$ restricted to some neighborhood of $\ex F^{\s}$. In particular, for each $\ex X^{\s_{i}}$, there must exist  an open neighborhood $U_{i}\subset \ex N_{\ex F^{\s}}$ of $\ex F^{\s}\cap\psi^{-1}(\ex X^{\s_{i}})$ and a commutative diagram
 \[\begin{tikzcd}[column sep=small]\ex N_{\ex F^{\s}} \dar{r_{\ex F^{\s}}}&\lar[hook] U_{i}\rar&\ex N_{\ex X^{\s_{i}}}\dar{r_{\ex X^{\s_{i}}}}
 \\ \ex F\ar{rr}{\psi}&&\ex X\end{tikzcd}\]
where the top line consists of fiberwise rigid maps, and the image if $U_{i}$ in $\ex N_{\ex X^{\s_{i}}}$ is contained in a neighborhood of $\ex X^{\s_{i}}$ on which $r_{\ex X^{\s_{i}}}$ is injective. The fiberwise rigid  map $U_{i}\longrightarrow \ex N_{\ex X^{\s_{i}}}$ is entirely determined by the existence of such a diagram, so $\psi\circ r_{\ex F^{\s}}$ is determined by the above commutative diagram on $U_{i}$. If $\s_{1}<\s_{2}$, then there is some diagram
\[\begin{tikzcd}\ex N_{\ex X^{\s_{1}}}\rar{\Phi_{\s_{1}<\s_{2}}} \dar{r} &\ex N_{\ex X^{\s_{2}}}\ar{dl}{r}
\\ \ex X  \end{tikzcd}\]
which is commutative and defined restricted to some neighborhood of $\ex X^{\s_{1}}$ within $\ex N_{\ex X^{\s_{1}}}$ intersected with the inverse image of some neighborhood of $\ex X^{\s_{2}}$ within $\ex N_{\ex X^{\s_{2}}}$. If $U_{i}$ are chosen small enough, they will have image contained inside these neighborhoods of $\ex X^{\s_{i}}$, so the conditions on $\psi\circ r_{\ex F^{\s}}$ on $U_{1}$ and $U_{2}$ will coincide on $U_{1}\cap U_{2}$. As the union of all such $U_{i}$ covers $\ex F^{\s}$, these compatible conditions determine $\psi\circ r_{\ex F^{\s}}$ on a neighborhood of $\ex F^{\s}\subset \ex N_{\ex F^{\s}}$.

\

With our conditions understood, we are now ready to locally construct a map which satisfies all conditions which $r_{\ex F^{\s}}$ must satisfy apart from  $G$-equivariance. We shall do this case by case depending on which of our conditions must hold.

\begin{claim}\label{r1c} The map $r_{\ex A^{\s}}$ considered as a map from an open subset of $\ex N_{\ex F^{\s}}$ to $\ex F$ obeys all conditions for $r_{\ex F^{\s}}$.
\end{claim}

Condition \ref{r1} holds for $r_{\ex A^{\s}}$ thought of in this way. We showed above that if condition \ref{r1} holds, conditions \ref{r2} and \ref{R3} automatically hold and the map automatically satisfies Definition \ref{normally rigid structure}. It remains to check that condition \ref{R4} also holds. 

 As $\ex A\longrightarrow \ex X$ is a normally rigid map, there is a commutative diagram
\[\begin{tikzcd}[column sep=small]\ex N_{\ex A^{\s}} \dar{r_{\ex A^{\s}}}&\lar[hook] U'_{i}\rar&\ex N_{\ex X^{\s_{i}}}\dar{r_{\ex X^{\s_{i}}}}
 \\ \ex A\ar{rr} &&\ex X\end{tikzcd}\]
therefore, where condition \ref{r1} holds, there is a diagram
\[\begin{tikzcd}[column sep=small]\ex N_{\ex F^{\s}}\dar{r_{\ex F^{\s}}}&\lar \ex N_{\ex A^{\s}} \dar{r_{\ex A^{\s}}}&\lar[hook] U'_{i}\rar&\ex N_{\ex X^{\s_{i}}}\dar{r_{\ex X^{\s_{i}}}}
 \\\ex F\ar[bend right]{rrr}{\psi}& \lar \ex A\ar{rr} &&\ex X\end{tikzcd}\]
 which commutes when restricted to appropriate open subsets, so condition \ref{R4} automatically holds where condition \ref{r1} holds. Therefore, Claim \ref{r1c} is true.
 
 \begin{claim}\label{r2c} Every point in the image of $\ex B^{\s'}$ in  $\ex F^{\s}$ has a neighborhood in $\ex N_{\ex F^{\s}}$ with a map satisfying conditions \ref{r2} and  \ref{R4} and  items (\ref{xc1}) and (\ref{xc2}) from Definition \ref{normally rigid structure}. \end{claim} 
 
Suppose that our chosen point  $p$ in $\ex F^{\s}$ has image in $\ex X^{\s_{1}}$. To prove Claim \ref{r2c}, consider the  commutative diagram
\[\begin{tikzcd}[column sep=small]\ex N_{\ex F^{\s}} &\lar \ex N_{\ex B^{\s'}} \dar{r_{\ex B^{\s'}}}&\lar[hook] U''_{1}\rar&\ex N_{\ex X^{\s_{1}}}\dar{r_{\ex X^{\s_{1}}}}
 \\\ex F\ar[bend right]{rrr}{\psi}& \lar \ex B\ar{rr} &&\ex X\end{tikzcd}\] 
As noted in our discussion of condition \ref{r2} above, we may fill in this diagram with a map $r_{0}$ from a neighborhood $U_{0}$ of  our point $p$ in $\ex N_{\ex F^{\s}}$ to $\ex F$ so that $r_{0}$ satisfies items (\ref{xc1}) and (\ref{xc2}) from Definition \ref{normally rigid structure} and so that the following diagram  commutes:
\[\begin{tikzcd}[column sep=small]U_{0}\dar{r_{0}} &\lar U'''\rar[hook] &\ex N_{\ex B^{\s'}} \dar{r_{\ex B^{\s'}}}&\lar[hook] U''_{1}\rar&\ex N_{\ex X^{\s_{1}}}\dar{r_{\ex X^{\s_{1}}}}
 \\\ex F\ar[bend right]{rrrr}{\psi}& & \ar{ll} \ex B\ar{rr} &&\ex X\end{tikzcd}\]
 In the above,  $U'''$ is the inverse image of $U_{0}\subset \ex N_{\ex F^{\s}}$ in $\ex N_{\ex B^{\s'}}$. Note that in order to extend the map $r_{\ex B^{\s'}}$ to $r_{0}$, we usually will need to choose $U_{0}$ so that the image of $\ex B^{\s'}$ will be proper in $U_{0}$. 
 
So long as $U_{0}$ is chosen small enough, there exists a unique fiberwise rigid map $U_{0}\longrightarrow \ex N_{\ex X^{\s_{1}}}$ so that the diagram
\[\begin{tikzcd} U_{0}\dar{r_{0}}\rar &\ex N_{\ex X^{\s_{1}}}\dar{r_{\ex X^{\s_{1}}}}
\\ \ex F\rar{\psi}&\ex X\end{tikzcd}\]
commutes to first order at $\ex F^{\s}\cap U_{0}$. The same uniqueness property for the fiberwise rigid map $U_{1}''\longrightarrow \ex N_{\ex X^{\s_{1}}}$ implies that the top loop of the diagram
\[\begin{tikzcd}[column sep=small]U_{0}\dar{r_{0}}\ar[bend left]{rrrr} &&\ar{ll} U'''\cap U''_{1} \dar{r_{\ex B^{\s'}}}\ar{rr}&&\ex N_{\ex X^{\s_{1}}}\dar{r_{\ex X^{\s_{1}}}}
 \\\ex F\ar[bend right]{rrrr}{\psi}& & \ar{ll} \ex B\ar{rr} &&\ex X\end{tikzcd}\]
 commutes, therefore  the outer loop also  commutes when $U_{0}$ is  restricted to the image of $U'''\cap U''_{1}$. This image of $U'''\cap U''_{1}$ is some exploded submanifold. As $\psi$ is a submersion, $r_{0}$ restricted to some smaller neighborhood $U$ of $p$ may be modified off the image of $U'''\cap U''_{1}$ to a map $r$ so that the diagram
\[\begin{tikzcd}[column sep=small]U\dar{r}\ar[bend left]{rrrr} &\lar U'''\rar[hook] &\ex N_{\ex B^{\s'}} \dar{r_{\ex B^{\s'}}}&\lar[hook] U''_{1}\rar&\ex N_{\ex X^{\s_{1}}}\dar{r_{\ex X^{\s_{1}}}}
 \\\ex F\ar[bend right]{rrrr}{\psi}& & \ar{ll} \ex B\ar{rr} &&\ex X\end{tikzcd}\]
commutes, and so that $r$ is equal to $r_{0}$ to first order on $U\cap \ex F^{\s}$, and is equal to $r_{0}$ on the intersection of $U$ with the image of $\ex N_{\ex B^{\s'}}$.  Our map $r$ now satisfies both conditions \ref{r2} and \ref{R4}, and items (\ref{xc1}) and (\ref{xc2}) of Definition \ref{normally rigid structure}. This completes the proof of Claim \ref{r2c}.

\

\begin{claim}\label{r3c}The map 
\[r_{\ex F^{\s'}}\circ\Phi_{\s<\s'}:\ex N_{\ex F^{\s}}\rvert_{r(U^{\s'})}\longrightarrow \ex F\]
satisfies all the conditions required of $r_{\ex F^{\s}}$. \end{claim}

In our discussion of condition \ref{R3} above, we have already shown that $r_{\ex F^{\s'}}\circ\Phi_{\s<\s'}$ satisfies Definition \ref{normally rigid structure}, and conditions \ref{r1}, \ref{r2} and \ref{R3}. It remains to verify that condition \ref{R4} also holds.
 As the maps $r_{\ex F^{\s'}}$ have been chosen compatibly with the normally rigid structure on $\ex X$, for $\s'>\s$ there are diagrams
\[\begin{tikzcd}\ex N_{\ex F^{\s}}\rvert_{r(U^{\s'})}\rar \ar{dr}{r_{\ex F^{\s}}} &\ex N_{\ex F^{\s'}}\dar{r_{\ex F^{\s'}}}
&\lar[hook] U'''_{i}\rar&\ex N_{\ex X^{\s_{i}}}\dar{r_{\ex X^{\s_{i}}}}
\\ & \ex F\ar{rr}{\psi}&&\ex X\end{tikzcd}\]
which commute when restricted to the appropriate open subsets, so condition \ref{R4} holds automatically wherever condition \ref{R3} holds. This completes the proof of Claim \ref{r3c}.

\begin{claim}\label{r4c} Any point in $\ex F^{\s}$ has a neighborhood with a map satisfying condition \ref{R4} and items (\ref{xc1}) and (\ref{xc2}) of Definition \ref{normally rigid structure}.
\end{claim}

Claim \ref{r4c} is easier to prove than Claim \ref{r2c}. A map $r$ satisfying items (\ref{xc1}) and (\ref{xc2}) of Definition \ref{normally rigid structure} locally exists. Such a map obeys condition \ref{R4} to first order at $\ex F^{\s}$, therefore it may be modified to a map which satisfies condition \ref{R4} and which is equal to $r$ to first order at $\ex F^{\s}$, and  therefore satisfies 
items (\ref{xc1}) and (\ref{xc2}) of Definition \ref{normally rigid structure}. This completes the proof of Claim \ref{r4c}.

\

The maps from Claims \ref{r1c}, \ref{r2c}, \ref{r3c} and \ref{r4c} all agree to first order on $\ex F^{\s}$, are $G$-equivariant to first order on $\ex F^{\s}$, and all satisfy condition \ref{R4}, so they agree and are $G$-equivariant on a neighborhood of $\ex F^{\s}$ when composed with $\psi:\ex F\longrightarrow \ex X$. The maps from Claims \ref{r1c}, \ref{r2c}, \ref{r3c} all agree and are $G$-equivariant on the image of $\ex N_{B^{\s'}}$ intersected with a neighborhood of $\ex F^{\s}$. The maps from Claims \ref{r1c} and \ref{r3c} agree and are $G$-equivariant. 

We may therefore patch together our maps from Claims \ref{r1c}, \ref{r2c}, \ref{r3c}, and \ref{r4c} to a single map
\[r:\ex N_{\ex F^{\s}}\longrightarrow \ex F\]
which satisfies items (\ref{xc1}) and (\ref{xc2}) of Definition \ref{normally rigid structure}, and which

 \begin{itemize}
 \item satisfies condition \ref{R3} and is $G$-equivariant for some slightly shrunken neighborhoods $U^{\s'}$, (for this, the maps from Claim \ref{r3c} must be used),
 \item satisfies condition \ref{r1} and is $G$-equivariant on a neighborhood of the image of  $C_{0}\subset \ex A$, (for this, the maps from Claim \ref{r1c} and \ref{r3c} must be used,)
  \item satisfies condition \ref{r2} on a neighborhood of the image of $ C_{1}\subset \ex B$, and is $G$-equivariant restricted to the image of $\ex N_{\ex B^{\s'}}$ intersected with this neighborhood, (for this, the maps from Claims \ref{r1c}, \ref{r2c}, and \ref{r3c} may be used,)
\item satisfies condition \ref{R4} and has $G$-equivariant composition with $\psi$, (all the maps from Claims \ref{r1c}, \ref{r2c}, \ref{r3c}, and \ref{r4c} may be used for this.)
\end{itemize}

Applying Lemma \ref{G modify} to  $r$, then $G$-equivariantly  reparametrizing outside of a neighborhood of $\ex F^{\s}\subset \ex N_{\ex F^{\s}}$ gives a $G$-equivariant map
\[ r_{\ex F^{\s}}\longrightarrow \ex F\] 
which agrees with $r$ to first order at $\ex F^{\s}$ and hence satisfies items (\ref{xc1}) and (\ref{xc2}) of Definition \ref{normally rigid structure}, and which
\begin{itemize}
\item agrees with $r$ on the shrunken neighborhoods $U^{\s'}$ and therefore obeys condition \ref{R3} and item (\ref{xc3}) of Definition \ref{normally rigid structure}.
\item  agrees with $r$ and hence satisfies condition \ref{r1} on a neighborhood of the image of $C_{0}\subset \ex A$, 
\item agrees with $r$ on the intersection of the image of $\ex N_{\ex B^{\s'}}$ with a neighborhood of $C_{1}\subset \ex B$, and hence satisfies condition \ref{r2} on this neighborhood
\item agrees with $r$ on a neighborhood of $\ex F^{\s}$ after composition with $\psi$, and therefore obeys condition \ref{R4}.
\end{itemize}

This is our required $r_{\ex F^{\s}}$. As the dimension of $\totl{\ex F^{\s}}$ is strictly greater than the dimension of $\ex F^{\s'}$ if $\s<\s'$, we may construct all $r_{\ex F^{\s}}$ for $\s$ so that the dimension of $\totl{\ex F^{\s}}$ is $0$, then inductively  construct all $r_{\ex F^{\s}}$ for $\s$ so the dimension of $\totl{\ex F^{\s}}$ is $k$. This completes the construction of a normally rigid structure on $\ex F$ compatible with the normally rigid structures on $\ex A$, $\ex B$, and $\ex X$.

\stop

\subsection{Normally rigid vector bundles}

\begin{defn}\label{nrvb}Suppose that $(\ex F,G)$ has a normally rigid structure. A normally rigid structure on a $G$-equivariant, complex vector bundle $V\longrightarrow \ex F$ is for each strata $\ex F^{\s}$ of $(\ex F,G)$ a $G$-equivariant map of complex vector bundles
\[\begin{tikzcd}V_{\ex F^{\s}}\dar\rar{l_{\ex F^{\s}}}& V\dar 
\\ \ex N_{\ex F^{\s}}\rar{r_{\ex F^{\s}}}&\ex F\end{tikzcd}\]
so that
\begin{itemize}
\item $ V_{\ex F^{\s}}$ is the pullback of the vector bundle $V$ over the map $\ex N_{\ex F^{\s}}\longrightarrow \totl{\ex F^{\s}}\subset\totl{\ex F}$. (Note that $V$ may be considered as vector bundle over $\totl{\ex F}$ instead of as a vector bundle over $\ex F$.)
\item Let  $V\rvert_{\ex F^{\s}}$ indicate the pullback of $V$ to $\ex F^{\s}$. There is a canonical inclusion of $V\rvert_{\ex F^{\s}}$ into $V_{\ex F^{\s}}$.  The map $l_{\ex F^{\s}}$ is the identity restricted to the image of $V\rvert_{\ex F^{\s}}$ in the sense that the diagram
\begin{equation}\label{lc1}\begin{tikzcd} V\rvert_{\ex F^{\s}}\rar[hook]\ar[hook, bend left]{rr} \dar &V_{\ex F^{\s}}\dar\rar{l_{\ex F^{\s}}}& V\dar 
\\ \ex F^{\s}\rar[hook]\ar[hook, bend right]{rr}& \ex N_{\ex F^{\s}}\rar{r_{\ex F^{\s}}}&\ex F\end{tikzcd}\end{equation}
commutes.
\item There is a neighborhood $U^{\s}$ of each $\ex F^{\s}\subset \ex N_{\ex F^{\s}}$ satisfying the requirements on $U^{\s} $ from Definition \ref{normally rigid structure} so that $l_{\ex F^{\s}}$ is an isomorphism onto its image when restricted to $U^{\s}$, and so that there is a commutative diagram  whenever $\s<\s'$:
\begin{equation}\label{lc2}\begin{tikzcd}
V_{\ex F^{\s}}\rvert_{U^{\s}\cap \Phi^{-1}_{\s<\s'}(U^{\s'})} \dar[hook]\ar[bend right=50]{ddd} \ar{drr}{l_{\ex F^{\s}}}
\\V_{\ex F^{\s}}\rar{\phi_{\s<\s'}}\rvert_{r(U^{\s'})}\dar  &V_{\ex F^{\s'}}\rar[swap]{l_{\ex F^{\s'}}} \dar & V\dar
\\ \ex N_{\ex F^{\s}}\rar{\Phi_{\s<\s'}}\rvert_{r(U^{\s'})}& \ex N_{\ex F^{\s'}}\rar{r_{\ex F^{\s'}}}& \ex F
\\ U^{\s}\cap \Phi^{-1}_{\s<\s'}(U^{\s'})\uar[hook]\ar{rru}{r_{\ex F^{\s}}}
\end{tikzcd}\end{equation}
where, as in Definition \ref{normally rigid structure}, $\Phi_{\s<\s'}$ is a fiberwise rigid map, and $\phi_{\s<\s'}$ is a  map of complex  vector bundles which is a lift of a map of vector bundles over  the map $\totl{\ex F^{\s}}\cap r(U^{\s'})\longrightarrow \totl{\ex F^{\s'}}$ induced from $\Phi_{\s<\s'}$. 

\end{itemize}
\end{defn}

Note that just as $\Phi_{\s<\s'}$ is canonically determined by $r_{\ex F^{\s'}}$, $\phi_{\s<\s'}$ is canonically determined by $l_{\ex F^{\s'}}$. In particular, there is a unique map $\phi_{0}$ so that the following diagram commutes.
\[\begin{tikzcd}V\rvert_{\ex F^{\s}\cap r(U^{\s'})}\dar{\phi_{0}}\rar[hook]& V
\\ V_{\ex F^{\s'}}\ar{ur}[swap]{l_{\ex F^{\s'}}} \end{tikzcd}\]
Our map $\phi_{s<s'}$ is then the unique map of vector bundles over $\Phi_{s<s'}$ so that the following diagram commutes.
\[\begin{tikzcd}V_{\ex F^{\s}}\dar{\phi_{s<s'}}\rar &\totl{V\rvert_{\ex F^{\s}\cap r(U^{\s'})}}\dar{\totl{\phi_{0}}}
\\ V_{\ex F^{\s'}}\rar& \totl{V_{\ex F^{\s'}}} \end{tikzcd}\]

\begin{defn}\label{nrvb map}Given  normally rigid vector bundles $V_{i}$ over $\ex F_{i}$, a map 
\[\begin{tikzcd}V_{1}\rar{\psi}\dar& V_{2}\dar
\\\ex F_{1}\rar{\Psi}&\ex F_{2}
\end{tikzcd}\] 
is normally rigid if the following holds:
\begin{enumerate}
\item $\Psi$ is normally rigid.
\item $\psi$ is an equivariant map of complex vector bundles.
\item For each $\ex F^{\s'}_{1}$ and $\ex F_{2}^{\s}$, there exists\begin{itemize}\item some neighborhood $U\subset \ex N_{\ex F^{\s'}_{1}}$ of $\ex F^{\s'}\cap \Psi^{-1}\ex F_{2}^{\s}$,
\item a fiberwise rigid map $U\longrightarrow \ex N_{\ex F_{2}^{\s}}$ satisfying the requirements of Definition \ref{normally rigid map}
\item and a map of complex vector bundles
\[\begin{tikzcd}(V_{1})_{\ex F_{1}^{\s'}}\rvert_{U}\rar\dar &(V_{2})_{\ex F_{2}^{\s}}\dar
\\ U \rar &\ex N_{\ex F^{\s}_{2}}
\end{tikzcd}\]
which is a lift of a map of vectorbundles over a map $\totl{U\cap \ex F_{1}^{\s'}}\longrightarrow \totl {\ex F^{\s}_{2}}$.
\end{itemize}
so that the following diagram commutes:

\[\begin{tikzcd}(V_{1})_{\ex F_{1}^{\s'}}\rvert_{U}\rar\dar{l_{\ex F_{1}^{\s'}}} &(V_{2})_{\ex F_{2}^{\s}}\dar{l_{\ex F_{2}^{\s}}}
\\ V_{1}\rar{\Psi} &V_{2}
\end{tikzcd}\]
\end{enumerate}
\end{defn}

Clearly, the property of being a normally rigid map of vectorbundles is preserved under composition. If we define a category of normally rigid vectorbundles with normally rigid maps as morphisms, then up to isomorphism, a normally rigid vectorbundle on $\ex F$ only keeps track of the germ of the maps $l_{\ex F^{\s}}$ around  $\ex F^{\s}\subset \ex N_{\ex F^{\s}}$.

\

The following lemma allows us to construct and extend normally rigid vector bundles. 

\begin{lemma}\label{nrv construction} Given
\begin{itemize}
\item a normally rigid structure on $(\ex F,G)$,
\item a $G$-equivariant complex vector bundle $V$ over $\ex F$
\item a normally rigid structure on $V$ restricted to an open subset $\ex A\subset \ex F$
\item a $G$-equivariant section $\dbar:\ex F\longrightarrow V$
\item a complex, $G$-equivariant sub vector bundle $V'$ of $V$ restricted to an open subset $O\subset \ex F$ so that $\dbar$ is transverse to $V'$, 
\item a normally rigid structure on 
\[\ex B:=\dbar\cap V'\subset O\subset \ex F\] which is compatible with the normally rigid structure on $\ex F$,
\item a normally rigid structure on $V'$ restricted to $\ex B$ which is compatible with the normally rigid structure on $\ex A$ in the sense that the inclusion of $V'\rvert_{\ex A\cap \ex B}$  into $V\rvert_{\ex A}$ is normally rigid,
\item and subsets $C_{0}\subset \ex A$ and $C_{1}\subset \ex B$ which are proper as subsets of $\ex F$, 
\end{itemize}
there exists a normally rigid structure on $V$ which is compatible with the given normally rigid structure on $V\rvert_{\ex A}$ restricted to an open neighborhood of $C_{0}\subset \ex A$, and compatible with the normally rigid structure on $V'$ restricted to an open neighborhood of $C_{1}\subset \ex B$. 
\end{lemma}
\pf

The proof is analogous to the proof of Lemma \ref{nr construction}. 

Suppose that the maps  $l_{\ex F^{\s'}}$ have been chosen compatibly for all $\s' >\s$. We shall now construct $l_{\ex F^{\s}}$.

Apart from $G$-equivariance, there are $4$  conditions which we want $l_{\ex F^{\s}}$ to satisfy:

\begin{enumerate}
\item\label{l1} In a neighborhood of $C_{0}\subset \ex A$, we wish our normally rigid structure to be compatible with the already chosen normally rigid structure. Where this condition applies, all other conditions hold.
\item\label{l2} $l_{\ex F^{\s}}$ must be compatible with the previously constructed maps $l_{\ex F^{\s'}}$ for all $\s' > \s$ (so the diagram (\ref{lc2}) from Definition \ref{nrvb} must commute.) Where this condition applies, all other conditions hold. 
\item\label{l3} In a neighborhood of $C_{1}\subset \ex B$, we wish our normally rigid structure to be compatible with the normally rigid structure on $V'$ over $\ex B$. 
\item \label{l4}  The diagram (\ref{lc1}) from Definition \ref{nrvb} must commute.
\end{enumerate}

As in the proof of Lemma \ref{nr construction}, it is easy to locally construct a map $l_{\ex F^{\s}}$ satisfying condition \ref{l1} around any point in $\ex A$.

Now consider condition \ref{l2}. Assume that our maps $l_{\ex F^{\s'}}$ are all compatible with each other when neighborhoods $U^{\s'}$ are used as in Definition \ref{nrvb}. We may then use the commutative diagram \ref{lc2} to define $l_{\ex F^{\s}}$ restricted to $r(U^{\s'})$. 

\begin{claim} \label{sl}The map $l_{\ex F^{\s'}}\circ \phi_{\s<\s'}$  satisfies condition \ref{l2}.

\end{claim}


The uniqueness property of $\phi_{\s<\s'}$ implies that for any $\s''>\s'$, the following diagram commutes
\[\begin{tikzcd}V_{\ex F^{\s}}\rvert_{r(U^{\s'})\cap r(U^{\s''}) }\rar{\phi_{\s<\s'}}\ar{dr}[swap]{\phi_{\s<\s''}}&V_{\ex F^{\s}}\rvert_{r(U^{\s''})}\dar{\phi_{\s'<\s''}}
\\ & V_{\ex F^{\s''}}
\end{tikzcd}\]
If $\s<\s''<\s'$, then a similar diagram commutes, and if $\s''$ is not greater than, less than, or equal to  $\s'$, then we may choose our neighborhoods so that $r(U^{\s'})$ and $r(U^{\s''})$ do not intersect. Diagram (\ref{lc2}) using $\s'$ and $\s''$ then implies that $l_{\ex F^{\s'}}\circ \phi_{\s<\s'}$ is compatible with $l_{\ex F^{\s''}}$ for all $\s''>\s$. This completes the proof of Claim \ref{sl}.

\

Now consider condition \ref{l3}. Suppose that $\ex B^{\s'}$ has nonempty intersection with $\ex F^{\s}$. Each connected component of $\ex B^{\s'}$ is either disjoint from $\ex F^{\s}$ or is sent isomorphically onto a connected component of $\ex F^{\s}\cap\ex B$. As in the proof of Lemma \ref{nr construction},  we shall abuse notation and use $\ex B^{\s'}$ to refer to those components which are sent to $\ex F^{\s}$. 

There is a canonical inclusion $V'\rvert_{\ex B}\hookrightarrow V$ which induces a canonical inclusion $V'_{\ex B^{\s'}}\hookrightarrow V_{\ex F^{\s}}$. Condition \ref{l3} requires that the following diagram commutes,
\[\begin{tikzcd} V'_{\ex B^{\s'}}\rar[hook]\dar{l_{\ex B^{\s'}}}&V_{\ex F^{\s}}\dar{l_{\ex F^{\s}}}
\\ V'\rvert_{\ex B}\rar[hook]& V
\end{tikzcd}\]
which specifies $l_{\ex F^{\s}}$ on the image of $V'_{\ex B^{\s'}}$. 

Condition \ref{l4} specifies that diagram (\ref{lc1}) commutes. As the analogous diagram commutes for  $l_{\ex B^{\s'}}$, around any point in $\ex B^{\s'}$, there exists a map $l_{\ex F^{\s}}$ locally satisfying conditions \ref{l4} and \ref{l3}.

\

Choose neighborhoods $U_{0}^{\s'}$ of $\ex F^{\s'}\subset \ex N_{\ex F^{\s'}}$ with closure contained in $U^{\s'}$.  We may  create a map $l_{\ex F^{\s}}$ satisfying conditions \ref{l1}, \ref{l2}, \ref{R3} and \ref{R4} by patching together
\begin{enumerate}\item maps which satisfy condition \ref{l1} on $\ex A$, \item maps which satisfy condition \ref{l2} on the union of $r(U^{\s'})$ for all $\s'>\s$,
\item maps which satisfy condition \ref{l3} and \ref{l4} on $O$ minus a neighborhood of $C_{0}\subset \ex A$ and minus $ U_{0}^{\s'}$ for all $\s'>\s$ 
\item and maps which satisfy condition \ref{l4} on $\ex F$ minus neighborhoods of $ C_{0}$ and $ C_{1}$ and minus $U_{0}^{\s'}$ for all $\s'>\s$.
 \end{enumerate}

Such a map $l_{\ex F^{\s}}$ obeys condition \ref{l1} on a neighborhood of $C_{0}$, condition \ref{l2} on the union of $r(U^{\s'}_{0})$ for all $\s'>\s$, condition \ref{l3} on a neighborhood of $C_{1}$, and condition \ref{l4} everywhere. Our map $l_{\ex F_{s}}$ must be $G$-equivariant where conditions \ref{l1} and \ref{l2} hold, and be $G$-equivariant restricted to $V'\rvert_{\ex B'}$ where condition \ref{l3} holds, and also be $G$-equivariant restricted to  $V\rvert_{\ex F^{\s}}\subset V_{\ex F^{\s}}$ because of condition \ref{l4}.  By averaging, we may therefore assume that $l_{\ex F^{\s}}$ is $G$-equivariant, and still obeys the same conditions. 

This $l_{\ex F^{\s}}$ obeys all the required compatibility conditions. We may therefore continue the construction until $l_{\ex F^{\s}}$ has been chosen compatibly for all $\s$. 

\stop

\subsection{Normally rigid structure on the stack $\ex F/G$}

\

We wish to construct a normally rigid structure on an embedded Kuranishi structure. Our Kuranishi charts use exploded manifolds $\ex F$ with an action of a finite group $G$, however the intrinsic thing being described is the stack $\ex F/G$. If $\ex F/G$ is equivalent to $\ex F'/G'$, then $\ex F$ and $\ex F'$ may be related by  taking a fiber product
\[\begin{tikzcd}\ex F'\times_{\ex F/G}\ex F \rar\dar &\ex F'\dar
\\ \ex F\rar & \ex F/G\end{tikzcd}\]
where $\ex F'\times_{\ex F/G}\ex F$ is an exploded manifold with an action of $G\times G'$, and may be regarded as a $G$-equivariant $G'$ bundle over $\ex F$, or a $G'$ equivariant $G$-bundle over $\ex F'$. The fiber product $\ex F'\times_{\ex F/G}\ex F$ is perhaps easiest to understand in the case at hand that $\ex F$ and $\ex F'$ are families of curves--- then a curve in $\ex F'\times_{\ex F/G}\ex F$ is equivalent a given curve with a choice of isomorphism with a curve in $\ex F$ and a curve in $\ex F'$.

Suppose that $(\ex F_{0},G\times G')$ is a $G$-equivariant $G'$ bundle over $(\ex F,G)$. Below, we shall show that a choice of normally rigid structure on $(\ex F_{0},G\times G')$ is equivalent to a choice normally rigid structure on $(\ex F,G)$.  As a corollary, we shall then show that if $\ex F_{1}/G_{1}=\ex F_{2}/G_{2}$ as stacks, then a choice of normally rigid structure on $(\ex F_{1},G_{1})$ is equivalent to a choice of normally rigid structure on $(\ex F_{2},G_{2})$. 
\

\begin{lemma}\label{isotropy preserve} Suppose that $\pi:\ex F_{0}\longrightarrow \ex F$ is a $G$-equivariant $G'$-bundle. Then $\pi$ preserves weak isotropy groups in the following sense: The homomorphism $G\times G'\longrightarrow G$ sends the weak isotropy group $I_{x}$ of $x\in \ex F_{0}$ isomorphically onto the weak isotropy group $I_{\pi(x)}$ of $\pi(x)\in \ex F$.
\end{lemma}

\pf 

Observe that the fact that $\pi$ is $G$-equivariant and $G'$-invariant implies that the image of $I_{x}$ must be contained in $I_{\pi(x)}$. As $\ex F_{0}$ is a $G'$-fold cover of $\ex F$, only the trivial element of $G'$ is contained in $I_{x}$, so the map $I_{x}\longrightarrow I_{\pi(x)}$ is injective. 

It remains to prove that the map $I_{x}\longrightarrow I_{\pi(x)}$ is surjective. The action  of $(g,\id)$ for  $g\in I_{\pi(x)}$  sends $x\in \ex F_{0}$ to some point topologically equivalent to some $x'$ in  $\pi^{-1}(\pi(x))$. As $\pi$ is a $G'$-fold cover, there exists some element $h\in G'$ which sends $x'$ to $x$. Therefore, $(g,h)\in I_{x}$. It follows that $I_{x}\longrightarrow I_{\pi(x)}$ is surjective. 

\stop
 
 \

The fact that $\pi:\ex F_{0}\longrightarrow \ex F$ is a covering map implies that a neighborhood of $x$ with the action of $I_{x}$ is isomorphic to a neighborhood of $\pi(x)$ with the action of $I_{\pi(x)}$. Therefore $\pi$ sends $(\ex F_{0},G\times G')^{\s}$ to $(\ex F,G)^{\s}$, and 
\[\pi:\ex F_{0}^{\s}\longrightarrow \ex F^{\s}\]
is a $G$-equivariant $G'$-bundle which extends canonically to a fiberwise rigid map
\[\ex N_{\ex F_{0}^{\s}}\longrightarrow \ex N_{\ex F^{\s}}\]
which is again a $G$-equivariant $G'$-bundle.

\

We can now state how normally rigid structures on $(\ex F_{0},G\times G')$ and $(\ex F,G)$ are related:

\begin{lemma}\label{qequiv} If $(\ex F_{0},G\times G')$ is a $G$-equivariant $G'$-bundle over $(\ex F,G)$, then given any choice of normally rigid structure on $\ex F_{0}$, there exists a unique choice of normally rigid structure on $\ex F$ so that the following diagrams commute for each $\s$:
\begin{equation}\label{rcom}\begin{tikzcd}\ex N_{\ex F_{0}^{\s}}\dar\rar{r_{\ex F^{\s}_{0}}}&\ex F_{0}\dar
\\ \ex N_{\ex F^{\s}}\rar{r_{\ex F^{\s}}}&\ex F\end{tikzcd}\end{equation}
Conversely, given any choice of normally rigid structure on $\ex F$, there exists a unique choice of normally rigid structure on $\ex F_{0}$ so that the above diagrams commute.

\end{lemma}

\pf

As the downward pointing arrows in the diagram (\ref{rcom}) are $G'$-fold covering maps, given  any choice of $G\times G'$-equivariant map $r_{\ex F^{\s}_{0}}:\ex N_{\ex F^{\s}_{0}}\longrightarrow \ex F_{0}$, there exists a unique $G$-equivariant map $r_{\ex F^{\s}}$ so that the diagram (\ref{rcom}) commutes. These induced maps automatically satisfy  conditions (\ref{xc1}) and (\ref{xc2})  of Definition \ref{normally rigid structure} if $r_{\ex F^{\s}_{0}}$ does, and the uniqueness condition implies that they satisfy condition (\ref{xc3}) if the maps $r_{\ex F^{\s}_{0}}$ do.

Conversely, given any $G$-equivariant map $r_{\ex F^{\s}}:\ex N_{\ex F^{\s}}\longrightarrow \ex F$ which is the identity restricted to $\ex F^{\s}\subset \ex N_{\ex F^{\s}}$, there exists a unique $G\times G'$-equivariant map $r_{\ex F^{\s}_{0}}$ which is the identity restricted to $\ex F^{\s}_{0}\subset \ex N_{\ex F^{\s}_{0}}$ so that diagram \ref{rcom} commutes. This lifted map satisfies condition (\ref{xc1}) of Definition \ref{normally rigid structure}, and satisfies condition (\ref{xc2}) if $r_{\ex F^{\s}}$ does. Again, the uniqueness of these lifted maps implies that they satisfy condition (\ref{xc3}) of Definition \ref{normally rigid structure} if the maps $r_{\ex F^{\s}}$ do.

\stop

\begin{remark}Similarly, if $V$ is a $G$-equivariant vector bundle over $\ex F$, and $(\ex F_{0},G\times G')$ is a $G$-equivariant $G'$-bundle over $(\ex F,G)$, then given any normally rigid structure on $V$ or the pullback $V'$ of $V$ to $\ex F_{0}$,  there exists a unique normally rigid structure on $V'$ or $V$ respectively  so that $V'\longrightarrow V$ is a normally rigid map of vector bundles.\end{remark}

 Suppose that as stacks, $\ex F_{1}/G_{1}=\ex F_{2}/G_{2}$. The definition of the stack $\ex F_{2}/G_{2}$ implies that the quotient map $\ex F_{1}\longrightarrow \ex F_{1}/G_{1}=\ex F_{2}/G_{2}$ is equivalent to a $G_{2}$ bundle $\ex F$ over $\ex F_{1}$ with a $G_{2}$-equivariant map $\ex F\longrightarrow \ex F_{2}$. In fact, $\ex F$ may be regarded as the fiber product of $\ex F_{1}$ with $\ex F_{2}$ over $\ex F_{1}/G_{1}$ 
\[\begin{tikzcd}\ex F\rar\dar&\ex F_{2}\dar
\\ \ex F_{1}\rar&\ex F_{1}/G_{1}=\ex F_{2}/G_{2}\end{tikzcd}\]
The $G_{1}$ action on $\ex F_{1}$ lifts to a $G_{1}$ action on $\ex F$ which commutes with the $G_{2}$ action on $\ex F$, so that the top map above is a $G_{1}$-invariant, $G_{2}$ equivariant $G_{1}$-fold cover, and the left hand map above is a $G_{2}$-invariant, $G_{1}$-equivariant $G_{2}$-fold cover.  In particular, \[(\ex F,G_{1}\times G_{2})/G_{2}=(\ex F_{1},G_{1})\text{ and }(\ex F,G_{1}\times G_{2})/G_{1}=\ex F_{2}\]

With the fiber product $\ex F$ understood, the following is now an immediate corollary of Lemma \ref{qequiv}

\begin{cor} If $\ex F_{1}/G_{1}=\ex F_{2}/G_{2}$ as stacks, then a choice of normally rigid structure on $(\ex F_{1},G_{1})$ is equivalent to a choice of normally rigid structure on $(\ex F_{2},G_{2})$. In particular, if $\ex F$ indicates the fiber product of $\ex F_{1}$ with $\ex F_{2}$ over $\ex F_{1}/G_{1}$, then given any choice of normally rigid structure on $\ex F_{1}$ or $\ex F_{2}$, there exists a unique choice of normally rigid structure on $\ex F$ and $\ex F_{2}$ or $\ex F_{1}$ respectively so that the following diagrams commute for all $\s$.

\[\begin{tikzcd}\ex N_{\ex F_{2}^{\s}}\rar{r_{\ex F^{\s}_{2}}}&\ex F_{2}
\\ \ex N_{\ex F^{\s}}\uar\dar\rar{r_{\ex F^{\s}}}&\ex F\uar\dar
\\ \ex N_{\ex F^{\s}_{1}}\rar{r_{\ex F^{\s}_{1}}}&\ex F_{1}
\end{tikzcd}\]
\end{cor}

As in \cite{evc}, we shall use `exploded orbifold' to mean a Deligne-Mumford exploded stack, which is a stack over the category of exploded manifolds which is locally equivalent to $\ex F/G$ for some exploded manifold $\ex F$ with the action
of a finite group $G$.

\begin{defn}A normally rigid structure on an exploded orbifold $\ex X$ is a choice of (equivalence class of) normally rigid structure on $(U,G)$ for any identification of $U/G$ with an open substack of $\ex X$ so that in all fiber product diagrams,
\[\begin{tikzcd}\ex U_{1}\times_{\ex X}\ex U_{2}\rar\dar&\ex U_{2}\dar
\\\ex U_{1}\rar&\ex X\end{tikzcd}\]
the left hand and top arrows are normally rigid maps.

\

Say that a map $h:\ex X\longrightarrow \ex Y$ between exploded orbifolds with normally rigid structures is normally rigid if $\ex X$ has a cover by open substacks $U/G$ with images contained in open substacks $U'/G'$ of $\ex Y$ so that $h$ lifts to a $G$-equivariant, normally rigid map $(U,G)\longrightarrow (U',G')$. 
\end{defn}

A normally rigid structure can be constructed on any exploded orbifold $\ex X$ using Lemma \ref{nr construction}. The idea is as follows: First, choose a countable,  locally finite cover of $\ex X$ by open substacks $U_{i}/G_{i}$. A normally rigid structure may be chosen on each $(U_{i},G_{i})$ in turn, and this normally rigid structure may be chosen compatible with previous choices on $U_{i'}$ at the expense of shrinking $U_{i'}$ a little. If $U/G$ is isomorphic to an open substack of $\ex X$, there is then a choice of normally rigid structure on $(U,G)$ which is compatible with the normally rigid structures on the (shrunken) $(U_{i},G_{i})$. This normally rigid structure on $(U,G)$ is unique up to a normally rigid isomorphism.

\

\subsection{Normally rigid structures on Kuranishi charts}

\

In \cite{evc}, a Kuranishi structure is constructed to cover the moduli stack of holomorphic curves in a family of exploded manifolds $\hat{\ex B}\longrightarrow \ex B_{0}$.  This Kuranishi structure comes naturally embedded in the moduli stack $\Msw$ of $\C\infty1$ curves on $\ex B$. The Kuranishi structure consists of charts $(\mathcal U_{i},V_{i},\hat f_{i}/G_{i})$ where \begin{itemize}
\item $\mathcal U_{i}$ is an open substack of $\Msw$.
\item $V_{i}$ is a nice obstruction bundle over $\mathcal U_{i}$. (In particular, any family $\hat h$ of curves in $\mathcal U_{i}$ is functorially assigned a  vector bundle $V_{i}(\hat h)$ over $\hat h$ which is a sub vector bundle of the sheaf $\mathcal Y(\hat h)$ over $\hat h$ which $\dbar\hat h$ is a section of.) 
\item $\hat f_{i}$ is a family of curves in $\mathcal U_{i}$,
\item and $G_{i}$ is a group of automorphisms of $\hat f_{i}$
\end{itemize}
so that $\hat f_{i}/G_{i}$ represents the moduli stack of curves $h$ in $\mathcal U_{i}$ with $\dbar h$  in $V_{i}(h)$. (In particular, $\hat f_{i}$ comes with a naturally defined section $\dbar $ of the vector bundle $V_{i}(\hat f_{i})$.)

 Two different Kuranishi carts $(\mathcal U_{i},V_{i},\hat f_{i}/G_{i}), \ (\mathcal U_{j},V_{j},\hat f_{j}/G_{j}) $ must be compatible in the sense that on $\mathcal U_{i}\cap \mathcal U_{j}$, $V_{i}$ is a subbundle of $V_{j}$ or $V_{j}$ is a subbundle of $V_{i}$. If $V_{i}$ is a subbundle of $V_{j}$, then $\hat f_{i}/G_{i}\cap \mathcal U_{j}$ is equal to the substack of $\hat f_{j}/G_{j}\cap \mathcal U_{i}$ where $\dbar$ is in  $V_{i}$. In particular, this implies that the fiber product of $\hat f_{i}$ with $\hat f_{j}$ over $\Msw$ is nice:
  \[\begin{tikzcd} \hat f_{i}\times_{\Msw}\hat f_{j} \dar\rar&\hat f_{j}\dar
\\ \hat f_{i}\rar& \Msw\end{tikzcd}\]
The defining property of $\hat f_{i}\times_{\Msw}\hat f_{j}$ is that it is a family of curves in $\Msw$ with the property that given any other family of curves $\hat h$ in $\Msw$, a map $\hat h\longrightarrow \hat f_{i}\times_{\Msw}\hat f_{j}$ is equivalent to a map $\hat h\longrightarrow \hat f_{i}$ and a map $\hat h\longrightarrow \hat f_{j}$. (Note that as defined in \cite{iec}, a map of families of curves is by definition an isomorphism restricted to each individual curve.) The group $G_{i}\times G_{j}$ acts on $\hat f_{i}\times_{\Msw}\hat f_{j}$ so that the left hand map above is a $G_{i}$-equivariant $G_{j}$-bundle over $\hat f_{i}\cap \mathcal U_{j}$, and so that the top map is a $G_{j}$-equivariant $G_{i}$-fold cover of the subset of $\hat f_{j}\cap \mathcal U_{i}$ where $\dbar\hat f_{j}$ intersects $V_{i}\subset V_{j}$. 

\

 A family of curves $\hat f$ corresponds to the diagram 
\[\begin{tikzcd}\ex C(\hat f)\dar\rar{\hat f} &\hat {\ex B}\dar
\\ \ex F(\hat f)\rar &\ex B_{0}\end{tikzcd}\]
We shall be putting normally rigid structures on the exploded manifolds $\ex F(\hat f_{i})$. Such normally rigid structures will need to be compatible in the sense that there exists a normally rigid structure on $\ex F(\hat f_{i}\times_{\Msw}\hat f_{j})$ so that the natural maps 
\[\begin{tikzcd} \ex F(\hat f_{i}\times_{\Msw}\hat f_{j})\dar\rar&\ex F(\hat f_{j})
\\ \ex F(\hat f_{i})\end{tikzcd}\]
are normally rigid. These compatibility conditions may be `composed' in the following sense: Suppose that $\dim V_{i}\leq \dim V_{j}\leq \dim V_{k}$ and the normally rigid structure on $\ex F(\hat f_{i})$ is compatible with the normally rigid structure on $\ex F(\hat f_{j})$, which is compatible with the normally rigid structure on $\ex F(\hat f_{k})$. Then consider the following commutative diagram:
\[\begin{tikzcd}  & \ex F(\hat f_{i}\times_{\Msw}\hat f_{k}) \ar{dl}\ar{dr}
\\ \ex F(\hat f_{i})& \ex F(\hat f_{i}\times_{\Msw}\hat f_{j}\times_{\Msw}\hat f_{k})\uar \ar{dl}\ar{dr}&\ex F(\hat f_{k})
\\  \ex F(\hat f_{i}\times_{\Msw}\hat f_{j})\uar\rar & \ex F(\hat f_{j})&\ex F(\hat f_{j}\times_{\Msw}\hat f_{k})\lar
\uar \end{tikzcd}\]
The lower internal square of the above diagram is a fiber product diagram, with the left pointing maps being  $G_{k}$-bundles over the open subsets defined by intersecting with $\mathcal U_{k}$. It follows that $\ex F(\hat f_{i}\times_{\Msw}\hat f_{j}\times_{\Msw}\hat f_{k} )$ has a normally rigid structure so that the two downward pointing maps from it are normally rigid. The upward pointing map identifies $\ex F(\hat f_{i}\times_{\Msw}\hat f_{j}\times_{\Msw}\hat f_{k} )$ as a $G_{i}\times G_{k}$-equivariant $G_{j}$-bundle over an open subset $\mathcal U_{j}\cap \ex F(\hat f_{i}\times_{\Msw}\hat f_{k})$, so there is a normally rigid structure on this open subset so that this upwards pointing map is normally rigid.  It follows that the two maps on the top row of this diagram are normally rigid when restricted to the open subset $\mathcal U_{j}\cap \ex F(\hat f_{i}\times_{\Msw}\hat f_{k})$. In other words, the normally rigid structure on $\ex F(\hat f_{i})$ is compatible with the normally rigid structure on $\ex F(\hat f_{k})$ when restricted to the open subset $\mathcal U_{j}$.

 \
 
 In order to construct compatible normally rigid structures on a given set of Kuranishi charts, we shall need to successively shrink the size of the Kuranishi charts we use. To do this, we shall use the notion from \cite{evc} of an extension $(\mathcal U_{i}^{\sharp},V_{i},\hat f^{\sharp}_{i}/G_{i})$ of a Kuranishi chart $(\mathcal U_{i},V_{i},\hat f_{i}/G_{i})$. The important aspects of extensions to keep in mind are  that $\hat f_{i}^{\sharp}$ contains the closure of $\hat f_{i}$ within $\Msw$, and that any extension may be made smaller if necessary. Any embedded Kuranishi structure by definition consists  of a set  Kuranishi charts with compatible extensions which give a locally finite cover of the moduli stack of holomorphic curves.  We shall use the notation 
 \[\hat f\exte \hat f^{\sharp} \ \ \text{ to mean that $\hat f^{\sharp}$ is an extension of $\hat f$}\]

The following lemma is Claim 8.11 from \cite{evc}. It allows us to construct things such as normally rigid structures chart by chart.

\begin{lemma}\label{well order} There is a well-ordering $\prec$ of the Kuranishi charts $(\mathcal U_{i},V_{i},\hat f_{i}/G_{i})$  so that
$j\prec i$ if $\dim V_{j}<\dim V_{i}$,  and so that for any fixed $i$,  $j\prec i$ for only a finite number of $j$ with $\dim V_{j}=\dim V_{i}$.

There are compatible extensions $(\mathcal U_{i,k},V_{i},\hat f_{i,k}/G_{i})$  and $(\mathcal U_{i,\underline k},V_{i},\hat f_{i,\underline k}/G_{i})$ of $(\mathcal U_{i},V_{i},\hat f_{i}/G_{i})$ for all $ k\succeq i$  so that\begin{itemize}
 \item if $k'\prec k$, then $\hat f_{i,k'}$ is an extension of $\hat f_{i,k}$. 
 \[\hat f_{i,\underline k}\exte\hat f_{i,k'}\ \ \ \ \ \ \ \text{ if }\ \ \  k'\prec k\]
 \item \[\hat f_{i,k}\exte\hat f_{i,\underline k}\]
  \item and the intersection of  $\hat f_{i,k}$ for all $k\succeq i$ contains an extension of $\hat f_{i}$.  
\end{itemize}
\end{lemma}

The following lemma allows us to inductively construct compatible normally rigid structures on Kuranishi charts.

\begin{lemma}\label{nri}
Suppose that compatible normally rigid structures have been chosen on $\ex F(\hat f_{i,\underline j})$ for all $i\prec j$. In this case, compatible means that if $i\prec k\prec j$ there is a normally rigid structure on $\ex F(\hat f_{i,\underline j}\times _{\Msw}\hat f_{k,\underline j})$ so that the natural maps 
\[\begin{tikzcd} \ex F(\hat f_{i,\underline j})&\lar \ex F(\hat f_{i,\underline j}\times_{\Msw}\hat f_{k,\underline j})\rar&\ex F(\hat f_{k,\underline j})
\end{tikzcd}\]
are normally rigid.

Then there exists a normally rigid structure on $\ex F(\hat f_{j,j})$ compatible with the normally rigid structure on $\ex F(\hat f_{i,j})\exte \ex F(\hat f_{i,\underline j})$ for all  $i\prec j$ in the sense that  the diagrams
\[\begin{tikzcd}\ex F(\hat f_{i, j})&\lar \ex F(\hat f_{i, j}\times_{\Msw}\hat f_{j, j})\rar&\ex F(\hat f_{j, j})
\end{tikzcd}\]
are normally rigid.

\

Suppose  that there is a $G$-invariant closed subset  $C$ of $\ex F(\hat f_{j,j})$ with an open neighborhood $U$ with a normally rigid structure so that the maps
\[\begin{tikzcd} \ex F(\hat f_{i, j})&\lar \ex F(\hat f_{i, j})\times_{\Msw}U\rar&U
\end{tikzcd}\]
are normally rigid. Then the normally rigid structure on $\ex F(\hat f_{j,j})$ may be chosen to coincide with the normally rigid structure on $U$ in a neighborhood of $C$.

\

Suppose further that there is a surjective submersion $\Msw\longrightarrow \ex X$ so that $\ex X$ is an exploded orbifold with a normally rigid structure and  the maps $\ex F(\hat f_{i})\longrightarrow \ex X$ are submersions for all $i$,  and are normally rigid for all $i\prec j$. Then the normally rigid structure on $\ex F(\hat f_{j,j})$ may be chosen so that $\ex F(\hat f_{j,j}) \longrightarrow \ex X$ is normally rigid.  
\end{lemma}

\pf

As any embedded Kuranishi structure is locally finite, there are only a finite number of indices $i$ so that $\hat f_{j,j}\times_{\Msw}\hat f_{i,\underline j}$ is nonempty. We may therefore choose extensions $ \hat f_{i,j_{k}}\rexte \hat f_{i}$ for all $k\prec j$ with $\hat f_{j,j}\times_{\Msw}\hat f_{k,\underline j}$  nonempty, so that 
\[ \hat f_{i,j}\exte\hat f_{i,j_{k}}\exte \hat f_{i,j_{k'}}\exte \hat f_{i,\underline j}\]
\[\text{ if  } k\prec k'\prec j \]

Fix an $i\prec j$ for which $\hat f_{j,j}\times_{\Msw}\hat f_{i,\underline j}$ is nonempty. Let $i'$ be next largest index with $\hat f_{j,j}\times_{\Msw}\hat f_{i',\underline j}$  nonempty.  Suppose that there is a normally rigid structure on $\ex F(\hat f_{j,j})$ so that 
\begin{itemize} \item for all $i\prec k \prec j$, the diagram
\[\begin{tikzcd}\ex F(\hat f_{k,j_{i'}})&\lar \ex F(\hat f_{k,j_{i'}}\times _{\!\!{}_{{}_{\Msw}}}\hat f_{j,j})\rar&\ex F(\hat f_{j,j} )\end{tikzcd}\] is normally rigid.
\item There exists an open neighborhood $U'$ of $C\subset \ex F(\hat f_{j,j})$ so that the normally rigid structure of $\ex F(\hat f_{j,j})$ coincides with the normally rigid structure on $U$  when restricted to $U'$.
\item The submersion $\ex F(\hat f_{j,j})\longrightarrow \ex X$ is normally rigid.
\end{itemize}
We shall show that the above conditions hold for all $i\prec j$ with $\hat f_{j,j}\times_{\Msw}\hat f_{i,\underline j}$ nonempty  by an induction which proceeds in reverse order of $i$. For the base case $i'=j$, the first condition above is empty and  Lemma \ref{nr construction} implies that such a normally rigid structure exists.

\

Let $\ex A\subset \ex F(\hat f_{j,j})$ be defined by
\[\ex A:=\ex F(\hat f_{j,j})\cap\lrb{U'\bigcup_{i\prec k\prec j}\mathcal U_{k,j_{i'}}}\] 
(Recall that $\mathcal U_{k,j_{i'}}\subset \Msw$ is an open neighborhood of $\hat f_{k,j_{i'}}$ from the definition of the Kuranishi chart $(\mathcal U_{k,j_{i'}},V_{k}, \hat f_{k,j_{i'}}/G_{k})$.)
Let $O\subset \ex F(\hat f_{j,j})$ be defined by 
\[O:=\mathcal U_{i,j_{i'}}\cap \ex F(\hat f_{j,j})\]
and let $\ex B\subset O\subset \ex F(\hat f_{j,j})$ be the exploded submanifold of $O$ corresponding to $\ex F(\hat f_{i,j_{i'}})$. Note that because the normally rigid structures on $\ex F(\hat f_{i,\underline j})$ and $\ex F(\hat f_{k,\underline j})$ are compatible, the normally rigid structures on $\ex B$ and $\ex A$ are compatible,
 so we may apply Lemma \ref{nr construction}. Lemma \ref{nr construction} states that there exists a normally rigid structure on $\ex F(\hat f_{j,j})$ so that the map to $\ex X$ is normally rigid,  the inclusion $\ex A\longrightarrow \ex F(\hat f_{j,j})$ is normally rigid away from the boundary of $\ex A$ within $\ex F(\hat f_{j,j})$, and similarly the inclusion $\ex B\longrightarrow \ex F(\hat f_{j,j})$ is normally rigid away from the boundary of  $\ex B$ within $\ex F(\hat f_{j,j})$. In particular $C$ stays away from the boundary of $\ex A$ within $\ex F(\hat f_{j,j})$, and because $\hat f_{k,j_{i}}\exte\hat f_{k,j_{i'}}$,   the image of $ \ex F(\hat f_{j,j}\times_{\Msw}\hat f_{k,j_{i}})$ within $\ex A$ stays away from the boundary of $\ex A$ within $\ex F(\hat f_{j,j})$. Similarly,  the image of $ \ex F(\hat f_{j,j}\times_{\Msw}\hat f_{i,j_{i}})$ within $\ex B$ stays away from the boundary of $\ex B$ within $\ex F(\hat f_{j,j})$. Therefore the normally rigid structure on $\ex F(\hat f_{j,j})$ may be chosen so that
\begin{itemize}
\item for all $i\preceq k\prec j$ the diagram 
\[\ex F(\hat f_{k,j_{i}})\longleftarrow\ex F(\hat f_{j,j}\times_{\Msw}\hat f_{k,j_{i}})\longrightarrow \ex F(\hat f_{j,j} )\] is normally rigid. 
\item There exists an open neighborhood $U''$ of $C\subset \ex F(\hat f_{ j,j})$ on which the normally rigid structure coincides with the normally rigid structure from $U$.
\item The submersion $\ex F(\hat f_{j,j})\longrightarrow \ex X$ is normally rigid.
\end{itemize}
Inductively continuing this construction for the finite number of $i\prec j$ with $\hat f_{j,j}\times_{\Msw}\hat f_{i,\underline j}$  nonempty, we may construct a normally rigid structure on $\ex F(\hat f_{j,j})$ so that
\begin{itemize}
\item for all $i\prec j$, the diagram
\[\ex F(\hat f_{i,j})\longleftarrow\ex F(\hat f_{i,j}\times_{\Msw}\hat f_{j,j})\longrightarrow \ex F(\hat f_{j,j} )\]
is normally rigid.
\item The normally rigid structure on $\ex F(\hat f_{j,j})$ coincides with the normally rigid structure on $U$ on some neighborhood of $C$.
\item The submersion $\ex F(\hat f_{j,j})\longrightarrow \ex X$ is normally rigid.
\end{itemize}

\stop

The following lemma allows us to compatibly choose a normally rigid structure on the vector bundles $V_{i}$ over $\ex F(\hat f_{i})$.

\begin{lemma}\label{nrvi}Suppose the following:
\begin{itemize}
\item
A normally rigid structure has been chosen on $\ex F(\hat f_{i,i})$ for all $i$ so that if $i\prec j$, then the diagram
\[ \ex F(\hat f_{i,j})\longleftarrow \ex F(\hat f_{i,j}\times_{\Msw}\hat f_{j,j})\longrightarrow \ex F(\hat f_{j,j})\]
is normally rigid.
 \item For all $i\prec j$, a normally rigid structure  has been chosen on the vector bundle $V_{i}$ over $\ex F(\hat f_{i,i})$ so that if $i\prec k \prec j$, the map of normally rigid vectorbundles \[V_{i}(\hat f_{i,k}\times_{\Msw}\hat f_{k,k})\longrightarrow V_{k}(\hat f_{k,k})\] is normally rigid in the sense that  the following diagram
  \[\begin{tikzcd}V_{i}(\hat f_{i,k})\dar&\lar V_{k}(\hat f_{i,k}\times_{\Msw}\hat f_{k,k})\rar\dar&V_{k}(\hat f_{k,k})\dar
  \\ \ex F(\hat f_{i,k})& \ex F(\hat f_{i,k}\times_{\Msw}\hat f_{k,k})\lar\rar&\ex F(\hat f_{k,k}) \end{tikzcd}\]
   is normally rigid.
 \item There is a (possibly empty) $G_{j}$-invariant closed subset $C\subset \ex F(\hat f_{j,j})$ with a $G_{j}$-invariant open neighborhood $U$ and chosen normally rigid structure on $V_{j}$ restricted to $U$ so that for any $i\prec j$, the map  $V_{i}(\hat f_{i,j}\times_{\Msw}U)\longrightarrow V_{j}(U)$  is normally rigid. 
  \end{itemize}

Then there exists a normally rigid structure on $V_{j}$ over $\ex F(\hat f_{j,j})$ which is compatible with the already chosen normally rigid structure on a neighborhood of $C$, and which is compatible with the normally rigid structures on $V_{i}$ restricted to $\hat f_{i,j}\exte\hat f_{i,i}$ for all $i\prec j$ in the sense that  $V_{i}(\hat f_{i,j}\times_{\Msw}\hat f_{j,j})\longrightarrow V_{j}(\hat f_{j,j})$ is normally rigid.
\end{lemma}

\pf

The proof of this lemma is analogous to the proof of  Lemma \ref{nri}, but uses Lemma \ref{nrv construction} in place of  Lemma \ref{nr construction}. In particular, use the same notation $\hat f_{i,j_{k}}$ for the extra extensions used in the proof of Lemma \ref{nri}, so we have
\[ \hat f_{i,j}\exte\hat f_{i,j_{k}}\exte \hat f_{i,j_{k'}}\exte \hat f_{i,\underline j} \ \ \ \text{ if  }\  k\prec k'\prec j \]

 Again, fix an $i\prec j$ for which $\hat f_{i,\underline j}\times_{\Msw}\hat f_{j,j}$ is nonempty, and let $i'$ be next largest index so that $\hat f_{i',\underline j}\times_{\Msw}\hat f_{j,j}$ is nonempty.

The inductive assumption in this case is the following:

  Suppose that there is a normally rigid structure on $V_{j}$ over $\ex F(\hat f_{j,j})$ so that 
\begin{itemize} \item for all $i\prec k \prec j$, the map $V_{k}(\hat f_{k,j_{i'}}\times_{\Msw}\hat f_{j,j})\longrightarrow V_{j}(\hat f_{j,j})$  is normally rigid.
\item There exists a $G_{j}$-invariant open neighborhood $U'$ of $C\subset \ex F(\hat f_{j,j})$ so that the normally rigid structure of $V_{j}$ over $\ex F(\hat f_{j,j})$ coincides with the previously chosen normally rigid structure on $V_{j}$  when restricted to $U'$.
\end{itemize}
Again, for the base case $i'=j$, the first condition above is empty, and Lemma \ref{nrv construction} implies that such a normally rigid structure exists.

\

As in the proof of Lemma \ref{nrv construction}, let $\ex A\subset \ex F(\hat f_{j,j})$ be defined by
\[\ex A:=\ex F(\hat f_{j,j})\cap\lrb{U'\bigcup_{i\prec k\prec j}\mathcal U_{k,j_{i'}}}\] 
Let $O\subset \ex F(\hat f_{j,j})$ be defined by 
\[O:=\mathcal U_{i,j_{i'}}\cap \ex F(\hat f_{j,j})\]
 let $V'$ be the vector bundle over  $O$ obtained by pulling back the vectorbundle  $V_{i}$ over  the map $O\longrightarrow \mathcal U_{i}$. The section $\dbar$ of $V_{j}$ intersects $V'$ transversely over  $\ex B\subset O\subset \ex F(\hat f_{j,j})$, which is the exploded submanifold of $O$ corresponding to $\ex F(\hat f_{i,j_{i'}})$. Because $\hat f_{i,\underline j}\exte\hat f_{i,k}$ and $\hat f_{k,\underline j}\exte \hat f_{k,k}$,  the normally rigid structures on $V_{i}$ over $\ex F(\hat f_{i,\underline j})$ and $V_{k}$ over $\ex F(\hat f_{k,\underline j})$ are compatible. Therefore, the normally rigid structures of $V'$ over $\ex B$ and $V_{j}$ over $\ex A$ are compatible, so we may apply Lemma \ref{nrv construction}. Lemma \ref{nrv construction} states that there exists a normally rigid structure on $V$ over $\ex F(\hat f_{j,j})$   compatible with the previous choice of normally rigid structure on $\ex A$ away from the boundary of $\ex A$ within $\ex F(\hat f_{j,j})$, and compatible with the normally rigid structure on $V'$ over the inclusion $\ex B\longrightarrow \ex F(\hat f_{j,j})$  away from the boundary of  $\ex B$ within $\ex F(\hat f_{j,j})$. In particular, $C$ stays away from the boundary of $\ex A$ within $\ex F(\hat f_{j,j})$, and because $\hat f_{k,j_{i}}\exte \hat f_{k,j_{i'}}$, the image of $ \ex F(\hat f_{k,j_{i}}\times_{\Msw}\hat f_{j,j})$ within $\ex A$ stays away from the boundary of $\ex A$ within $\ex F(\hat f_{j,j})$. Similarly,  the image of $\ex F(\hat f_{i,j_{i}}\times_{\Msw}\hat f_{j,j})$ within $\ex B$ stays away from the boundary of $\ex B$ within $\ex F(\hat f_{j,j})$. Therefore the normally rigid structure on $V_{j}$ over $\ex F(\hat f_{j,j})$ may be chosen so that
\begin{itemize}
\item for all $i\preceq k\prec j$, the map  
\[ V_{k}(\hat f_{k,j_{i}}\times_{\Msw}\hat f_{j,j})\longrightarrow V_{j}(\hat f_{j,j} )\] is normally rigid. 
\item There exists a $G_{j}$-invariant open neighborhood $U''$ of $C\subset \ex F(\hat f_{ j,j})$ on which the normally rigid structure on $V'$ coincides with the previously chosen normally rigid structure.
\end{itemize}
Continuing this construction for the finite number of $i\prec j$ with $\hat f_{i,\underline j}\times_{\Msw}\hat f_{j,j}$  nonempty, we may construct a normally rigid structure on $V_{j}$ over  $\ex F(\hat f_{j,j})$ so that
\begin{itemize}
\item for all $i\prec j$, the map
\[V_{i}(\hat f_{i,j}\times_{\Msw}\hat f_{j,j})\longrightarrow V_{j}(\hat f_{j,j} )\]
is normally rigid.
\item The normally rigid structure on $V_{j}$ over $\ex F(\hat f_{j,j})$ coincides with the previously chosen normally rigid structure  on some neighborhood of $C$.
\end{itemize}

\stop

\subsection{Relatively complex normally rigid structures}

\

In this section, we shall construct a kind of complex structure on our normally rigid structure on $\ex F(\hat f_{i})$. Recall that $\hat f_{i}$ is a family of curves in some family of exploded manifolds. \[\hat {\ex B}\longrightarrow \ex B_{0}\]
In \cite{evc}, we constructed $\hat f_{i}$ so that the map $\ex F(\hat f_{i})\longrightarrow \ex B_{0}$ is a submersion. Use the notation $T\ex F(\hat f_{i})\ov{\ex B_{0}}$ to indicate the vertical tangent space of $\ex F(\hat f_{i})$ over $\ex B_{0}$.

  Nice evaluation  maps $\Phi$ from $\Msw(\hat {\ex B})$ have targets equal to families of exploded manifolds or orbifolds with a family of almost complex structures
 \[\begin{tikzcd}\Msw(\hat {\ex B})\rar{\Phi}\ar{dr}&\ex X\dar
 \\ & \ex B_{0}\end{tikzcd}\]
 The best kind of evaluation maps are holomorphic submersions  in the sense that at any holomorphic curve $f$, the tangent map
 \[T_{f}\Phi :T_{f}\Msw(\hat{\ex B})\ov{\ex B_{0}}\longrightarrow T_{\Phi(f)}\ex X\ov{\ex B_{0}}\]
 is complex linear and surjective.  (For example, $\ex X$ may be the product of Deligne-Mumford space with some fiber product of $\hat{\ex B}$ with itself over $\ex B_{0}$. Then an appropriate $\Phi$ records the complex structure of the domain of a curve as well as the image in $\hat{\ex B}$ of some marked points on the curve.) Given a choice of holomorphic submersion $\Phi$, we also constructed our $\hat f_{i}$ in \cite{evc} so that $\Phi:\ex F(\hat f_{i})\longrightarrow \ex X$ was  a submersion.

\begin{remark}\label{J_{0} conditions}  In \cite{evc}, a homotopy of $D\dbar$ to a complex operator was used to construct a  complex structure $J_{0}$ on $T_{f}\ex F(\hat f_{i})\ov{\ex B_{0}}$ at holomorphic curves $f$ satisfying the following conditions 
\begin{itemize}
\item If $V_{i}(f)\subset V_{j}(f)$, there is a canonical complex linear short exact sequence
 \[0\longrightarrow T_{f}\ex F(\hat f_{i})\ov{\ex B_{0}}\longrightarrow T_{f}\ex F(\hat f_{j})\ov{\ex B_{0}}\xrightarrow {D\pi_{V_{i}}\dbar} V_{j}(f)/V_{i}(f)\longrightarrow 0\] 
\item For a chosen complex submersion $\Phi:\Msw\longrightarrow \ex X$, the map
 \[T_{f}\Phi: T_{f}\ex F(\hat f_{i})\ov{\ex B_{0}}\longrightarrow T_{\Phi(f)}\ex X\ov{\ex B_{0}} \]
is complex linear.
\item $J_{0}$ restricted to the subspace of $\mathbb R$-nil vectors in $T_{f}F(\hat f_{i})\ov{\ex B_{0}}$ is the canonical complex structure.
\item Although $J_{0}$ is only defined at holomorphic curves, it extends to a $\C\infty1$ almost complex structure on $\ex F(\hat f_{i})\ov{\ex B_{0}}$ on a neighborhood of the holomorphic curves
\end{itemize}
Any two embedded  Kuarnishi structures on $\Msw(\hat {\ex B})$  with  complex structures $J_{0}$   constructed using the method from \cite{evc} are joined by a cobordism consisting of an embedded Kuranishi structure on  $\Msw(\hat {\ex B}\times \mathbb R)$ with a complex structure obeying the same conditions as above.
\end{remark}

\

\begin{defn}A complex normally rigid structure on $(\ex F,G)$ is a normally rigid structure along with a $G$-invariant choice of complex structure on the fibers of $\ex N_{\ex F^{\s}}\longrightarrow \totl{\ex F^{\s}}$ so that
\begin{itemize}
\item each fiber of $\ex N_{\ex F^{\s}}$ is isomorphic (via a rigid map) to some $\mathbb C^{n}\times \et mP$ with the standard product complex structure
\item the transition maps from Definition \ref{normally rigid structure}
\[\Phi_{\s'<\s}:\ex N_{\ex F^{\s'}}\rvert_{r(U^{\s})}\longrightarrow \ex N_{\ex F^{\s}}\]
are fiberwise complex when the neighborhood  $U^{\s}$ of $\ex F^{\s}\subset\ex N_{\ex F^{\s}}$ is chosen small enough.
\end{itemize}

A normally rigid map $\ex F\longrightarrow \ex X$ is complex normally rigid if the locally defined  fiberwise rigid maps $\ex N_{\ex F^{\s'}}\longrightarrow \ex N_{\ex X^{\s}}$ from Definition \ref{normally rigid map} are fiberwise complex where defined.
 \end{defn}

\begin{remark}\label{induced J}  The extra information involved in upgrading a normally rigid structure on $\ex F$ to a complex normally rigid structure is a choice of complex structure $J$ on the $I_{x}$-anti-invariant part of $T_{x}\ex F$ for all $x$. If $x$ is contained in $\ex F^{\s}$, the image under $r_{\ex F^{\s}}$ of the tangent space to the fiber of $\ex N_{\ex F^{\s}}$ is equal to the span of the $I_{x}$-anti invariant part of $T_{x}\ex F$ and the $\mathbb R$-nil vectors in $T_{x}\ex F$. As the action of $I_{x}$ preserves the complex structure and both these subspaces, the complex structure on $\ex N_{\ex F^{\s}}$ induces a complex structure on the $I_{x}$-anti-invariant part of $T_{x}\ex F$. (The induced  complex structure on the subspace of $\mathbb R$-nil vectors is just the canonical complex structure on the space of $\mathbb R$-nil vectors, which carries no extra information.)
  
\end{remark}

\

\begin{lemma}\label{cn} Let $J$ be the standard complex structure on $\mathbb C^{n}$. There exists a constant $c_{n}>0$ with the following significance:
Given any two linear complex structures $J_{0}$, and $J_{1}$ on $\mathbb C^{n}$ so that 
\[\abs{J_{i}-J}\leq \epsilon <c_{n}\]
let 
\[W:=(1-J_{1}J_{0})^{-1}(1+J_{1}J_{0})\]
and let 
\[J_{\rho}:=(1+\rho W)^{-1}(1-\rho W)J_{0}\]
Then for $\rho\in[0,1]$, $J_{\rho}$ is also a linear complex structure on $\mathbb C^{n}$ which satisfies
\[\abs{J_{\rho}-J_{0}}\leq \epsilon\]
\end{lemma}

\pf

This is a well known method for averaging almost complex structures. For the space of endomorphisms $X:\mathbb C^{n}\longrightarrow C^{n}$  sufficiently close to $J_{0}$, the invertible  transformation 
 \[X\mapsto (1-XJ_{0})^{-1}(1 +X J_{0}):=W\] changes the condition of $X$ being a complex structure into a linear condition $J_{0}^{-1}WJ_{0}=-W$. $J_{\rho}$ is then the linear path between $J_{1}$ and $J_{0}$ in these new coordinates. Consider   $\abs {X-J}^{2}$ as a function of $X$ in these linear linear coordinates. If $J_{0}=J$,  then for $X$ in some neighborhood of $J$, the second derivative of  $\abs{X-J}^{2}$  is positive definite in these linear coordinates, therefore the same holds for $J_{0}$ in a neighborhood of $J$. It follows that for $J_{1}$ and $J_{0}$ close enough to $J$, the function $\abs {J_{\rho}-J}^{2}$ is a convex function of the parameter $\rho$. 
 
 \stop

Note that the constants in Lemma \ref{cn} are bounded above by $1$ and may be chosen so that $c_{n+1}\leq c_{n}$.

\begin{defn}\label{Xcnrs}A complex normally rigid structure on a family of almost complex  exploded orbifolds  $(\ex X,J)\longrightarrow \ex B_{0}$ is nice if there exists some $J$-invariant metric on $\ex X$ so that the following estimates hold:

Given any point $x\in X^{\s}$, $I_{\s}$ acts on $T_{x}\ex X\ov{\ex B_{0}}$, which splits into a $I_{s}$-invariant subspace and an orthogonal $I_{s}$-anti-invariant subspace, which is equal to the space of $I_{s}$-anti-invariant vectors in $T_{x}\ex X$ because $\ex B_{0}$ has no orbifold structure. $J$ therefore acts on this $I_{s}$-anti-invariant subspace of $T_{x}\ex X$. The complex structure on the fibers of $\ex N_{\ex X^{\s}}\longrightarrow \totl{\ex X^{\s}}$ induces a second complex structure $J'$ on the space of $I_{s}$-anti-invariant vectors in $T_{x}\ex X$. Given any nonzero $I_{s}$-anti-invariant vector  $v$, we require that 
\begin{equation}\label{JJ'}\abs{(J-J')v}<c_{n}\abs v\end{equation}
where $c_{n}$ is the constant from Lemma \ref{cn}, and $2n$ the relative  dimension of $\ex X\longrightarrow \ex B_{0}$ at $x$.

\end{defn}

\

\begin{remark} Lemma \ref{nr construction} implies that a  normally rigid structure on $\ex X$ exists so that $\ex X\longrightarrow \ex B_{0}$ is fiberwise rigid. If $\ex X\longrightarrow \ex B_{0}$ is fiberwise complex, a complex structure on this normally rigid normally rigid structure on $\ex X$ locally exists which obeys inequality (\ref{JJ'}) with $c_{n}$ any positive constant. Lemma \ref{cn} implies that such complex normally rigid structures may be averaged, so we may construct a nice complex normally rigid structure on $\ex X$ which obeys inequality (\ref{JJ'}) with any positive constant $c_{n}$. Moreover, any two such nice complex normally rigid structures are cobordant in the sense that they are restrictions of a nice complex normally rigid structure on $\ex X\times \mathbb R\longrightarrow\ex B_{0}\times \mathbb R$. (The argument for this is given is more detail in the proof of Theorem \ref{rks thm}.)
\end{remark}

\begin{defn}\label{metric def}Given an embedded Kuranishi structure $\{(\mathcal U_{i},V_{i},\hat f_{i}/G_{i})\}$ with compatible extensions $\hat f_{i}^{\sharp}$, and  a complex structure $J_{0}$ on $T_{f}\ex F(\hat f_{i}^{\sharp})\ov{\ex B_{0}}$ satisfying the conditions of Remark \ref{J_{0} conditions}, a compatible system of metrics is
\begin{itemize}
\item
 a choice of  $G_{i}$-invariant metric on $\ex F(\hat f^{\sharp}_{i})$ for all $i$ so that for all $i$ and $j$, there exists a metric on $\ex F(\hat f_{i}^{\sharp}\times_{\Msw}\hat f_{j}^{\sharp})$ so that the maps in the following diagram are isometries:
\[\begin{tikzcd}\ex F(\hat f_{i}^{\sharp}\times_{\Msw}\hat f_{j}^{\sharp})\rar\dar&\ex F(\hat f_{j}^{\sharp})
\\ \ex F(\hat f_{i}^{\sharp})\end{tikzcd}\]
and so that $J_{0}$ acts on $T_{f}\ex F(\hat f_{i}^{\sharp})\ov{\ex B_{0}}$ by isometries. 
\item a choice of $G_{i}$-invariant metric on each $V_{i}(\hat f_{i}^{\sharp})$  which is preserved by the complex structure and so that if $\dim V_{i}\leq \dim V_{j}$, then there is a metric on $V_{i}(\hat f_{i}\times_{\Msw}\hat f_{j})$ so that the natural maps
\[\begin{tikzcd}V_{i}(\hat f_{i}^{\sharp}\times_{\Msw}\hat f_{j}^{\sharp})\rar\dar&V_{j}(\hat f_{j}^{\sharp})
\\ V_{i}(\hat f_{i}^{\sharp})\end{tikzcd}\]
are isometries
\end{itemize}
\end{defn}

Such metrics exist locally, and may be averaged, so such a compatible system of metrics always exists.

\begin{defn}\label{rks def}Given a complex submersion $\Phi:\Msw(\hat {\ex B})\longrightarrow \ex X$ and a complex normally rigid structure on $\ex X$ which is nice in the sense of Definition \ref{Xcnrs}, a rigidified embedded Kuranishi structure compatible with $\Phi$ is:
\begin{itemize}\item an embedded, $\Phi$-submersive Kuranishi structure $(\mathcal U_{i},V_{i},\hat f_{i}/G_{i})$  with compatible extensions $(\mathcal U_{i}^{\sharp},V_{i},\hat f^{\sharp}_{i}/G_{i})$
\item   a complex structure $J_{0}$ on $T_{f}\ex F(\hat f_{i}^{\sharp})\ov{\ex B_{0}}$ at holomorphic curves $f$  in $\hat f_{i}^{\sharp}$ satisfying the conditions in Remark \ref{J_{0} conditions}, 
\item  a compatible system of metrics as in Definition \ref{metric def},
\item a complex normally rigid structure on  $(\ex F(\hat f_{i}^{\sharp}),G_{i})$,
\item a normally rigid vector bundle structure on  $V_{i}(\hat f_{i}^{\sharp})\longrightarrow  \ex F(\hat f_{i}^{\sharp})$
\end{itemize}
so that the following conditions hold:
\begin{enumerate}
\item \label{r3}If $V_{i}\subset V_{j}$ on $\mathcal U^{\sharp}_{i}\cap\mathcal U_{j}^{\sharp}$, then recall that in  the fiber product diagram,  
\[\begin{tikzcd} \ex F(\hat f_{i}^{\sharp}\times_{\Msw}\hat f_{j}^{\sharp})\dar\rar&\ex F(\hat f_{j}^{\sharp})
\\ \ex F(\hat f_{i}^{\sharp})\end{tikzcd}\]
 the downward pointing map is equal to a $G_{i}$-equivariant $G_{j}$ bundle over the open subset $\mathcal U^{\sharp}_{j}\cap \ex F(\hat f_{j})$ of $\ex F(\hat f_{j})$, and that the righthand map  is a $G_{j}$-equivariant, $G_{i}$-fold cover of the sub exploded manifold of $\ex F(\hat f_{j}^{\sharp})\cap \mathcal U^{\sharp}_{i}$ where $\dbar$ is in $V_{i}\subset V_{j}$.
  
\begin{enumerate}\item $ \ex F(\hat f_{i}^{\sharp}\times_{\Msw}\hat f_{j}^{\sharp})$ has a complex normally rigid structure so that the maps in the diagram above
are  complex normally rigid maps,
\item and the inclusion of $V_{i}$ as a subvectorbundle of $V_{j}$ over the righthand map is normally rigid. 
\item\label{rk3} Moreover, if $f$ is any curve in $\hat f_{j}^{\sharp}$ which is isomorphic to a curve in $\hat f_{i}^{\sharp}$,
\[D\pi_{V_{i}}\dbar:T_{f}\ex F(\hat f^{\sharp}_{j})\longrightarrow V_{j}(f)/V_{i}(f)\]
is surjective and is approximately holomorphic when restricted to $I_{f}$-anti-invariant vectors in the following sense:
 Given any $I_{f}$-anti-invariant vector $v$ in $T_{f}\ex F(\hat f_{j}^{\sharp})$, 
\[\abs{J\nabla_{v}\dbar-\nabla_{Jv}\dbar}_{V_{j}/V_{i}}\leq \epsilon \abs{\nabla_{v}\dbar}_{V_{j}/V_{i}}\]
\begin{itemize}\item where $\abs{\cdot}_{V_{j}/V_{i}}$ indicates size of the orthogonal projection of a vector in $V_{j}$ to $V_{j}/V_{i}$ using our compatible system of metrics,
\item $J$ indicates both the induced complex structure on the $I_{f}$-anti-invariant part  of $T_{f}\ex F(\hat f)$ (see Remark \ref{induced J}), and also the complex structure on $V_{j}$,
\item and $\epsilon<\frac{c_{n}}{2n}$, where the constant  $c_{n}$ is from Lemma \ref{cn} and  $n$ is the maximal complex dimension of $V_{k}$ so that $f$ is in $\hat f_{k}^{\sharp}$. \end{itemize}
\end{enumerate}
\item \label{r4}The maps $\ex F(\hat f_{i}^{\sharp})\longrightarrow \ex X$ are complex normally rigid maps
\item \label{r5}  at any holomorphic curve $f$, given any nonzero $I_{f}$-anti-invariant vector $v$ in $T_{f}\ex F(\hat f_{i}^{\sharp})\ov{\ex X}$,
\[\abs{(J-J_{0})v}<c_{n}\abs v\]
where $J$ is the complex structure induced by our complex normally rigid structure as in Remark \ref{induced J},   $c_{n}$ is the constant from Lemma \ref{cn}, and $2n$ is the largest relative dimension at $f$ of the maps $\Phi:\ex F(\hat f_{j}^{\sharp})\longrightarrow \ex X$. 
\end{enumerate}
\end{defn}

To understand condition \ref{rk3} above, note that if $f$ is in $\hat f_{j}^{\sharp}$ and $\mathcal U_{i}^{\sharp}$, then $f$ is isomorphic to a curve in $\hat f_{i}^{\sharp}$ if and only if  $\dbar f\in V_{i}$. If $f$ is holomorphic, then $D\dbar$ is transverse to $V_{i}$ at $f$, so $D\dbar \pi_{V_{i}}:T_{f}\Msw\longrightarrow \Y(f)/V_{i}(f)$ is surjective and has kernel which may be identified with $T_{f}\ex F(\hat f_{i}^{\sharp})$. In particular, this implies that for any $f$ so that $\dbar f\in V_{i}$ and $f$ is close enough to a holomorphic curve,  \[D\pi_{V_{i}}\dbar:T_{f}\ex F(\hat f^{\sharp}_{j})\longrightarrow V_{j}(f)/V_{i}(f)\]
is surjective, and has kernel equal to the image of $T_{f}\ex F(\hat f_{i}^{\sharp})$. Condition \ref{rk3} roughly states that $D\pi_{V_{i}}\dbar$ is approximately holomorphic when restricted to $I_{\s}$-anti-invariant vectors. We shall see that the inequality used to state this in condition \ref{rk3} is preserved by our method of averaging complex structures. 


\begin{thm}\label{rks thm}Given any complex submersion $\Phi:\Msw(\hat {\ex B})\longrightarrow \ex X$  and a complex normally rigid structure on $\ex X$ satisfying the conditions of Definition \ref{Xcnrs}, there exists a rigidified embedded Kuranishi structure on $\Msw(\hat{\ex B})$ compatible with $\Phi$. 

Any two such rigidified embedded Kuranishi structures are cobordant in the sense that given two such structures compatible with possibly different complex normally rigid structures on $\ex X$,  there exists  a complex normally rigid structure on $\ex X\times \mathbb R$ which restricts to the two given structures at $\{0,1\}\subset \mathbb R$, and
there exists a rigidified embedded Kuranishi structure on $\Msw(\hat {\ex B}\times \mathbb R)$ compatible with the complex submersion $\Msw(\hat {\ex B}\times \mathbb R)\longrightarrow \ex X\times \mathbb R$ and this rigidified embedded Kuranishi structure restricts to $\{0,1\}\subset\mathbb R$ to be the two original rigidified embedded Kuranishi structures.\end{thm}

\pf

It is proved in \cite{evc} that there exists a $\Phi$-submersive embedded Kuranishi structure $(\mathcal U_{i},V_{i},\hat f_{i})$ on $\Msw$, and that any two such $\Phi$-submersive embedded Kuranishi structures are cobordant. It is also proved in \cite{evc} that  there exists a complex structure $J_{0}$  satisfying the conditions of Remark \ref{J_{0} conditions}, and any choice of such $J_{0}$ defined over $\{0,1\}$ in our cobordism may be extended to a choice of $J_{0}$ on the entire cobordism. We may then choose a system of metrics satisfying Definition \ref{metric def}.

Lemma \ref{nri} implies that at the expense of shrinking our extensions $\hat f_{i}^{\sharp}$ a little, we may choose compatible normally rigid structures on $\ex F(\hat f^{\sharp}_{i})$, and that we may choose these normally rigid structures so that the maps $\ex F(\hat f^{\sharp}_{i})\longrightarrow \ex X$ are normally rigid. Lemma \ref{nri} also implies that cobordisms between embedded Kuranishi structures lift to cobordisms between such normally rigid structures.  

Lemma \ref{nrvi} implies that at the cost of reducing the size of our extensions $\hat f^{\sharp}_{i}$, we may choose compatible normally rigid structures on the complex vector bundles $V_{i}$ over $\ex F(\hat f_{i}^{\sharp})$. Lemma \ref{nrvi} also implies that any cobordism of embedded Kuranishi structures and normally rigid structures defined so far lifts to a cobordism between the choice of normally rigid structure on these vector bundles $V_{i}$.

It remains to construct complex structures on our already chosen normally rigid structures on $\ex F(\hat f_{i}^{\sharp})$. This step may require us to shrink our Kuranishi charts, but our shrunken Kuranishi charts shall always have the same intersection with the moduli stack of holomorphic curves. To keep notation simple, we shall not use separate notation for these smaller Kuranishi charts.

Choose a well ordering $\prec$ on our Kuranishi charts as in Lemma \ref{well order}, and also choose extensions $\hat f_{i,j}$ and $\hat f_{i,\underline j}$ of $\hat f_{i}$ satisfying the conditions of Lemma \ref{well order}, so for $i\preceq j\prec k$
\[\hat f_{i}\exte \hat f_{i,k}\exte \hat f_{i,\underline k}\exte\hat f_{i,j}\]

\begin{claim}\label{Jextend}Suppose that \begin{itemize}\item$f$ is a holomorphic curve in $\hat f_{j,\underline j}$,
 \item $f$ is a curve in $\hat f_{i,\underline j}$ for some $i\prec j$,
  \item for all $i\prec j$
  so that $f$ is in $\hat f_{i,\underline j}$, a complex structure on the normally rigid structure has been chosen on a $G_{i}$-equivariant neighborhood of $f$ in $\ex F(\hat f_{i,\underline j})$, and that this complex structure satisfies conditions \ref{r4} and \ref{r5}, and satisfies condition $\ref{r3}$ for the maps
\[\ex F(\hat f_{i,i'}\times_{\Msw}\hat f_{i',i'})\longrightarrow \ex F(\hat f_{i',i'})\]
 for all  $i\prec i'\prec j$ so that $f$ is in both $\hat f_{i,i'}$ and $\hat f_{i',i'}$.
 \end{itemize}
  Then there exists a complex structure on the normally rigid structure on a $G$-equivariant neighborhood of $f$ in $\ex F(\hat f_{j,\underline j})$ which satisfies conditions \ref{r4} and \ref{r5} and satisfies condition \ref{r3} for the maps
 \[\ex F(\hat f_{i,j}\times_{\Msw}\hat f_{j,j})\longrightarrow \ex F(\hat f_{j,j})\]
   for all $i\prec j$ so that $f$ is in $\hat f_{i,\underline j}$.\end{claim}

We shall now prove Claim \ref{Jextend}. Choose an $i$ so that $V_{i}$ has the maximal dimension so that $i\prec j$ and $f$ is in $\hat f_{i,\underline j}$. Suppose that $f$ is in the strata $\ex F(\hat f_{j,\underline j})^{\s}$ and $\ex F(\hat f_{i,\underline j}\times_{\Msw}\hat f_{j,\underline j})^{\s'}$. We must choose a complex structure $J$ on the fibers of $\ex N_{\ex F(\hat f_{j,\underline j})^{\s}}$ in a neighborhood of $f$. Let $W$ indicate the $I_{\s}$-anti-invariant part of $V_{i}/V_{j}$ (restrict these vector bundles to the substack represented by  $\ex F(\hat f_{i,\underline j}\times_{\Msw}\hat f_{j,\underline j})^{\s'}/G_{i}\times G_{j}$ so that $I_{\s}$ acts on each fiber). There exists a $G_{i}$-invariant, $G_{j}$-equivariant map of rigid bundles 
\[\begin{tikzcd}W\dar\rar& \ex N_{\ex F(\hat f_{j,j})^{\s}}\dar
\\\totl{\ex F(\hat f_{i,\underline j}\times_{\Msw}\hat f_{j,\underline j})^{\s'}}\rar&\totl{\ex F(\hat f_{j,\underline j})^{\s}}
\end{tikzcd}\]
defined in a neighborhood of $f$ so that 
\begin{itemize} \item the image of $W$ is the subspace of the kernel of $T\Phi:T\ex F(\hat f_{j,\underline j})\rvert_{\ex F(\hat f_{j,\underline j} f)^{\s}}\longrightarrow T\ex X$ which is orthogonal to the image of $T\ex N_{\ex F(\hat f_{i,\underline j}\times_{\Msw}\hat f_{j,\underline j})^{\s'}}\ov{\totl{\ex F(\hat f_{i,\underline j}\times_{\Msw}\hat f_{j,\underline j})^{\s'}}}$.
\item this map composed with the derivative of $\pi_{V_{i}}\dbar:\ex F(\hat f_{j,\underline j})\longrightarrow V_{j}/V_{i}$ is the canonical identification of $W$ as a subspace of $V_{j}/V_{i}$,
\end{itemize}
 (Such a map locally exists around $f$ because part of the definition of Kuranishi charts being $\Phi$-submersive is that $D\dbar$ restricted to the kernel of $T_{f}\Phi$ is transverse to $V_{i}(f)$.) Using the above map, the following is a pullback diagram:

\[\begin{tikzcd}\ex N_{\ex F(\hat f_{i,\underline j}\times_{\Msw}\hat f_{j,\underline j})^{\s'}}\oplus W\dar \rar&\ex N_{\ex F(
\hat f_{j,\underline j})^{\s}}\dar
\\ \totl{\ex F(\hat f_{i,\underline j}\times_{\Msw}\hat f_{j,\underline j})^{\s'}}\rar &\totl{\ex F(\hat f_{j,\underline j})^{\s}}
 \end{tikzcd}\]
The product complex structure on $\ex N_{\ex F(\hat f_{i,\underline j}\times_{\Msw}\hat f_{j,\underline j})^{\s'}}\oplus W$ is $G_{i}\times G_{j}$-invariant, and the top map above is $G_{i}$-invariant and $G_{j}$-equivariant, therefore we may use this product complex structure to define $J$ on the image of the top map above. At $f$, the inclusion of $W$ is automatically $J_{0}$-holomorphic, therefore  any extension of this $J$ will then obey condition \ref{r5} in a neighborhood of $f$. Also note that on a neighborhood of $f$, any extension of $J$ will automatically  obey condition \ref{r3} for the map $\ex F(\hat f_{i,\underline j}\times_{\Msw}\hat f_{j,\underline j})\longrightarrow \ex F(\hat f_{j,\underline j})$. Because $V_{i}$ has maximal dimension of any $i\prec j$ with $f\in\hat f_{i,\underline j}$ condition \ref{r3} will also hold on a neighborhood of $f$ for the maps $\ex F(\hat f_{k,\underline j}\times_{\Msw}\hat f_{j,\underline j})\longrightarrow \ex F(\hat f_{k,\underline j})$ for any $k\prec j$ with $f\in \hat f_{k,\underline j}$

Because $W$ is in the kernel of $\Phi$, our complex structure obeys condition \ref{r4} on the image of $\ex F(\hat f_{i,\underline j}\times_{\Msw}\hat f_{j,\underline j})^{\s'}$. We may choose our extension to also satisfies condition \ref{r4}.

Once this complex structure is defined on a neighborhood of $f$ in $\ex F(\hat f_{j,\underline j})^{\s}$, it may be extended uniquely to a neighborhood of $f$ in $\ex F(\hat f_{j,\underline j})^{\s_{0}}$ for $\s_{0}<\s$ using the condition that the maps $\Phi_{\s_{0}<\s}$ are complex. On a small enough neighborhood,  the analogous maps for $\ex F(\hat f_{i,\underline j})^{\s'}$ are also complex, therefore condition \ref{r3} will still hold in a neighborhood which is chosen small enough that the inequality from condition \ref{rk3} still holds.
Similarly, condition \ref{r5} will hold in a neighborhood, and condition \ref{r4} will automatically be satisfied. 

 This completes the proof of Claim \ref{Jextend}.

\begin{claim}\label{J2}Around any holomorphic curve $f$ in $\hat f_{j}^{\sharp}$, there exists a complex structure on our normally rigid structure which obeys conditions \ref{r4} and \ref{r5}
\end{claim}

 At the holomorphic curve $f$, we may choose a complex structure on $\ex N_{\ex F(\hat f_{j}^{\sharp})^{\s(f)}}$ which obeys condition \ref{r4}, and which agrees with $J_{0}$ on the $I_{f}$-anti-invariant part of the kernel of $T_{f}\Phi$. At $f$, this complex structure obeys condition \ref{r5}, and therefore any extension obeys condition \ref{r5} in a neighborhood. We may choose an extension which obeys condition \ref{r4}. 

This complex structure on $\ex N_{\ex F(\hat f_{j}^{\sharp})^{\s}}$ on a neighborhood of $f$ may be extended for all $\s'<\s$ to a complex structure on $\ex N_{\ex F(\hat f_{j}^{\sharp})^{\s'}}$ on a neighborhood of  $f$ so that the maps $\Phi_{\s'<\s}$ are all complex. The resulting complex normally rigid structure obeys conditions \ref{r4} and \ref{r5}, but does not necessarily  obey condition \ref{r3}. This completes the proof of Claim \ref{J2}.

\

We are now ready to construct a complex  normally rigid structure on a neighborhood of the holomorphic curves within $\hat f_{j,j}$. Recall that our previously constructed charts are on extensions of $\hat f_{i,\underline j}$, and $\hat f_{i,j}\exte\hat f_{i,\underline j}$. For any holomorphic curve $f$ in $\hat f_{j,\underline j}$, which is also in some previously considered chart,  use the complex structure from Claim \ref{Jextend} on a neighborhood which does not intersect $\hat f_{i,j}$ if $f$ is not in $\hat f_{i,\underline j}$.  For other holomorphic curves $f$ use the complex structure from Claim \ref{J2} on a neighborhood of $f$ which does not intersect $\hat f_{i,j}$ for any $i\prec j$. Each of these complex structures obeys all the required conditions on their domains of definition when condition \ref{r3} is only required to hold for $\hat f_{i,j}\exte\hat f_{i,\underline j}$. Using the above domains, we may pick a countable cover of the holomorphic curves in $\hat f_{j,\underline j}$, then piece together the corresponding complex structures  using normally rigid cutoff functions and the averaging method from Lemma \ref{cn}. Lemma \ref{cn} implies that averaging complex structures in this way will produce a complex normally rigid structure which still obeys all the required conditions, but which is now defined on a neighborhood of all the holomorphic curves in $\hat f_{j,j}$.

The construction of a normally rigid structure may then continue by transfinite induction to give a complex normally rigid structure on a neighborhood of the holomorphic curves in $\hat f_{i,i}$ for all $i$ so condition \ref{r3} holds for the maps 
$\hat f_{i,j}\times_{\Msw}\hat f_{j,j}\longrightarrow \hat f_{j,j}$, and all other conditions hold. Then we may restrict our Kuranishi charts and extensions to subsets on which our complex normally rigid structure is defined, defining the required rigidified embedded Kuranishi structure. 

Because the complex structures on our normally rigid structures may be averaged and still satisfy the required conditions, it follows that any two rigidified embedded Kuranishi structures are cobordant.

\stop

\section{Normally complex sections}

Let $U:=\mathbb C^{n}\times \et mP$ be a standard coordinate chart with a rigid, holomorphic action of $G$, and $V$ be a finite dimensional complex vector space with a $\mathbb C$-linear action of $G$.

 Recall from \cite{iec} that a smooth monomial on $\et mP$ is a complex function in the form $\zeta=\totl{c\tilde z^{\alpha}}$. We may pick a basis $\zeta_{1},\dotsc,\zeta_{N}$ of smooth monomials on $\et mP$. Given a complex vector space $V$, say that a map $\et mP\longrightarrow V$ is a  polynomial map if it may be written as a complex polynomial map in the coordinates $\zeta_{1},\dotsc,\zeta_{N}$.    
 We shall need a notion of degree for polynomials on $\et mP$, however there may be monomial relations between the $\zeta_{i}$ which do not respect degree. (For example, maybe $\zeta_{1}\zeta_{2}=\zeta_{3}^{5}$.)
 
 \begin{defn}\label{degree def} Given a point $b$ in the interior of the polytope  $P$,  define the $b$-degree of a smooth monomial $\zeta=\totl{c\e a \tilde z^{\alpha}}$ on $\et mP$ to be $\alpha\cdot b +a$, which may be intrinsically regarded as the tropical part of $c\e a\tilde z^{\alpha}$ evaluated at the point $b\in P$.
 
 Given a monomial on $\mathbb C^{n}\times \et mP$ in the form $z^{\beta}\zeta$ where $\zeta$ is a smooth monomial on $\et mP$, and $z^{\beta}$ is a monomial on $\mathbb C^{n}$, define the $b$-degree of $z^{\beta}\zeta$ to be the usual degree of $z^{\beta}$ plus the $b$-degree of $\zeta$.
   \end{defn}
   
   We shall also call a polynomial  $b$-homogeneous if it is a sum of monomials with some fixed $b$-degree. Any polynomial breaks up canonically into a  sum of $b$-homogeneous  polynomials. We shall see in Lemma \ref{degree part} below that similarly, any holomorphic function has a well defined notion of a $b$-degree $d$ part which is a $b$-homogeneous polynomial of degree $d$. 
   
    In the case that there is a rigid  action of  $G$  on $\mathbb C^{n}\times \et mP$, we shall always choose $b\in P$ to be fixed by $G$. Then the action of $G$ preserves the set of $b$-homogeneous polynomials with  any fixed $b$-degree. In particular, the $b$-degree $d$ part of any $G$-invariant polynomial is still $G$-invariant.

 \begin{remark} The $b$-degree of a product of monomials is the sum of the $b$-degrees of the monomials involved. Because $b$ is in the interior of $P$, the $b$-degree of each nonconstant monomial is strictly positive. It follows that the vector space of polynomials with $b$-degree less than  $d$ is finite dimensional. Note also that there exists some constant $c>0$ so that the polynomials of with $b$-degree less than $d$ contain all the polynomials in $\zeta_{1},\dotsc,\zeta_{N}$ with degree less than $cd$.  \end{remark}

We can interpret the $b$-degree $d$ polynomials as polynomials that scale in a certain way as follows: Embed  $\totl{\et mP}$ inside $\mathbb C^{N}$ using the map $(\zeta_{1},\dotsc,\zeta_{N})$. There exists  a unique smooth vector field $v_{b}$ on $\mathbb C^{N}$ defined by by  requiring that $v_{b}\zeta_{i}$ is $\zeta_{i}$ times the $b$-degree of $\zeta_{i}$. Because the relations between the $\zeta_{i}$ are all $b$-homogeneous, $v_{b}$ is tangent to $\totl{\et mP}$. The defining property of the corresponding vectorfield $v_{b}$ on $\et mP$ is that given any monomial $\zeta$, $v_{b}\zeta$ is equal to $\zeta$ times the $b$-degree of $\zeta$. We may similarly define an analogous vectorfield $v_{b}$ on $\mathbb C^{n}\times \totl{\et mP}$. Define
\[\Phi_{b,t}:\mathbb C^{n}\times \totl{\et mP}\]
to be the flow of $v$ for time $t$.  The map $\Phi_{b,t}$ scales every monomial with $b$-degree $d$ by a factor of $e^{td}$.

\begin{lemma}\label{degree part} Given any holomorphic function $f$ defined on a neighborhood of  $0\times \et m{P^{\circ}}\subset\mathbb C^{n}\times \et mP$, and a $b$-degree $d$ there exists a unique polynomial $p$ on $\mathbb C^{n}\times \et mP$ which approximates $f$ up to $b$-degree $d$ in the following sense:

\[\lim_{t\to -\infty}e^{-t d}(f-p)\circ\Phi_{b,t}=0\]
(where the limit may be interpreted pointwise).
\end{lemma}

 \pf
 
 Note that for all $x\in \mathbb C^{n}\times \totl{\et mP}$, the limit as $t\rightarrow -\infty$ of $\Phi_{b,t}(x)$ is $0\times \totl{\et m{P^{\circ}}}$, so $f$  is defined at $\Phi_{b,t}(x)$ for $t$ negative enough. 
  As $f$ may be written as a holomorphic function of $z_{1},\dotsc,z_{n}, \zeta_{1},\dotsc,\zeta_{N}$, $f$ may be written as a power series which converges absolutely in some neighborhood $U$ of $0$ in $\mathbb C^{n}\times \mathbb C^{N}$. Let $p$ consist of all the terms in this powerseries with $b$ degree less than or equal to $d$. Composition with $\Phi_{b,t}$ multiplies the coefficient of each monomial with $b$-degree $d$  by $e^{td}$. For any $x$, $\Phi_{b,t}(x)$ is contained in $U$ for $t$ negative enough, so the effect of composing $f$ with  $\Phi_{b,t}$  may be read off from this powerseries for $t$ negative enough. Let $d'>d$ be the smallest $b$-degree for which there is some $b$-degree $d'$ monomial. As $t\to -\infty$, it follows that $e^{-t d}(f-p)\circ\Phi_{b,t}(x)$ is bounded by some constant times $e^{(d'-d)t}$, therefore $p$ satisfies the required conditions for being an approximation of $f$ up to $b$-degree $d$.
  
  The uniqueness of $p$ follows from the observation that if  $\lim_{t\to -\infty}e^{-t d}p'\circ\Phi_{b,t}=0$ and $p'$ is a polynomial with $b$-degree  at most $d$, then $p'=0$. 

\stop

 \

Let $\poly {}{}UV$ denote the complex vector space of polynomial maps from $U$ to $V$. We may multiply such maps by polynomials in $\poly{}{}U{\mathbb C}$. Thought of this way, $\poly{}{}UV$ is a finitely generated $\poly{}{}U{\mathbb C}$ module. If $G$ acts on $U$ and $V$, there is an action of $G$ on $\poly{}{}UV$ by conjugation. Let $\poly{}GUV$ denote the polynomials which are invariant under this action--- these are the $G$ equivariant polynomial maps.  As $G$ is finite, $\poly{}GUV$ is a Noetherian $\poly{}GU{\mathbb C}$-module, and $\poly{}GU{\mathbb C}$ is a finitely generated complex algebra. In fact, Proposition 6.8 from \cite{schwarz} implies that given any $G$-equivariant open subset $U'$ of $U$, the $G$-equivariant holomorphic maps $U'\longrightarrow V$ are generated by $\poly{}GUV$ as a module over the ring of $G$-invariant holomorphic functions on $U'$.

Note that if $U=\mathbb C^{n}\times \et mP$, then any complex, rigid isomorphism $U\longrightarrow U$  or complex linear isomorphism $V\longrightarrow V$ pulls back $\poly{}{}UV$ to $\poly{}{}UV$. Therefore, the following definition makes sense:

\begin{defn}\label{polydef} Let $\ex N\longrightarrow M$ be a complex rigid bundle  with a fiber preserving, complex and rigid action of $G$, and let $V$ be a complex vector bundle over $M$ with an action of $G$. Then define $\poly {}G{\ex N}V$ to be the bundle over $M$ with sections $p$ consisting of smooth, fiber preserving maps $\ex N\longrightarrow V$ which restricted to each fiber are $G$-equivariant polynomial maps. 

Define $\poly dG{\ex N}V\subset \poly{}G{\ex N}V$ to be the sub bundle consisting of maps which are fiberwise polynomial maps of $b$-degree less than or equal to $d$.
 
\end{defn}

Note that $\poly{}G{\ex N}V$ is a finitely generated $\poly{}G{\ex N}{\mathbb C}$-module, and that it is generated by $\poly dG{\ex N}V$ for $d$ large enough. Note also that there is a natural evaluation map \[\poly{}G{\ex N}V\times_{M}\ex N\longrightarrow V\] which sends $(p,x)$ to $p(x)$. 

\begin{defn}\label{zdef} Using the notation of Definition \ref{polydef}, define 
\[Z_{d}\subset \poly dG{\ex N}V\times_{M}\ex N\]
to be the inverse image of $\{0\}$ under the natural evaluation map, so 
\[Z_{d}:=\{(p,x)\text{ so that }p(x)=0\}\subset \poly dG{\ex N}V\times_{M}\ex N\]
 
 When no ambiguity is possible, we shall drop the subscript $d$ and use  $Z:=Z_{d}$.
\end{defn}

\

Denote the fiber of $\ex N$ and $V$ over a point  by $U$ and $V$ respectively. The smooth part of fiber of $Z$ over this point is an algebraic subvariety of $\poly dGUV\times \totl{U}$. In particular, this implies that the singularities of $\totl Z$  are modeled on complex algebraic singularities. 

We shall encode the `niceness' of singularities of $\totl{Z}$  using a Whitney stratification. A Whitney stratification  of $\totl{Z}$ is a stratification of $\totl{Z}$ into smooth manifolds with two extra regularity conditions referred to as $(a)$ and $(b)$ regularity. As proved in \cite{trotman},  the weaker $(a)$-regularity condition is equivalent to requiring that every  smooth map to $\poly dGUV\times \totl{U}$ is transverse to the strata of $\totl Z$ on an open subset.  Because the fiber of $\totl{Z}$ over a point is a complex algebraic variety, it admits a Whitney stratification with each strata a variety (this is Theorem 19.2 of \cite{whitney}).

\begin{defn}\label{nws} A nice Whitney stratification of  $\totl{Z}\subset \poly dG{\ex N}V\times_{M}\ex N$ is a $G$-invariant Whitney stratification of $\totl{Z}$ so that
\begin{itemize}
\item  any $G$-equivariant smooth vectorfield $v$ on $\poly dG{\ex N}V\times _{M}\ex N$ is tangent to every strata of $\totl{Z}$ if $v$ is tangent to $Z$ in the sense that the flow of $v$ sends $Z$ to itself,
\item and under the projection map $\totl Z\longrightarrow \totl{\ex N}$, every strata of $\totl Z$ is sent inside a strata of $(\totl {\ex N},G)$ from Definition \ref{strata def}.
\end{itemize}
\end{defn}

\begin{lemma}There exists a nice Whitney stratification of $\totl{Z}$.
\end{lemma}

\pf

Recall that $\totl{Z}$ is a bundle of complex algebraic varieties over a smooth manifold $M$. We shall consider the case that $M$ is a point, because the general case follows easily. 

The construction starting on page 210 of \cite{mather} gives  a canonical minimal Whitney stratification of any real semi-algebraic subvariety.  This minimal Whitney stratification is preserved by any diffeomorphism of the ambient space which preserves the subvariety. In our case,  $\totl Z$ is a complex algebraic subvariety, so it may also be considered to be a real algebraic subvariety by forgetting some structure.  The corresponding minimal Whitney stratification of $\totl Z$ may not be one of our nice stratifications simply because some strata may not project to be contained within strata of $(\totl{\ex N},G)$. We shall get around this by artificially introducing singularities where $\totl Z$ intersects strata of $(\totl{\ex N},G)$. 

Denote by $ Z^{\s}$ the intersection of $\totl Z$ with $\poly dG{\ex N}V\times \totl{\ex N}^{\s}$. Recall that the closure of $\ex N^{\s}$ consists of the union $\ex N^{\s'}$ for all $\s'\geq \s$. The closure of $\totl{\ex N}^{\s}$ is a complex variety, so   $\bigcup _{\s'\geq \s}Z^{\s'}$ is also a complex subvariety. Choose some finite dimensional vector space $X$ with independent subspaces $E_{0}$, $E_{1}$, and $E_{\s}$ for every nonempty strata $\ex N^{\s}$ of $(\ex N,G)$. Then define the subvariety $Z'\subset X\times \poly dG{\ex N}V\times \totl{\ex N}$ by
\[Z':=\lrb{(E_{0}\cup E_{1})\times \totl{Z}}\bigcup_{\s}\lrb{\oplus_{\s'\geq \s}E_{\s'}}\times Z^{\s}\]

Consider the minimal Whitney stratification of $Z'$. Every strata of $Z'$ which intersects $(E_{0}\setminus 0)\times \totl Z$ must be contained in $E_{0} \times \totl Z$. It follows that $E_{0}\times \totl Z$ must be a union of strata. Similarly, $E_{1}\times \totl Z$ must be a union of strata, so their intersection $\{0\}\times \totl Z$ is a union of strata. The corresponding  Whitney stratification of $\totl Z=\{0\}\times \totl{Z}$ shall be our nice Whitney stratification of $\totl Z$. 

Any smooth $G$-equivariant vectorfield which is tangent to $Z$ is also tangent to all $Z^{\s}$. Such a vectorfield lifts to  a vectorfield on $X\times \poly dG{\ex N}V\times \totl{\ex N}$ which is tangent to $Z'$ and also in the kernel of the projection to $X$. Because mimimal Whitney stratifications are invariant under diffeomorphisms of the ambient space, this lifted vectorfield must be tangent to all strata of $Z'$. It follows that our original $G$-equivariant vectorfield must be tangent to all strata of $\totl Z$. Similarly, our Whitney stratification of $\totl Z$ is $G$-invariant because the action of $G$ on $Z'$ preserves the minimal Whitney stratification.

It remains to check that every strata of $\totl Z$ is contained in $Z^{\s}$ for some $\s$. Any strata of $Z'$ which intersects $\oplus_{\s'\geq\s} E_{\s'}\setminus \oplus _{\s''>\s} E_{\s''}$ must be contained entirely within $\oplus_{\s'\geq \s}E_{\s'}\times \bigcup_{\s'\geq \s}Z^{\s'}$. It follows that $\oplus_{\s'\geq \s}E_{\s'}\times \bigcup_{\s'\geq \s}Z^{\s'}$ is a union of strata. As $\{0\}\times \totl Z$ is also a union of strata, it follows that $\{0\}\times \bigcup_{\s'\geq \s}Z^{\s'}$ is a union of strata for all $\s$, and therefore $Z^{\s}$ is a union of strata for all $\s$.

\stop

\

We shall consider a map from a smooth manifold to be transverse to $Z$ if it is transverse to every strata of some nice Whitney stratification of $Z$. As $\totl {\ex N}$ may not be smooth, we need to make a more careful strata-by-strata definition of a section which is transverse to $Z$.

\begin{defn}\label{Z transverse def} Say that a section $r:\ex N\longrightarrow \poly dG{\ex N}V\times _{M}\ex N$ is transverse to $Z$ at a point $x\in \ex N$ if in a neighborhood of $x$,  there exists a nice  Whitney stratification of $\totl{Z}$ so that the following holds

Let  the strata of $\totl Z$ containing $x$ be $\totl Z_{x}$. To be transverse to $Z$,  $\totl r$ restricted to $\totl{\ex N^{\s(x)}}$ must be transverse to $\totl Z_{x}$ as a submanifold of $\poly dG{\ex N}V\times _{M}\totl{\ex N^{\s(x)}}$. 
\end{defn}

For the above definition to make sense, it was important to require that the strata of $Z$ project to be contained inside the strata of $\ex N$. 

 \begin{lemma}\label{open Z transverse}For any smooth section $r:\ex N\longrightarrow \poly dG{\ex N}V\times _{M}\ex N$, the subset of $\ex N$ on which $r$ is transverse to $Z$ is open. 
\end{lemma}

\pf

Factorize the submersion $\totl {\ex N}\longrightarrow M$ through an embedding into a smooth manifold, $\totl{\ex N}\longrightarrow X\longrightarrow M$. The smooth section $r$ extends to a smooth section $r':X\longrightarrow \poly Gd{\ex N}V\times_{M}X$.

Suppose that $r$ is transverse to $Z$ at $x\in\ex N$. Denote  the strata of $\totl Z$ containing $\totl{r(x)}$ by $\totl Z_{x}$. We have that $r$ restricted to $\totl{\ex N^{\s(x)}}$ is transverse at $x$ to $\totl Z_{x}$ as a submanifold of $\poly dG{\ex N}V\times _{M}\totl{\ex N^{\s(x)}}$. It follows that $r'$ is transverse to $\totl Z_{x}$ at the image of $x$ in $X$. As being transverse to a Whitney stratified space is an open condition, $r'$ is transverse to every strata of $\totl Z$ on a neighborhood of $x$. The fact that every strata of $\totl Z$ is sent inside a strata of $\totl{\ex N}$ then implies that $r$ is transverse to $Z$ on a neighborhood of $x$.

\stop

As the set of ($G$-invariant) sections  transverse to $Z$ at $x$ is dense in the space of ($G$-invariant) sections,  Lemma \ref{open Z transverse} implies that the set of smooth ($G$-invariant) sections transverse to $Z$ is open and dense.

\begin{lemma}\label{ZT1}Suppose that $s_{1}$ and $s_{2}$ are two smooth maps $\ex N\longrightarrow \poly dG{\ex N}V$, and that the graph of $s_{1}-s_{2}$ is in $Z$. Then the graph of $s_{2}$ is transverse to $Z$ if the graph of $s_{1}$ is transverse to $Z$.
\end{lemma}

\pf Consider the smooth automorphism $\Psi$ of  $\poly dG{\ex N}V\times_{M}\ex N$ defined by
\[\Psi(p,x):=(p+s_{1}(x)-s_{2}(x),x)\] This automorphism $\Psi$ preserves any nice Whitney stratification of $\totl Z$ because it is the flow of some vectorfield which may be written as a sum of smooth functions times $G$-equivariant vectorfields tangent to $Z$. The graph of $s_{2}$ composed with $\Psi$ is equal to the graph of $s_{1}$, therefore the graph of $s_{2}$ is transverse $Z$ if and only if the graph of $s_{1}$ is transverse to $Z$.

\stop

\begin{lemma}\label{ZT2} Suppose that $d$ is large enough that $\poly dG{\ex N}V$ generates $\poly {}G{\ex N}V$, and let $s$ be a section $\ex N\longrightarrow \poly dG{\ex N}V$. Consider $s$ as a section $s'$  of $\poly {d+1}G{\ex N}V$. Then the graph of  $s'$ is transverse to $Z_{d+1}$ if and only if the graph of  $s$ is transverse to $Z_{d}$ 
\end{lemma}

\pf

We have chosen $d$ large enough that $\poly dG{\ex N}V$ generates $\poly {d+1}G{\ex N}V$ as a $\poly {}G{\ex N}{\mathbb C}$-module, so given any $p_{0}\in\poly{d+1}G{\ex N}V$, there exist $h_{i}\in\poly {}G{\ex N}{\mathbb C}$ and $p_{i}\in \poly dG{\ex N}V$ so that 
\[\sum_{i}h_{i}p_{i}-p_{0}\]
is a section of $Z_{d+1}$.
In particular, this implies that the family of  transformations 
\[\Psi_{t}(p,x):=(p+t(\sum_{i}h_{i}(x)p_{i}-p_{0}),x)\]
preserves every strata of any nice Whitney stratification of $\totl{Z_{d+1}}$.  Therefore every strata of $\totl{Z_{d+1}}$ is transverse to $\poly dG{\ex N}V\times_{M}\ex N\subset \poly {d+1}G{\ex N}V\times_{M}\ex N$. The intersection of a nice Whitney stratification of $\totl {Z_{d+1}}$ with $\poly dG{\ex N}V\times_{M}\ex N$ therefore defines a Whitney stratification of $\totl {Z_{d}}$. Any $G$-equivariant vectorfield tangent to $Z_{d}$ extends to a $G$-equivariant vectorfield tangent to $Z_{d+1}$, therefore the resulting stratification of $\totl{Z_{d}}$ is also tangent to all vectorfields tangent to $Z_{d}$, and is therefore a nice Whitney stratification. It follows that if $s'$ is transverse to $Z_{d+1}$, then $s$ is transverse to $Z_{d}$.

The fact that $\poly dG{\ex N}V$ generates $\poly {d+1}G{\ex N}V$ implies  that there exists a fiber preserving projection 
\[\Phi:\poly {d+1}G{\ex N}V\times_{M}\ex N \longrightarrow \poly dG{\ex N}V\times_{M}\ex N\]
which is an algebraic submersion restricted to each fiber, and which is the identity restricted to $\poly dG{\ex N}V\times_{M}\ex N$, and which is compatible with the natural evaluation maps to $V$, so $\Phi^{-1}(Z_{d})=Z_{d+1}$. It follows that $\Phi^{-1}$ pulls back any nice Whitney stratification of $\totl{Z_{d}}$ to a nice Whitney stratification of $\totl{Z_{d+1}}$. 
 Therefore if $s$ is transverse to a nice Whitney stratification of $\totl{Z_{d}}$, then $s'$ is transverse to the pulled back nice Whitney stratification of $\totl{Z_{d+1}}$.

\stop

\

\

Recall  Definition \ref{nrvb} of a normally rigid vectorbundle $V$ over $\ex F$. We shall use the notation $\poly{}{I_{s}}{\ex N_{\ex F^{\s}}}{V}$ to mean  Definition \ref{polydef} applied with the normally rigid bundle $\ex N_{\ex F^{\s}}$ over $\totl{\ex F^{\s}}$ and the vector bundle $V$ restricted to $\totl{\ex F}$. There is a natural $G$-equivariant evaluation map
  \[\poly{}{I_{s}}{\ex N_{\ex F^{\s}}}V\times_{\totl{\ex F^{\s}}}\ex N_{\ex F^{\s}}\longrightarrow V_{\ex F^{\s}}\]
where as in Definition \ref{nrvb}, the notation $V_{\ex F^{\s}}$ indicates the pullback of $V$ along $\ex N_{\ex F^{\s}}\longrightarrow \totl{\ex F^{\s}}$.

\

\begin{defn}[normally complex section] \label{nc section def}
 Suppose that $(\ex F,G)$ has a complex normally rigid structure and a normally rigid, $G$-equivariant vectorbundle $V$. A $G$-equivariant section  \[s:\ex F\longrightarrow V\] is normally complex if the following holds:
 
For each strata $\ex F^{\s}$,  there exists some fiber preserving, smooth, $G$-equivariant map
 \[s^{\s}:\ex N_{\ex F^{\s}}\longrightarrow \poly{d}{I_{\s}}{\ex N_{\ex F^{\s}}}{V}\]
so that by restricting to an open neighborhood $U^{\s}$ of $\ex F^{\s}$ in $\ex N_{\ex F^{\s}}$, the following diagram commutes:
\begin{equation}\label{ss diag}\begin{tikzcd}
U^{\s}\rar{(s^{\s},\id)}\dar{r_{\ex F^{\s}}}&\poly{d}{I_{s}}{\ex N_{\ex F^{\s}}}{V}\times_{\totl{\ex F^{\s}}}\ex N_{\ex F^{\s}}\rar&V_{\ex F^{\s}}\dar{l_{\ex F^{\s}}}
\\ \ex F\ar{rr}{s}&& V
\end{tikzcd} \end{equation}

\

Say that $s$ is transverse to $0$ at a point $f\in \ex F^{\s}$ if the graph of  $s^{\s}$ is transverse to $Z$ at $f$ in the sense of Definition \ref{Z transverse def}.
\end{defn}

As noted in Remark \ref{smooth remark}, a choice of normally rigid structure on $\ex F$ defines a smooth structure on $\ex F$ even if $\ex F$ only started off with a $\C\infty1$ structure. In the above definition $\C\infty1$ could be used in place of `smooth'.

\begin{remark}
Note that     $s^{\s}$ corresponds to the zero section $s$ if and only if the graph of $s^{\s}$ is in $Z\subset \poly{d}{I_{s}}{\ex N_{\ex F^{\s}}}{V}\times_{\totl{\ex F^{\s}}}\ex N_{\ex F^{\s}}$. Lemmas \ref{ZT1} and \ref{ZT2} then imply that the definition of a normally complex section being transverse to zero is independent of the choice of section $s^{\s}$ and the choice of $d$. 
\end{remark}

 \begin{lemma}\label{psi}Suppose that $s'<s$, and consider the following commutative diagram from Definition \ref{nrvb}
  \[\begin{tikzcd}
V_{\ex F^{\s'}}\rvert_{U^{\s'}\cap \Phi^{-1}_{\s'<\s}(U^{\s})} \dar[hook]\ar[bend right=50]{ddd} \ar{drr}{l_{\ex F^{\s'}}}
\\V_{\ex F^{\s'}}\rar{\phi_{\s'<\s}}\rvert_{r(U^{\s})}\dar  &V_{\ex F^{\s}}\rar[swap]{l_{\ex F^{\s}}} \dar & V\dar
\\ \ex N_{\ex F^{\s'}}\rar{\Phi_{\s'<\s}}\rvert_{r(U^{\s})}& \ex N_{\ex F^{\s}}\rar{r_{\ex F^{\s}}}& \ex F
\\ U^{\s'}\cap \Phi^{-1}_{\s'<\s}(U^{\s})\uar[hook]\ar{rru}{r_{\ex F^{\s'}}}
\end{tikzcd}\]
There is an open neighborhood $O\subset \ex N_{\ex F^{\s'}}$ of $\ex F^{s'}\cap r(U^{\s})$ so that there exist $G$-equivariant maps of complex vector bundles over   $O$ 
\[ O\times_{\totl{\ex F^{\s}}}\poly d{I_{\s}}{\ex N_{\ex F^{\s}}}V\xrightarrow{\ \ \ \ \psi\ \ \ \ } O\times_{\totl{\ex F^{\s'}}}\poly{d'}{I_{\s'}}{\ex N_{\ex F^{\s'}}}V \]
\[ O\times_{\totl{\ex F^{\s}}}\poly d{I_{\s}}{\ex N_{\ex F^{\s}}}V\xleftarrow{\ \ \ \ \psi'\ \ \ \ } O\times_{\totl{\ex F^{\s'}}}\poly{d'}{I_{\s'}}{\ex N_{\ex F^{\s'}}}V \]
so that the following diagrams commute
\begin{equation}\label{psicd}\begin{tikzcd}
O\ar{dd}{\Phi_{\s'<\s}}    &
O\times_{\totl{\ex F^{\s'}}} \poly {d'}{I_{\s'}}{\ex N_{\ex F^{\s'}}}V\rar\lar&V_{\ex F^{\s'}}\rvert_{r(U^{\s})} \ar{dd}{\phi_{\s'<\s}}
  \\ &\uar{\psi} O\times_{\totl{\ex F^{\s}}}\poly d{I_{\s}}{\ex N_{\ex F^{\s}}}V\dar\ar{ul}
\\ \ex N_{\ex F^{\s}}& \ex N_{\ex F^{\s}}\times_{\totl{\ex F^{\s}}}\poly{}{I_{\s}}{\ex N_{\ex F^{\s}}}V \lar\rar    &V_{\ex F^{\s}}
\end{tikzcd}\end{equation}
\begin{equation}\label{psi'cd}\begin{tikzcd}
O\ar{dd}{\Phi_{\s'<\s}}    &
O\times_{\totl{\ex F^{\s'}}} \poly {d'}{I_{\s'}}{\ex N_{\ex F^{\s'}}}V\rar\lar\dar{\psi'}&V_{\ex F^{\s'}}\rvert_{r(U^{\s})} \ar{dd}{\phi_{\s'<\s}}
  \\ & O\times_{\totl{\ex F^{\s}}}\poly d{I_{\s}}{\ex N_{\ex F^{\s}}}V\dar\ar{ul}
\\ \ex N_{\ex F^{\s}}& \ex N_{\ex F^{\s}}\times_{\totl{\ex F^{\s}}}\poly{}{I_{\s}}{\ex N_{\ex F^{\s}}}V \lar\rar    &V_{\ex F^{\s}}
\end{tikzcd}\end{equation}

 \end{lemma}

\pf Note that all the arrows in the right hand side of diagrams (\ref{psicd}) and (\ref{psi'cd}) may be interpreted as maps of complex vector bundles. It follows that maps $\psi$ or $\psi'$ of complex vector bundles so that diagram (\ref{psicd}) or (\ref{psi'cd}) commute may be averaged, so we may construct $\psi$ or $\psi'$ locally. 

 Locally choose  a basis of  $G$-equivariant sections   $p_{j}$ of $\poly d{I_{\s}}{\ex N_{\ex F^{\s}}}V$. Each of these sections  pull back using $\phi_{\s'<\s}$ to a  $G$-equivariant, fiberwise holomorphic  map 
\[p_{j}':O\longrightarrow V_{\ex F^{\s'}}\rvert_{r(U^{\s})}\]

Proposition 6.8 from \cite{schwarz} implies that $p_{j}'$ may be written as a finite sum of $G$-invariant, fiberwise holomorphic functions on $O$ times sections of $\poly{d'}{I_{\s'}}{\ex N_{\ex F^{\s'}}}{V}$. If $q_{i}$ is locally a basis of $G$-invariant sections of $\poly{d'}{I_{\s'}}{\ex N_{\ex F^{\s'}}}{V}$, then we may write
\[p_{j}'=\sum_{i}h_{i,j}q_{i}\]
We may locally define a map $\psi$ of complex vector bundles using the coordinates given by $p_{j}$ and $q_{i}$. The map with matrix entries $h_{i,j}$ satisfies the required commutativity condition. As locally defined maps satisfying the required commutativity condition may be averaged and patched together, there exists a smooth $G$-invariant map $\psi$
satisfying the required commutativity condition so long as $O$ is contained in  the open subset $U^{\s'}\cap \Phi^{-1}_{\s'<\s}(U^{\s})$ where the commutativity conditions on $\phi_{\s'<\s}$ are required to hold. Note for later use that $\psi$ sends fiberwise holomorphic sections to fiberwise holomorphic sections.

\

We shall now locally construct $\psi'$ around a point $x$ in $\ex F^{\s}\cap r(U^{\s})$.  In what follows, we shall use $Z_{d'}$ to refer to the subset of $\ex N_{\ex F^{\s'}}\times_{\totl{\ex F^{\s'}}}\poly{d'}{I_{s'}}{\ex N_{\ex F^{\s'}}}V$ constructed as in Definition \ref{zdef}.

\begin{claim}\label{kc}Given any section $s$ of $\ex N_{\ex F^{\s'}}\times_{\totl{\ex F^{\s'}}}\poly{d'}{I_{\s'}}{\ex N_{\ex F^{\s'}}}V$, there exists a section $s'$ of $ \ex N_{\ex F^{\s'}}\times_{\totl{\ex F^{\s}}}\poly{d}{I_{\s}}{\ex N_{\ex F^{\s}}}V$ and a section $s_{0}$ of $Z_{d'}$ so that   $s=\psi(s')+s_{0}$ on a neighborhood of $x$.
\end{claim}

The proof of Claim \ref{kc} shall take some time, but first note that the existence of $\psi'$ follows easily from Claim \ref{kc}. In particular, we may choose locally a basis of the vector bundle $\ex N_{\ex F^{\s'}}\times_{\totl{\ex F^{\s'}}}\poly{d'}{I_{\s'}}{\ex N_{\ex F^{\s'}}}V$ consisting of sections $s$, and define $\psi'$ locally to be the linear map  which sends each $s$ to $s'$. Such a $\psi'$ will make diagram \ref{psi'cd} commute because diagram $\ref{psicd}$ commutes. As with the construction of $\psi$, we may then patch together and average such locally defined maps to create a $G$-invariant,  $\psi'$ defined on a neighborhood $O$ of $\ex F^{\s'}\cap r(U^{\s})$ which satisfies the required commutativity condition. 

\

To prove Claim \ref{kc}, we shall use that $\poly {}{I_{\s'}}{\ex N_{\ex F^{\s'}}}V$ is finitely generated over $\poly {}{I_{\s'}}{\ex N_{\ex F^{\s'}}}{\mathbb C}$. Let $h_{1},\dotsc,h_{n}$ be fiberwise holomorphic generators around $x$.  As the action of $I_{\s'}$ preserves $b$-degrees, it follows that we may choose  the $h_{i}$ to be  $b$-homogeneous around $x$.  We may also assume that that the $b$-degree of $h_{i+1}$ is at least the $b$-degree of $h_{i}$. Assume further that all the $h_{i}$ of a given $b$-degree are $\mathbb C$-linearly independent at $x$, and that $h_{i}$ is not in the module generated by all the other $h_{i'}$ when restricted to any neighborhood of $x$. Note that $h_{i}$ do not necessarily generate the vector bundle $\ex N_{\ex F^{\s'}}\times_{\totl{\ex F^{\s'}}}\poly{d'}{I_{\s'}}{\ex N_{\ex F^{\s'}}}V$, but $\{h_{i}\}$ does generate this vector bundle mod $Z_{d'}$ in the sense that any section of this vector bundle may be written as a section of $Z_{d'}$ plus a section of the sub-bundle generated by  $\{h_{i}\}$. 
 
 \begin{claim}\label{approx} Given any   section $h$ in $\poly {}{I_{\s'}}{\ex N_{\ex F^{\s'}}}V$,  there exists a   section $h'$ of $\poly {}{I_{\s}}{\ex N_{\ex F^{\s}}}V$ so that the following diagram commutes up to $b$-degree $d'$ at the point $x$
 \[\begin{tikzcd}V_{\ex F^{\s'}}\rar{\phi_{\s'<\s}}&V_{\ex F^{\s}}
 \\ \ex N_{\ex F^{\s'}}\rar{\Phi_{\s'<\s}}\uar{h}&\ex N_{\ex F^{\s}}\uar{h'}\end{tikzcd}\]
 More precisely, $\phi_{\s'<\s}\circ h$ and $h'\circ\Phi_{\s'<\s}$ restricted to the fiber containing $x$ are holomorphic maps which Lemma \ref{degree part} allows to be approximated at $x$ up to $b$-degree $d'$. These two approximations are equal. 
  \end{claim}
 
Claim \ref{approx} may be easily proved after noting that polynomials maps may be arbitrarily specified up to a given degree  at a finite set of points, and that invariant polynomial maps may be obtained by averaging.  
 
 \
 
In the proofs that follow, we shall need a notion of the $b$-degree $d$ approximation of a fiberwise holomorphic section of $\ex N_{\ex F^{\s'}}\times_{\totl{\ex F^{\s'}}}\poly{d'}{I_{\s'}}{\ex N_{\ex F^{\s'}}}V$. Suppose that
\[s=\sum f_{i}p_{i}\]
 is a fiberwise holomorphic section of $\ex N_{\ex F^{\s'}}\times_{\totl{\ex F^{\s'}}}\poly{d'}{I_{\s'}}{\ex N_{\ex F^{\s'}}}V$, where each $p_{i}$ is a degree $d_{i}$ section of $\poly{d'}{I_{\s'}}{\ex N_{\ex F^{\s'}}}V$, and each $f_{i}$ is a fiberwise holomorphic function on $\ex N_{\ex F^{\s'}}$. The fibewise $b$-degree $d$ approximation of $s$ is the section
 \[\sum f_{i}'p_{i}\]
where on each fiber, $f_{i}'$ is the $b$-degree $(d-d_{i})$ approximation of the fiberwise holomorphic function $f_{i}$. Say that $s$ and $s'$ are equal at $x$ up to $b$-degree $d$ if their fiberwise $b$-degree $d$ approximations agree on the fiber containing $x$. 

\begin{claim}\label{approx2} Given any fiberwise holomorphic  section $s:\ex N_{\ex F^{\s'}}\longrightarrow\ex N_{\ex F^{\s'}}\times_{\totl{\ex F^{\s'}}}\poly{d'}{I_{\s'}}{\ex N_{\ex F^{\s'}}}V$ defined in a neighborhood of $s$, there exists a fiberwise holomorphic  section $s'$ of $ \ex N_{\ex F^{\s'}}\times_{\totl{\ex F^{\s}}}\poly{d}{I_{\s}}{\ex N_{\ex F^{\s}}}V$ and a fiberwise holomorphic section $s_{0}$ of $Z_{d'}$ defined in a neighborhood of $x$ so that    $s$ is equal to $\psi(s')+s_{0}$ up to $b$-degree $d'$ at $x$.

\end{claim}

 
 We shall prove Claim \ref{approx2} from Claim \ref{approx}. Around the point $x$, the section $s$ is locally equal to $\sum_{i} f_{i}p_{i}$ where $f_{i}$ are fiberwise holomorphic $\mathbb C$-valued functions on  $\ex N_{\ex F^{\s'}}$ and $p_{i}$ are $b$-homogeneous sections of $\poly{d'}{I_{\s'}}{\ex N_{\ex F^{\s'}}}V$. Claim \ref{approx} allows us to choose sections $p_{i}'$ of $\poly {}{I_{\s}}{\ex N_{\ex F^{\s}}}V$  equal to $p_{i}$ at $x$ to order $d'$ in the sense of Claim \ref{approx}. Let $s'$ be the section of $\ex N_{\ex F^{\s'}}\times_{\totl{\ex F^{\s}}}\poly{d}{I_{\s}}{\ex N_{\ex F^{\s}}}V$ defined by  $s':=\sum_{i} f_{i}p_{i}'$. Note that $\psi(s')$ is a fiberwise holomorphic section, so around $x$, $\psi(s')$ defines a holomorphic map from the fiber containing $x$ to $V$.  Lemma \ref{degree part} allows us to approximate both $s$ and $\psi(s')$ considered as holomorphic maps to $V$ from this fiber containing $x$. These two approximations are equal up to $b$-degree $d'$ at $x$, because of Claim \ref{approx} and the commutativity of diagram (\ref{psicd}).

 Restricted to the fiber of $\ex N_{\ex F^{\s'}}\longrightarrow \totl{\ex F^{\s'}}$ containing $x$, the $b$-degree $d'$ part of the fiberwise holomorphic section $s-\psi(s')$  is therefore a holomorphic section of $Z_{d'}$, which may be extended to a fiberwise holomorphic section $s_{0}$ of $Z_{d'}$. Claim \ref{approx2} has now been proven because $s$ is equal to $\psi(s')+s_{0}$ up to $b$-degree $d'$  at $x$.

 \

 Try to approximate $h_{1}$ by $\psi(s')$. Claim \ref{approx2} and the observation that $\{h_{i}\}$ generates the vector bundle $\ex N_{\ex F^{\s'}}\times_{\totl{\ex F^{\s'}}}\poly{d'}{I_{\s'}}{\ex N_{\ex F^{\s'}}}V$ mod $Z_{d'}$ implies that there exists a fiberwise holomorphic  section $s'$ of $\ex N_{\ex F^{\s'}}\times_{\totl{\ex F^{\s}}}\poly{d}{I_{\s}}{\ex N_{\ex F^{\s}}}V$ and a section $s_{0}$ of $Z_{d'}$, both defined on a neighborhood of $x$ so that  
 \[\psi(s')+s_{0}=f_{1}h_{1}+ f_{2}h_{2}+\dotsc \]
where $f_{1}(x)=1$, and all other terms $f_{i}h_{i}$ vanish up to the degree of $h_{1}$ at $x$. We may multiply $s'$ and $s_{0}$ by $f_{1}^{-1}$ on a neighborhood of $x$ to obtain new fiberwise holomorphic  sections $s_{1}'$ and $s_{1,0}$ defined on a neighborhood of $x$ so that 
\[\psi(s_{1}')+s_{1,0}=h_{1}+f'_{2}h_{2}+\dotsb \] 
 on a neighborhood of $x$, where  each term $f'_{i}h_{i}$ vanishes at $x$ up to the degree of $h_{1}$ on the fiber of $\ex N_{\ex F^{\s'}}\longrightarrow \totl{\ex F^{\s'}}$ containing $x$.
 
 Now suppose that for all $i<k$ we have on some neighborhood of $x$ a fiberwise holomorphic section  $s'_{i}$ of $\ex N_{\ex F^{\s'}}\times_{\totl{\ex F^{\s}}}\poly{d}{I_{\s}}{\ex N_{\ex F^{\s}}}V$ and a fiberwise holomorphic section $s_{i,0}$ of $Z_{d'}$ so that
 \[\psi(s_{i}')+s_{i,0}=h_{i}+f_{i+1}h_{i+1}+\dotsb \]
 on a neighborhood of $x$.
 Use Claim \ref{approx2} to find fiberwise holomorphic  sections $s'$ and $s_{0}$ defined on a neighborhood of x so that $\psi(s')+s_{0}$ approximates $h_{k}$ at $x$ up to the $b$-degree of $h_{k}$ and so that 
 \[\psi(s')+s_{0}=\sum f_{i}h_{i}\] 
The assumption that $h_{k}$ is not in the module generated by all the other $h_{i}$ implies that  $\sum_{i<k}f_{i}h_{i}$ vanishes at $x$ up to the $b$-degree of $h_{k}$, because otherwise we could arrive at a contradiction by replacing each term with its  approximation at $x$ up to the $b$-degree of $h_{k}$.
   By subtracting off  fiberwise holomorphic functions times the already constructed $s_{i}'$ and $s_{i,0}$, we may therefore obtain new fiberwise holomorphic sections also approximating $h_{k}$ to leading order in the sense of Claim \ref{approx2}, but so that  
 \[\psi(s_{i}')+s_{i,0}=f_{k}h_{k}+f_{k+1}h_{k+1}+
 \dotsb\]
 on a neighborhood of $x$,
 where $f_{k}(x)=1$ and all other terms vanish up to the $b$-degree of $h_{k}$ at $x$. We may obtain sections  $s_{k}'$ and $s_{k,0}$ with the required properties by multiplying by  $f_{k}^{-1}$ on a neighborhood of $x$.
 
 Continuing in this fashion, we obtain sections $s_{n}'$ and $s_{n,0}$ so that 
 \[\psi(s'_{n})+s_{n,0}=h_{n}\] 
 in a neighborhood of $x$, then may obtain the other $h_{i}$ as $\psi(s')+s_{0}$ in a neighborhood of $x$ by starting with $s_{i}'$ and $s_{i,0}$, then subtracting off multiples of $s_{k}'$ and $s_{k,0}$ for $k>i$. As $\{h_{i}\}$ generates the vector bundle $\ex N_{\ex F^{\s'}}\times_{\totl{\ex F^{\s'}}}\poly{d'}{I_{\s'}}{\ex N_{\ex F^{\s'}}}V$ mod $Z_{d'}$, it follows that on a neighborhood of $x$, the image of $\psi$ also  generates the vector bundle $\ex N_{\ex F^{\s'}}\times_{\totl{\ex F^{\s'}}}\poly{d'}{I_{x}}{\ex N_{\ex F^{\s'}}}V$ mod $Z_{d'}$, which is Claim \ref{kc}.
 
 As noted above, now that we have proved Claim \ref{kc}, Lemma \ref{psi} has been proved.

 \stop
%
%
%
%

\

\begin{lemma} The subset of $\ex F$  on which a normally complex section $s:\ex F\longrightarrow V$   is transverse to $0$ is open.
\end{lemma}

\pf

Suppose that $s$ is transverse to $0$ at $f\in \ex F^{\s}$. By definition,  the section $s^{\s}$ of $\ex N_{\ex F^{\s}}\times_{\totl{\ex F^{\s}}}\poly d{I_{\s}}{\ex N_{\ex F^{\s}}}V$ is transverse to $Z_{d}$ at $f$. As the property of being transversal to a Whitney stratification is open, there exists a neighborhood $O$ of $f$ so that $s$ is  transverse to $0$ on $O\cap \ex F^{\s}$. We shall now verify that if $O$ is chosen small enough, then for all $\s'<\s$, the section $s$ is transverse to $0$ on $O\cap \ex F^{\s'}$.

 Let $\psi$  and $\psi'$ be as in Lemma \ref{psi}. One way to construct a section $s^{\s'}$ satisfying the conditions of definition \ref{nc section def} is by locally setting 
 \[s^{\s'}:=\psi(s^{\s}\circ \Phi_{\s'<\s})\]
Note that $\psi'(s^{\s'})$  is equal to the restriction of $s^{\s}$ plus a section of $Z_{d}$, so Lemma \ref{ZT1} implies that  $\psi'(s^{\s'})$ is transverse to $Z_{d}$ when restricted to some neighborhood of $f$.  The effect of $\psi'\circ \psi$ or $\psi\circ\psi'$ on sections is to add some section of $Z_{d}$ or $Z_{d'}$ respectively. It follows that $(\psi')^{-1}Z_{d}=Z_{d'}$. It also follows that $\psi'$ is transverse to every strata of any nice Whitney stratification of $Z_{d}$, so $\psi'$ pulls back any nice Whitney stratification of $\totl{Z_{d}}$ to a nice Whitney stratification of $\totl{Z_{d'}}$. 
 As $\psi'(s^{\s'})$ is transverse to the original nice Whitney stratification of $\totl{Z_{d}}$, it follows that $s^{\s'}$ is transverse to the pulled back nice Whitney stratification of $\totl{Z_{d'}}$.

\stop

 \
 
 \
 
 Note that the normally complex sections of $V$ are a module over the sheaf of $G$-equivariant, smooth, $\mathbb C$-valued functions on $\ex F$.  In particular, multiplying $s$ by  a $G$-equivariant, smooth, $\mathbb C$-valued function $h$ corresponds to multiplying $s^{\s}$ by $h\circ r_{\ex F^{\s}}$. The existence of smooth $G$-equivariant bump functions implies that normally complex sections may be constructed locally. 
 Because $\poly{}{I_{s}}{\ex N_{\ex F^{\s}}}{V\rvert_{\totl{\ex F^{\s}}}}$ is a Noetherian $\poly{}{I_{s}}{\ex N_{\ex F^{\s}}}{\mathbb C}$-module, it follows that the normally complex sections of $V$ are a finitely generated module over the sheaf of $G$-equivariant, smooth, $\mathbb C$-valued functions.

\subsection{Normally complex perturbations of $\dbar$}
 
 \

\begin{defn}\label{c def} A normally complex perturbation of $\dbar$ is a collection of  $G_{i}$-equivariant sections 
\[s_{i}:\ex F(\hat f_{i}^{\sharp})\longrightarrow V_{i}(\hat f_{i}^{\sharp})\]
where $\{\hat f_{i}\}$ comes from some compatible set of extensions, $(\mathcal U_{i}^{\sharp},V_{i},\hat f_{i}^{\sharp}/G_{i})$ of a  rigidified embedded Kuranishi structure 
 (definition \ref{rks def}), 
so that 
\begin{enumerate}
\item\label{c1} There exist open substacks $\mathcal O_{i}\subset\mathcal O_{i}^{\sharp}\mathcal U^{\sharp}_{i}\subset\Msw$ containing $\hat f_{i}$ so that \begin{itemize}
\item the intersection of $\hat f_{i}^{\sharp}$ with  $\mathcal O_{i}$  is $\hat f_{i}$ and the intersection of $\hat f_{i}^{\sharp}$ with $\mathcal O_{i}^{\sharp}$ is $\hat f^{\sharp}_{i}$

\item the closure within $\hat f_{j}^{\sharp}$ of $\hat f^{\sharp}_{j}\cap \mathcal O_{i}$ is contained within $\mathcal O_{i}^{\sharp}$
 \item and  at any curve $f$ in $\hat f_{i}^{\sharp}$ but not in  $\mathcal O_{j}$ for any $j$, 
\[\abs{\dbar f-s_{i}(f)}< \frac 12\abs{\dbar f}\]
\end{itemize}
\item\label{c2} If $\dim V_{i}\leq \dim V_{j}$, then there exists a $G_{i}\times G_{j}$-equivariant  section $s_{i}$ of $V_{i}(\hat f_{i}^{\sharp}\times_{\Msw}\hat f_{j}^{\sharp})$  so that the following diagram commutes:
\[\begin{tikzcd}V_{i}(\hat f_{i}^{\sharp})&\lar V_{i}(\hat f_{i}^{\sharp}\times_{\Msw}\hat f_{j}^{\sharp})\rar &V_{j}(\hat f_{j}^{\sharp})
\\ \ex F(\hat f_{i}^{\sharp})\uar{s_{i}}&\lar\uar{s_{i}} \ex F(\hat f_{i}^{\sharp}\times_{\Msw}\hat f_{j}^{\sharp})\rar&\uar{s_{j}}\ex F(\hat f_{j}^{\sharp})\end{tikzcd}\]
\item\label{c3}  If $\dim V_{i}\leq\dim V_{j}$, and $f$ is in $\hat f_{j}^{\sharp}$ and $\mathcal O_{i}^{\sharp}$,  
\[\abs{\dbar(f)-s_{j}(f)}_{V_{j}/V_{i}}\leq \epsilon \dim V_{j}\abs{\dbar f}_{V_{j}/V_{i}}\] 
where $\abs{\cdot}_{V_{j}/V_{i}}$ indicates the size of the orthogonal projection from $V_{j}(f)$ onto the orthogonal complement of $V_{i}(f)$ and $2\epsilon$ is strictly less than $\dim V_{k}$ for all $k$ so that $f\in \hat f_{k}^{\sharp}$.
\item\label{c4} $s_{i}$ is normally complex in the sense of Definition \ref{nc section def}.
\end{enumerate}

Say that a normally complex perturbation is transverse to $0$ if each $s_{i}$ is  transverse to zero in the sense of definitions \ref{nc section def} and  \ref{Z transverse def}. Say that a normally complex perturbation is compatible with a complex submersion $\Phi:\Msw\longrightarrow \ex X$ if the rigidified embedded Kuranishi structure is compatible with $\Phi$.

 \end{defn}
 
 Conditions \ref{c1} and \ref{c3} above are designed to play the role of a `compactly supported perturbation' condition to stop solutions to $s_{j}=0$  escaping from the region where we have control. Condition \ref{c1} ensures that any solution $f$ to the equation $ s_{j}(f)=0$ is contained in some $\mathcal O_{i}$. Condition \ref{c3} ensures that any solution $f$ to $s_{j}(f)=0$ which is contained in $\mathcal O_{i}$ must also be contained in $\hat f_{i}$. Condition \ref{r2} then implies that such an $f$ must also satisfy $s_{i}(f)=0$. 
 \begin{lemma}\label{cs}
 Given a normally complex perturbation of $\dbar$, the substack 
 \[ s^{-1}(0):=\bigcup_{i}\{f\in \hat f_{i}\text{ so that }s_{i}(f)=0\}/G_{i}\subset \Msw\]
 is a closed substack of $\Msw$.
 \end{lemma}
 
 \pf
 
 Because each Kuranishi structure is locally finite,  any curve $f$ in the closure of $s^{-1}(0)$ must be contained in the closure of some $s_{i}^{-1}(0)\subset \hat f_{i}$. One of the conditions of being an extension is that $\hat f_{i}^{\sharp}$ contains the closure of $\hat f_{i}$, therefore $f$ must be contained in $\hat f_{i}^{\sharp}$ and $s_{i}(f)=0$. As argued above, condition \ref{c1} then implies that $f$ is in $\mathcal O_{j}$ for some $j$, and then condition \ref{c3} implies that $f$ must in fact be contained in $\hat f_{j}$, and condition \ref{c2} implies that $s_{j}(f)=0$, therefore $f\in s^{-1}(0)$.
 
 \stop
 
 \begin{lemma}\label{open transverse} Given a normally complex perturbation of $\dbar$, suppose that $s_{i}$ is transverse to $0$ at a curve $f$ in $\hat f_{i}^{\sharp}$, and that $V_{j}(f)\supset V_{i}(f)$. Then $s_{j}$ is also transverse to $0$ at $f$.
 \end{lemma}

 \pf
 
  Consider the commutative diagram 
 \[\begin{tikzcd}V_{i}(\hat f_{i}^{\sharp})&\lar V_{i}(\hat f_{i}^{\sharp}\times_{\Msw}\hat f_{j}^{\sharp})\rar &V_{j}(\hat f_{j}^{\sharp})
\\ \ex F(\hat f_{i}^{\sharp})\uar{s_{i}}&\lar\uar{s_{i}} \ex F(\hat f_{i}^{\sharp}\times_{\Msw}\hat f_{j}^{\sharp})\rar&\uar{s_{j}}\ex F(\hat f_{j}^{\sharp})\end{tikzcd}\]
from item \ref{c2} of Definition \ref{c def}. As specified in Definition \ref{rks def}, the horizontal maps in the above diagram are normally rigid. As the left pointing maps are $G_{i}$-equivariant $G_{j}$-bundles (defined over an open subset of $\ex F(\hat f_{i}^{\sharp})$ or $V_{i}(\hat f_{i}^{\sharp})$ respectively), the middle map labeled $s_{i}$ is transverse to $0$ at $f$ if the original section $s_{i}$ is transverse to $0$ at $f$.

 Item \ref{c3} of Definition \ref{c def} implies that just like $\dbar$, $s_{j}$ restricted to $\mathcal O_{i}^{\sharp}$ intersects $V_{i}$ exactly  on the image of top righthand map above. (Recall that $V_{i}$  defines a sub vector bundle of $V_{j}(\hat f_{j}^{\sharp})$  which remains defined on some open neighborhood of $f$ in $\ex F(\hat f_{j}^{\sharp})$.) Item \ref{c3} of Definition \ref{c def} also gives that for any vector $v$ in $T_{f}\ex F(\hat f_{j}^{\sharp})$, 
 \[\abs{\nabla_{v} s_{j}-\nabla_{v} \dbar }_{V_{j}/V_{i}}\leq \frac 12\abs{\nabla_{v} \dbar }_{V_{j}/V_{i}}\]
As condition \ref{rk3} of Definition \ref{rks def} specifies that $\dbar\hat f_{j}^{\sharp}$ is transverse to $V_{i}$ at $f$,  the above inequality implies  that $s_{j}$ is transverse to $V_{i}$ at $f$. 
 
 Use the notation $\ex F$ to indicate a connected component of a small $G_{i}\times G_{j}$-equivariant neighborhood of $f$ in $\ex F(\hat f_{i}^{\sharp}\times_{\Msw}\hat f_{j}^{\sharp})$. (If this neighborhood is small enough, only $I_{f}$ will act on $\ex F$, but this does not affect our definition of transversality.)  Use the notation $\ex F^{\s}$ for the strata of $(\ex F,I_{f})$ containing $f$. Use the notation $V$ to indicate $V_{i}(\hat f_{i}^{\sharp}\times_{\Msw}\hat f_{j}^{\sharp})$ restricted to $\ex F$.

 As specified by Definition \ref{nrvb map}, we have fiberwise rigid maps
 \[\begin{tikzcd} V_{\ex F^{\s}}\dar \rar & V_{j}(\hat f_{j}^{\sharp})_{\ex F(\hat f_{j}^{\sharp})^{\s'}}\dar
 \\ \ex N_{\ex F^{\s}}\rar & \ex N_{\ex F(\hat f_{j}^{\sharp})^{\s'}}\end{tikzcd}\]
As the map $\ex N_{\ex F^{\s}}\longrightarrow \ex N_{\ex F(\hat f_{j}^{\sharp})^{\s'}}$ is a $I_{f}$-equivariant fiberwise rigid embedding, there exists some $I_{f}$-equivariant neighborhood $\ex N'$ of its image, where $\ex N'\longrightarrow M$ is some fiberwise rigid bundle with a $I_{f}$-equivariant, fiberwise rigid projection
\[\ex N_{\ex F^{\s}}\xleftarrow{\ \ \pi\ 
\ } \ex N' \]
covered by a $I_{f}$-equivariant fiberwise rigid vector bundle map
\[\begin{tikzcd}V_{\ex F^{\s}}\dar & \lar[swap]{\hat \pi} V'\dar
\\ \ex N_{\ex F^{\s}}&\lar[swap]{\pi} \ex N'\end{tikzcd}\]
where $V'$ is the restriction of $V_{i}(\hat f_{j}^{\sharp})_{\ex F(\hat f_{j}^{\sharp})}$ to $\ex N'\subset \ex N_{\ex F(\hat f_{j}^{\sharp})^{\s'}}$. Choose the map $\hat \pi$ to be a projection in the sense that the inclusion of $V_{\ex F^{\s}}$ into $V'$ followed by $\hat \pi$ is the identity. Moreover, choose $\hat \pi$ so that  $\hat \pi$ is an orthogonal projection when restricted to fibers over the image of $\ex N_{\ex F^{\s}}$.  By choosing $\ex N'$ small enough, we may assume that $\hat \pi$ is fiberwise surjective. Identify  $\pi^{*}V$ with the orthogonal complement of $\ker \hat \pi$, so there is a rigid splitting
\[V'=\ker\hat \pi\oplus \pi^{*}V\]
Note that $\pi^{*}V$ is the subvectorbundle corresponding to $V_{i}$ when restricted to the image of $\ex N_{\ex F^{\s}}$. In particular, this implies that $s_{j}:\ex N'\longrightarrow V'$ is transverse to $\pi^{*}V\subset V'$.

Choose $d$ large enough that the polynomials of $b$-degree $d$ generate $\poly{}{I_{f}}{\ex N_{\ex F^{\s}}}{V}$ and $\poly{}{I_{f}}{\ex N'}{V'}$. Consider the evaluation map
\[\poly{d}{I_{f}}{\ex N'}{V'}\longrightarrow V'=\ker{\hat \pi}\oplus \pi^{*}V\]
Choose a section $s'$ of $\poly{d}{I_{f}}{\ex N'}{V'}$ representing $s_{j}$.
\[\begin{tikzcd}\ex N'\rar{s'}\ar[bend left]{rr}{s_{j}}&\poly{d}{I_{f}}{\ex N'}{V'}\times_{M}\ex N'\rar&V' \end{tikzcd}\]

Define 
\[X\subset \poly{d}{I_{f}}{\ex N'}{V'}\times_{M}\ex N'\]
to be the inverse image of $\pi^{*}V\subset V'$ under the evaluation map. The fact that $s_{j}$ is transverse to $\pi^{*}V$ implies that $s'$ is transverse to $X$ and that $X$ is smooth on some neighborhood of its intersection with the graph of $s_{j}$.

 
Extend $\pi:\ex N'\longrightarrow \ex N_{\ex F^{\s}}$ to a map
\[\pi':\poly{d}{I_{f}}{\ex N'}{V'}\times_{M}\ex N'\longrightarrow \poly d{I_{f}}{\ex N_{\ex F^{\s}}}V\times _{\totl{\ex F^{\s}}} \ex N_{\ex F^{\s}}\]
 so that the following diagram commutes
 \[\begin{tikzcd}\poly{d}{I_{f}}{\ex N'}{V'}\times_{M}\ex N'\rar \dar{\pi'} & V'\dar{\hat \pi}
 \\ \poly d{I_{f}}{\ex N_{\ex F^{\s}}}V\times _{\totl{\ex F^{\s}}} \ex N_{\ex F^{\s}}\rar&V\end{tikzcd}\]
and so that at any point in the image of the inclusion $\ex N_{\ex F^{\s}}\longrightarrow \ex N'$, the map on polynomials is the natural one given by the inclusion $\ex N_{\ex F^{\s}}\longrightarrow \ex N'$ and the projection $\hat \pi:V'\longrightarrow V$.

 \
 
 Note that the intersection of $(\pi')^{-1}Z_{d}$ with $X$ is equal to $Z_{d}\subset \poly d{I_{f}}{\ex N'}{V'}\times_{M}\ex N'$. On a neighborhood of $(s'(f),f)$, $X$ is smooth and $\pi'$ restricted to $X$ is a submersion. It follows that any nice Whitney stratification of $Z_{d}$ pulls back under $\pi'$ to give a nice Whitney stratification of $Z_{d}$ in a neighborhood of $(s'(f),f)$.

 The inclusion $\ex N_{\ex F^{\s}}$ into $\ex N'$ followed by the graph of $s'$ followed by $\pi'$ defines a section $\ex N_{\ex F^{\s}}\longrightarrow \poly d{I_{f}}{\ex N_{\ex F^{\s}}}V$ which represents $s_{i}$, and which is therefore transverse at $f$ to some nice Whitney stratification of $Z_{d}\subset \poly d{I_{f}}{\ex N_{\ex F^{\s}}}V\times_{\totl{\ex F^{\s}}}\ex N_{\ex F^{\s}}$. The intersection of $s'$ with $X$ is therefore transverse to the pulled back nice Whitney stratification of $Z_{d}$ at $f$ when when $Z_{d}$ is considered as a subspace of $X$. As $s'$ is also transverse to $X$ at $f$, it follows that $s'$ is transverse $Z_{d}$ at $f$ in the sense of Definition \ref{Z transverse def}, so $s_{j}$ is transverse to $0$ at $f$.

\stop

\begin{lemma}\label{O choice}Given any rigidified embedded Kuranishi structure, $(\mathcal U_{i},V_{i},\hat f_{i}/G_{i})$ with extensions $(\mathcal U_{i}^{\sharp},V_{i},\hat f^{\sharp}_{i}/G_{i})$,   there exist open substacks $\mathcal O_{i}\subset\mathcal O_{i}^{\sharp}\subset\mathcal U^{\sharp}_{i}\subset\Msw$, containing $\hat f_{i}$ and     obeying the conditions of Definition \ref{c def} item \ref{c1}, so in particular \begin{itemize}
\item the intersection of $\hat f_{i}^{\sharp}$ with  $\mathcal O_{i}$  is $\hat f_{i}$, and the intersection of $\hat f_{i}^{\sharp}$ with $\mathcal O_{i}^{\sharp}$ is $\hat f_{i}^{\sharp}$
\item the closure within $\hat f_{j}^{\sharp}$ of $\hat f^{\sharp}_{j}\cap \mathcal O_{i}$ is contained within $\mathcal O_{i}^{\sharp}$ 
\end{itemize}
Additionally,    given a well ordering $\prec$ of Kuranishi charts, and extensions
\[\hat f_{i} \exte\hat f_{i}^{\sharp}\exte \hat f_{i,j}\exte\hat f_{i,\underline j}\]
for all $i,j$ so that
\[\hat f_{i,\underline k}\exte\hat f_{i,j} \text{ if }k\succ j\]
there exist open substacks $\mathcal O_{i,j}\subset \mathcal U_{i,j}$ and $\mathcal O_{i,\underline j}\subset\mathcal U_{i,\underline j}$ so that
\begin{itemize}
\item\[\mathcal O_{i}\subset\mathcal O_{i}^{\sharp}\subset\mathcal O_{i,\underline k}\subset \mathcal O_{i,j}\subset\mathcal O_{i,\underline j}\text { for all $k\succ j$}\]
\item For all $k\succeq j$,  
\[\mathcal O_{i,k}\cap\hat f_{i,\underline j}=\hat f_{i,k}\]
and
\[\mathcal O_{i,\underline k}\cap \hat f_{i,\underline j}=\hat f_{i,\underline k}\]
\item The closure within $\hat f_{j,\underline j}$ of $\mathcal O_{i,j}\cap \hat f_{j,\underline j}$ is contained in $\mathcal O_{i,\underline j}$
\end{itemize}

 \end{lemma}

\pf

 Note that the definition of an extension of a Kuranishi chart given in \cite{evc} tells us that $\mathcal U_{i}$ is an open substack of $\mathcal U_{i}^{\sharp}$ so that $\hat f_{i}$ is the intersection of $\hat f_{i}^{\sharp}$ with $\mathcal U_{i}$, and so that the closure within $\Msw$ of $\hat f_{i}$ is contained in $\hat f_{i}^{\sharp}$. 

 Similar statements hold for all the extensions
\[\hat f_{i}\exte \hat f_{i}^{\sharp}\exte \hat f_{i,k}\exte \hat f_{i,\underline k}\exte\hat f_{i,j}\exte \hat f_{i,\underline j} \text{ for }i\preceq j\prec k\]

 It follows that the image of $\hat f_{i,j}\times_{\Msw}\hat f_{j,\underline j}$ in $\hat f_{j,\underline j}$ has closure within $\hat f_{j,\underline j}$ contained in $\mathcal U_{i,\underline j}$. The curves in $\hat f_{j,\underline j}$ which are in the closure of $\mathcal U_{i,j}\cap \hat f_{j,\underline j}$ but are not in $\mathcal U_{i,\underline j}$ form a closed subset which does not intersect the closure of the image of $\hat f_{i,j}\times_{\Msw}\hat f_{j,\underline j}$. Therefore, there exists an  open substack $C_{i,j}$ of $\hat f_{j,\underline j}/G_{j}$ so that 
 \begin{itemize}
 \item  the closure $\bar C_{i,j}$ of $C_{i,j}$ within $\Msw$ does not intersect  $\hat f_{i,j}$. (This uses the fact  that there exists an extension of $\hat f_{j,\underline j}$  in order to know what the closure of $C_{i,j}$ outside of $\hat f_{j,\underline j}$ looks like.) 
  \item $C_{i,j}$ contains  every curve which is not in $\mathcal U_{i,\underline j}$ but which is in the closure within $\hat f_{j,\underline j}$ of $\hat f_{j,\underline j}\cap  \mathcal U_{i,j}$  
  
  \[C_{i,j}\supset \overline{\hat f_{j,\underline j}\cap  \mathcal U_{i,j}}\setminus\mathcal U_{i,\underline j}\] 
 \end{itemize}
  As every embedded Kuranishi structure is locally finite, we may choose $C_{i,j}$ to be empty for all but finitely many $j$. For each $k$ with $i\prec k\prec j$, $\bar C_{i,k}$ does not intersect  $\hat f_{i,k}$, which contains the closure of $\hat f_{i,j}$.  We may therefore  inductively require the following additional condition on $C_{i,j}$:
  \begin{itemize}
   \item $C_{i,j}$ contains $\hat f_{j,\underline j}\cap \bar C_{i,k}$ for all $i \prec k\prec j$
\end{itemize}

Now define 
\[\mathcal O_{i,j}:=\mathcal U_{i,j}\setminus \bigcup_{j\succeq k\succeq i}\bar C_{i,k}\]
\[\mathcal O_{i,\underline j}:=\mathcal U_{i,\underline j}\setminus\bigcup_{j\succ k\succeq i}\bar C_{i,k}\]

These open substacks $\mathcal O_{i,j}\subset\mathcal O_{i,\underline j}\subset \mathcal U_{i,\underline j}$ satisfy the required conditions. 

In order to have $O^{\sharp}_{i}$ a substack of $\mathcal O_{i,j}$ for all $j$, define
\[\mathcal O^{\sharp}_{i}:= \mathcal U^{\sharp}_{i}\setminus \bigcup_{j\succeq i} \bar C_{i,j}\]
We can similarly remove offending substacks from $\mathcal U_{i}$ to define $\mathcal O_{i}\subset\mathcal O_{i}^{\sharp}$ satisfying the required conditions.  In particular, we can choose an open substack $C_{i,j}'\subset \hat f_{j}^{\sharp}$ so that
\begin{itemize}
\item The closure $\bar C_{i,j}'$ of $C_{i,j}'$ within $\Msw$ does not intersect $\hat f_{i}$
\item $C_{i,j}'$ contains every curve in $\hat f_{j}^{\sharp}$ which is in the closure of $\mathcal U_{i}\cap\hat f_{j}^{\sharp}$ but which is not in $\mathcal O_{i}^{\sharp}$.
\end{itemize}
Again, we may choose $C'_{i,j}$ to be empty for all but finitely many $j$. Then define
\[\mathcal O_{i}:=\mathcal U_{i}\setminus \lrb{\bigcup_{j\succeq i}\bar C_{i,j}'\cup \bar C_{i,j}}\]

\stop

\begin{prop}\label{nc prop} Given any rigidified embedded Kuranishi structure on the holomorphic curves in $\Msw(\hat {\ex B})$ there exists a normally complex perturbation of  $\dbar$ in the sense of Definition \ref{c def} which is transverse to $0$. 

Any two such normally complex perturbations of $\dbar$ (defined using possibly different rigidified embedded Kuranishi structures) are cobordant in the sense that there exists a normally complex perturbation of $\dbar$ on some rigidified embedded Kuranishi structure on $\Msw(\hat{\ex B}\times {\mathbb R})$ which restricts to the given perturbations of $\dbar$ at $0$ and $1\in \mathbb R$, and which  is also transverse to $0$.

\end{prop}

 \pf  Choose a well  ordering $\prec$ on our Kuranishi charts and  extensions 
 \[\hat f_{i}\exte \hat f_{i}^{\sharp}\exte \hat f_{i,j}\exte \hat f_{i,\underline j}\]
 so that
\[   \hat f_{i,\underline k}\exte \hat f_{i,j} \text{ for all }k\succ j\]
Lemma \ref{well order} provides us with such extensions only for all $i\preceq j$, however the same argument gives extensions for all $i$ and $j$. 

 Then choose open substacks $\mathcal O_{i}\subset \mathcal O_{i}^{\sharp}\subset\mathcal O_{i,j}\subset\mathcal O_{i,\underline j}\subset \mathcal U_{i,\underline j}$  satisfying the conditions of  Lemma \ref{O choice}.

 \
 
Fix an index $j$. Suppose that a normally complex perturbation of $\dbar$ has been chosen on $\hat f_{i,\underline j}$ for all $i\prec j$, and that these normally complex perturbations of $\dbar$ are compatible in the sense of items \ref{c2} and \ref{c3} of Definition \ref{c def}, where $\hat f_{k,\underline j}$ is always used in place of $\hat f_{k}^{\sharp}$ as an extension of $\hat f_{k}$. Assume further that item \ref{c3} holds on the larger set $\mathcal O_{i,\underline j}\supset \mathcal O^{\sharp}_{i}$. 

Let $f$ be any curve in $\hat f_{j,\underline j}$. Suppose that $f\in \hat f_{i,\underline j}$, and that $ V_{i}$ has   the maximum dimension so that $i\prec j$ and $f\in \hat f_{i,\underline j}$. As the map $\hat f_{i,\underline j}\times_{\Msw}\hat f_{j,\underline j}\longrightarrow \hat f_{i,\underline j}$ is a $G_{i}$-equivariant $G_{j}$-bundle over an open subset of $\hat f_{i,\underline j}$, the section $s_{i}$ lifts uniquely to a $G_{i}\times G_{j}$-equivariant section $s_{i}$ of $V_{i}(\hat f_{i,\underline j}\times_{\Msw}\hat f_{j,\underline j})$ so that the following diagram commutes
\[\begin{tikzcd}V_{i}(\hat f_{i,\underline j})&\lar V_{i}(\hat f_{i,\underline j}\times_{\Msw}\hat f_{j,\underline j})
\\ \ex F(\hat f_{i,\underline j})\uar{s_{i}}&\lar\uar{s_{i}} \ex F(\hat f_{i,\underline j}\times_{\Msw}\hat f_{j,\underline j})\end{tikzcd}\]
The map $\hat f_{i,\underline j}\times_{\Msw}\hat f_{j,\underline j}\longrightarrow \hat f_{j,\underline j}$ is a $G_{j}$-equivariant $G_{i}$ bundle over its image, which is a sub exploded manifold of $\hat f_{j,\underline j}$. The fact that the maps 
\[\begin{tikzcd}V_{i}(\hat f_{i,\underline j}\times_{\Msw}\hat f_{j,\underline j})\rar\dar& V_{j}(\hat f_{j,\underline j})\dar
\\ \ex F(\hat f_{i,\underline j}\times_{\Msw}\hat f_{j,\underline j})\rar & \ex F(\hat f_{j,\underline j})\end{tikzcd}\]
 are normally rigid implies that on a neighborhood of $f$ in $\ex F(\hat f_{j,\underline j})$, there exists a normally complex section $s_{j}$ of $V_{j}(\hat f_{j,\underline j})$ so that the required diagram commutes:
\[\begin{tikzcd}V_{i}(\hat f_{i,\underline j}\times_{\Msw}\hat f_{j,\underline j})\rar& V_{j}(\hat f_{j,\underline j})
\\ \uar{s_{i}}\ex F(\hat f_{i,\underline j}\times_{\Msw}\hat f_{j,\underline j})\rar & \ex F(\hat f_{j,\underline j})\uar{s_{j}}\end{tikzcd}\]
At $f$, such a  section $s_{j}$ may be chosen so that $\nabla s_{j}$  is equal to an arbitrary $I_{f}$-invariant map which 
\begin{itemize}
\item is equal to $\nabla s_{i}$ on the image of $T_{f}\ex F(\hat f_{i,\underline j}\times_{\Msw}\hat f_{j,\underline j})$ within $T_{f}\ex F(\hat f_{j,\underline j})$
\item is complex on the $I_{f}$ anti-invariant vectors within $T_{f}\ex F(\hat f_{j,\underline j})$.
\end{itemize}
  As specified by condition \ref{rk3} of Definition \ref{rks def},  $D\pi_{V_{i}}\dbar$ at $f$ is surjective and approximately holomorphic in the sense that 
\[\abs{J\nabla_{v}\dbar -\nabla_{Jv}\dbar}_{V_{j}/V_{i}}\leq \epsilon \abs{\nabla_{v}\dbar}_{V_{j}/V_{i}}\]
 where $v$ indicates any $I_{f}$-anti-invariant vector in $T_{f}\ex F(\hat f_{j,\underline j})$ and $\epsilon$ is a constant strictly less than $\frac 1{2\dim V_{j'}}$ for all $j'$ so that $f\in \hat f_{j',\underline j}$. Similarly, for any $k\prec j$ so that $f$ is in  $\hat f_{k,\underline j}$, $D\pi_{V_{k}}\dbar$ is surjective at $f$ and the analogous inequality holds.
 \begin{equation}\label{k}\abs{J\nabla_{v}\dbar -\nabla_{Jv}\dbar}_{V_{j}/V_{k}}\leq \epsilon \abs{\nabla_{v}\dbar}_{V_{j}/V_{k}}\end{equation}
  Our freedom in specifying the derivative of $s_{j}$ normal to the image of $\hat f_{i,\underline j}\times_{\Msw}\hat f_{j,\underline j}$ and the observation that $\abs{\dbar}_{V_{j}/V_{i}}=0$ on this image implies that $s_{j}$ may be chosen so that at $f$, 
 \begin{equation}\label{ie}\abs{\nabla_{v}(\dbar-s_{j})}_{V_{j}/V_{i}}\leq \epsilon \abs{\nabla_{v}\dbar}_{V_{j}/V_{i}}\end{equation}  
 where the above inequality is now valid for any vector $v\in T_{f}\ex F(\hat f_{j,\underline j})$. Similarly, our freedom in choosing the derivative of $s_{j}$  normal to  the image of $T_{f}\ex F(\hat f_{i,\underline j}\times_{\Msw}\hat f_{j,\underline j})$ implies that we may choose $s_{j}$ so that at $f$, 
  \begin{equation}\label{ke}\abs{\nabla_{v}(\dbar-s_{j})}_{V_{j}/V_{k}}\leq \epsilon \abs{\nabla_{v}\dbar}_{V_{j}/V_{k}}\end{equation}
where the above inequality now holds for any vector $v\in T_{f}\ex F(\hat f_{j,\underline j})$ normal to  the image of $T_{f}\ex F(\hat f_{i,\underline j}\times_{\Msw}\hat f_{j,\underline j})$.
 
   In particular, as both $\dbar$ and $s_{j}$ are in $V_{i}$ when restricted to the image of $\hat f_{i,\underline j}\times_{\Msw}\hat f_{j,\underline j}$, inequality (\ref{ie}) implies that on some neighborhood of $f$,  our section $s_{j}$ satisfies the inequality required in condition \ref{c3} of Definition \ref{c def}. In particular, for all $f'$ in some neighborhood of $f$ within $\hat f_{j,\underline j}$, 
 \[\abs{\dbar(f')-s_{j}(f')}_{V_{j}/V_{i}}\leq \epsilon \abs{\dbar f'}_{V_{j}/V_{i}}\]
where $\epsilon$ is some new constant strictly less than $\frac 1{2\dim V_{j'}}$ for all $j$ so that $f'\in\hat f_{j',j}$. (Note that this condition previously involved $j$ so that $f\in \hat f_{j',\underline j}$. We now use the shrunken $\hat f_{j',j}\exte\hat f_{j',\underline j}$ so that our condition may hold in a neighborhood of $f$.)

 If $f$ is also in $\mathcal O_{k,\underline j}$, we need to verify that $s_{j}$ is also compatible  with $s_{k}$ in the above fashion. The only case in which this is a new inequality is if $\dim V_{k}<\dim V_{j}$. If $f$ is in $\hat f_{k,\underline j}$, then the fact that $s_{i}$ obeys condition \ref{c3} on $\mathcal O_{k,\underline j}$ together with inequality (\ref{ke}) imply that for all $f'$ in some neighborhood of $f\in\hat f_{j,\underline j}$,  
 \[\abs{\dbar(f')-s_{j}(f')}_{V_{j}/V_{k}}\leq \epsilon(\dim V_{k}+1) \abs{\dbar f'}_{V_{j}/V_{k}}\leq \epsilon(\dim V_{j}) \abs{\dbar f'}_{V_{j}/V_{k}} \]
where   $\epsilon<\frac 1{2\dim V_{j'}}$ for all $j'$ so that $f'\in\hat f_{j',j}$.

On the other hand, if $f$ is in $\mathcal O_{k,\underline j}$ but not in $\hat f_{k,\underline j}$, then $\dbar f$ is not in $V_{k}$, therefore the fact that $s_{i}$ obeys condition \ref{c3} on $\mathcal O_{k,\underline j}$ again implies that the above inequality holds at $f$, and therefore holds in a neighborhood of $f$.

 The section $s_{j}$ thus constructed will obey all the requirements of Definition \ref{c def} on some neighborhood of $f$. Note that if $f$ is not in $\mathcal O_{i,\underline j}$, item \ref{c3}  may not hold for the index $i$, however some neighborhood of $f$ in $\hat f_{j,\underline j}$ will then not intersect $\mathcal O_{i,j}$, therefore  item \ref{c3} will trivially hold on  $\mathcal O_{i,j}\subset\mathcal O_{i,\underline j}$. Similarly, if $f$ is not in $\hat f_{i,\underline j}$, item \ref{c2} may not hold for the index $i$, but then some neighborhood of $f$ will not intersect $\hat f_{i,j}\exte\hat f_{i,\underline j}$, so item \ref{c2} will hold trivially for $\hat f_{i,j}$ in a neighborhood of $f$.

\

We have locally constructed $s_{j}$ around any $f$ which is contained in some $\hat f_{i,\underline j}^{\sharp}$ for $i\prec j$. Now suppose that $f\in \hat f_{j,\underline j}$ is not contained in $\hat f_{i,\underline j}$ for any $i\prec j$. Then there exists a normally complex section of $V_{j}(\hat f_{j,\underline j})$ which is equal to $\dbar f$ at $f$. On some open neighborhood of $f$, such a section will obey condition \ref{c1} of Definition \ref{c def}, because either $\dbar f=0$, in which case $f$ is in some $\mathcal O_{k}$ on an open neighborhood, or   $\abs{\dbar f}>0$ and $\dbar f-s_{j}(f)=0$, in which case, for all $f'$ in a neighborhood of $f$, 
\[\abs{s_{j}(f')-\dbar(f')}<\frac 12\abs{\dbar f}\]
Similarly, such a section will obey condition \ref{c3} on $\mathcal O_{i, j}$ for any $i\prec j$,  because  the fact that $f$ is not in $\hat f_{i,\underline j}$ implies that $\dbar f$ is not in $V_{i}(f)$. (If $V_{i}$ is not defined at $f$, then some open neighborhood of $f$ within $\hat f_{j,\underline j}$ will not intersect $\mathcal O_{i,j}$ at all, therefore the compatibility condition \ref{c3} will hold trivially.) Such a section will also obey condition \ref{c2} trivially for the slightly shrunken extensions $\hat f_{i,j}\subset \hat f_{i,\underline j}$.

 Average the locally defined sections $s_{j}$ above using a partition of unity. Each of the conditions of Definition \ref{c def} is preserved by averaging, therefore the averaged section $s_{j}$ will obey all the required conditions of Definition \ref{c def}. By averaging, we  may also ensure that $s_{j}$ is $G_{j}$-equivariant. It remains to show that if each  $s_{i}$ is transverse to $0$, such a  $s_{j}$ may be chosen to be transverse to $0$ too.

 As $s_{i}$ is transverse to zero for all $i\prec j$, Lemma \ref{open transverse} implies that  $s_{j}$ is also transverse to $0$ on the closure of $\mathcal O_{i,j}\cap \hat f_{j,j}$ for all $i\prec j$. On the complement of all such $O_{i,j}$,  conditions \ref{c2} and \ref{c3} do not apply. The only remaining condition, condition \ref{c1}, is an open condition, because $\dbar$ is nonzero on the region on which it applies. Our $s_{j}$ may therefore be perturbed outside the closure of $\mathcal O_{i,j}$ to obtain a  section transverse to $0$ obeying all the required compatibility conditions on $\mathcal O_{i,j}$. 
 
 The construction of a perturbation of $\dbar$ transverse to $0$ may then proceed by transfinite induction. In the end, condition \ref{c3} will  hold on $\mathcal O_{i}$ which is contained in $\mathcal O_{i,j}$ for all $j$, and condition \ref{c2} will only hold on the  extension $\hat f^{\sharp}_{i}$ of $\hat f_{i}$ which is contained in $\hat f_{i,j}$ for all $j$.

  \

Given any two normally complex perturbations of $\dbar$, it is proved in \cite{evc} that the corresponding Kuranish structures are cobordant, and Theorem \ref{rks thm} states that all choices required to rigidify the two corresponding Kuranishi structures extend over the cobordism. We can also ensure that the cobordism is trivial in a  neighborhood of our two given normally complex perturbations of $\dbar$.  On this region where the cobordism is trivial, the neighborhoods $\mathcal O_{i}\subset\mathcal O_{i}^{\sharp}$ may then be chosen to be the trivial extensions of the corresponding neighborhoods used in defining the given normally complex perturbations of $\dbar$. (The reason that these neighborhoods $\mathcal O_{i}$ may chosen to be trivial extensions in a neighborhood is that their construction only requires the removal of offending  substacks coming from other $\hat f_{j}^{\sharp}$. All the new Kuranishi charts in our cobordism may be chosen to not intersect a neighborhood of each end of the cobordism, so nothing further needs be removed on some neighborhood of each end of the cobordism.) Then the two given normally complex perturbations of $\dbar$ extend to normally complex perturbations of $\dbar$ on a neighborhood.  The above argument  may then be applied to construct another normally complex perturbation of $\dbar$, and then an averaging procedure may be used to get normally complex perturbation of $\dbar$ which restricts to our two given normally complex perturbations of $\dbar$. We must now show that this normally complex perturbation of $\dbar$ may be chosen to be transverse to $0$.

In our construction of a normally complex perturbation of $\dbar$ above, if $V_{i}$ has the lowest dimension so that $f$ is in $\hat f_{i}^{\sharp}$, then the equivariant sections of $\poly d{I_{f}}{\ex F(\hat f_{i}^{\sharp})^{\s}}{V_{i}}$ which correspond to normally complex perturbations of $\dbar$ are open. It follows that any normally complex perturbation of $\dbar$ extends to finite dimensional family of normally complex perturbations of $\dbar$ which is transverse to $0$ at $f$.  It follows that a generic normally complex perturbation is transverse to $0$ at $f$. As our normally complex perturbation of $\dbar$ is already transverse to $0$ when restricted to the two given normally complex perturbations, it may be modified elsewhere to be transverse to $0$ without affecting the property that it restricts to our two given normally complex perturbations.

  \stop

\subsection{Components of $\Msw$}

\

To every curve $f$ in $\Msw(\hat{\ex B})$, we can naturally associate a polytope $P(f)$, and a weak isotropy group $I_{f}$. 

The polytope $P(f)$ should be regarded as parametrizing the moduli space of variations of $\totb f$. For example $\ex B$ is a smooth manifold, then $P(f)$ is equal to $[0,\infty)^{k}$, where $k$ is the number of nodes of the domain of $f$. 

\begin{defn}\label{Pdef} Given any curve $f\in\Msw(\hat {\ex B})$, define $P(f)$ as follows: Choose any family of curves $\hat f$ containing $f$ with universal tropical structure. Define $P(f)$ 
to be the tropical structure of $\ex F(\hat f)$ at $f$.
\end{defn} 

Such families with universal tropical structure were constructed in \cite{uts}--- the uniqueness of $P(f)$ follows from the universal property of universal tropical structures. Note that in particular, if $(\mathcal U, V,\hat f/G)$ is a Kuranishi chart containing $f$, then $\hat f$ has universal tropical structure, so $P(f)$ is the tropical structure of $\ex F(\hat f)$ at $f$. 

\begin{defn}\label{Ifdef}For any $f\in\Msw(\hat{\ex B})$, define the weak isotropy group $I_{f}$ of $f$ as follows: Given any family of curves $\hat f$ with universal tropical structure containing $f$, let $\hat f^{\circ}\subset \hat f$ indicate the subfamily consisting of curves sent to the same point as $f$ under the projection $\hat f\longrightarrow \totl{\ex F(\hat f)}$. Define $I_{f}$ to be the automorphism group of $\hat f^{\circ}$.
\end{defn}

In particular, the moduli stack of curves topologically indistinguishable from $f$ in $\Msw$ is equal to $\hat f^{\circ}/I_{f}$. Note that $I_{f}$ acts naturally on $T_{f}\Msw$.
If $(\mathcal U,V,\hat f/G)$ is a Kuranishi chart containing $f$, then the weak isotropy group $I_{f}$ of $f$ is the weak isotropy group of $f$ in $(\ex F(\hat f),G)$. In other words, $I_{f}$ is the subgroup of $G$ which preserves $\totl f\subset \totl {\ex F(\hat f)}$. If $f$ is a curve in a smooth manifold, then $I_{f}$ is equal to the group of automorphisms of $\totl{f}$.

\

\begin{defn} Define a component  $\Msw_{\co}\subset\Msw(\hat {\ex B})$ to be a connected component of the substack consisting of curves $f$ with some fixed $I_{f}=I_{\co}$ and $P(f)=P(\co)$. Denote the component containing $f$ to be $\Msw_{\co(f)}$.
\end{defn}

Quite a lot of information is hidden in $\co$. In the case that the target is a smooth manifold $X$, the number of smooth components of curves in $\Msw_{\co}(X)$ is equal to $\dim P(\co)$. Each of these components is some branched cover of some underlying curve in $X$. The component $\co$ contains all topological information about these underlying curves, the homotopy type of the branched covers, and the combinatorics of how these covers are glued together to form curves in $\Msw_{\co}$. Similar information is contained in $\co$ in the general case that the target is a family of exploded manifolds $\hat{\ex B}$, however now each smooth component will not necessarily be a branched cover of a curve in $\hat{\ex B}$.

Given holomorphic curve $f\in\Msw_{\co}(\ex B)$,   $D\dbar_{f}$ represents a class in real $I_{f}$-equivariant $K$ theory, and its complex linear part,  $D\dbar^{\mathbb C}_{f}$ represents a class in complex $I_{f}$-equivariant $K$ theory, which is equal to $[D\dbar_{f}]$ in real $I_{f}$-equivariant $K$-theory.  

\begin{lemma} Let $f\in\Msw(\ex B)$ be any holomorphic curve. Then a formula for the class in complex $I_{f}$-equivariant $K$ theory represented by $D\dbar_{f}^{\mathbb C}$ is as follows:

\begin{equation}\label{index}[D\dbar^{\mathbb C}_{f}]=H^{0}(T^{*}\ex C\otimes K)^{*}-H^{0}(T\ex C)+n(\mathbb C-H^{0}(K)^{*})+\sum_{C_{i}}\lrb{\frac{\abs {I_{i}}}{\abs{I_{f}}}\int_{C_{i}}f^{*}c_{1}}\mathbb C[I_{f}/I_{i}] \end{equation}

where
\begin{itemize}
\item $\ex C$ is the domain of $ f$, 
\item $T^{*}\ex C$ indicates the sheaf of holomorphic sections of $T^{*}\ex C$.
\item $K$ is the sheaf of holomorphic sections of $T^{*}\ex C$ which vanish on external edges of $\ex C$.
\item The dimension of $\ex B$ is $2n$.
\item The sum is over orbits $C_{i}$ of the action of $I_{f}$ on the set of smooth components of $\totl{\ex C}$
\item $I_{i}$ is the subgroup of $I_{f}$ which fixes $C_{i}$ pointwise - so the orbit of a generic point in $C_{i}$ has size $\abs{I_{f}/I_{i}}$.
\item $c_{1}$ is a differential form on $\ex B$ representing the first Chern class of the almost complex structure $J$. (As noted in \cite{dre}, the integral of $c_{1}$ over each component of $\totl{\ex C}$ is well defined, and an integer. $C_{i}$ factors through  a degree $\abs{I_{f}/I_{i}}$-fold map to an (incomplete) holomorphic curve, so $\frac{\abs {I_{i}}}{\abs{I_{f}}}\int_{C_{i}}f^{*}c_{1}$ is an integer.)
\item $\mathbb C[I_{f}/I_{i}]$ indicates the complex $I_{f}$-representation corresponding to the group ring of $I_{f}/I_{i}$.
\end{itemize}
\end{lemma}

\pf 
We shall use a version of Serre duality for $\ex C$ which follows from usual Serre duality applied to each component of $\totl{\ex C}$: Given any holomorphic vector bundle $E$ over $\ex C$, then
\[H^{1}(E)=H^{0}(E^{*}\otimes K)^{*}\]
Alternatively, the kernel of the adjoint of the $\dbar$ operator acting on sections of $E$ is equal to $H^{0}(E^{*}\otimes K)$. The isomorphism in this case is as follows: the range of $\dbar$ consists of sections of $ \lrb{T^{*}\ex C\otimes E}^{(0,1)}$ which vanish on edges of $\ex C$. Such a section pairs with  an element of $H^{0}(E^{*}\otimes K)$ to give a complex valued two-form  which may then be integrated to give a complex number. This pairing describes $H^{0}(E^{*}\otimes K)$ as a subspace of the dual of the range of $\dbar$. This pairing vanishes on  $\dbar$ of any section, and describes $H^{0}(E^{*}\otimes K)$ as the kernel of the adjoint   of $\dbar$ acting on sections of  $E$.

This version of Serre duality may be familiar to some readers as Serre duality for nodal curves:  $K$ corresponds to  the sheaf of holomorphic sections of the cotangent bundle of the nodal curve $\totl{\ex C}$ which allow poles  with opposite residues at each side of a node.

The part of $T_{f}\Msw$ corresponding to variation of complex structure of the domain and quotienting out by complex automorphisms is equal to $H^{1}(T\ex C)-H^{0}(T\ex C)$, (which is the negative of the index of the $\dbar$ operator acting on smooth sections of $T\ex C$). Serre duality then implies that  $H^{1}(T\ex C)=H^{0}(T^{*}\ex C\otimes K)^{*}$, so the first two terms in the right hand side of equation \ref{index} come from variations of complex structure and infinitesimal automorphisms of $\ex C$.

 It remains to prove that the last two terms in the above expression correspond to the index of $D\dbar_{f}^{\mathbb C}$ restricted to sections of $f^{*}T\ex B$.  Note that the index of $D\dbar_{f}^{\mathbb C}$ depends only on the $I_{f}$-equivariant homotopy class of the complex vector bundle  $f^{*}T\ex B$. Consider the topological space $\totl{\ex C}/I_{f}$. Our $I_{f}$-equivariant vector bundle $f^{*}T\ex B$ is the pullback of a complex vector bundle from $\ex C\longrightarrow \totl{\ex C}/I_{f}$, and any complex vectorbundle on $\totl{\ex C}/I_{f}$ pulls back to an $I_{f}$-equivariant vector bundle. Homotopy classes of complex vector bundles on $\totl{\ex C}/I_{f}$ are classified by their dimension and first Chern class $c_{1}\in H^{2}(\totl{\ex C}/I_{f})$. There is one generator for $H^{2}(\totl{\ex C}/I_{f})$ for each $I_{f}$-orbit $C_{i}$ on the set of smooth components of $\totl {\ex C}$. It follows that homotopy classes of equivariant  complex vector bundles on $\ex C$ are classified by a choice of dimension, and an integer $a_{i}$ for each $C_{i}$. (The bundle of interest, $f^{*}T\ex B$ is classified by the integers $a_{i}=\frac{\abs {I_{i}}}{\abs{I_{f}}}\int_{C_{i}}f^{*}c_{1}$.)
 
 Consider the case of a holomorphic  line bundle with all $a_{i}<0$. In particular, define a divisor $D$ on $\totl{\ex C}/I_{f}$ by choosing $\abs{a_{i}}$ distinct smooth points on the ith component, and let our vector bundle be the pullback of $O(D)$--- so our vector bundle is $O(D')$ where $D'$ is the pullback of $D$ to $\ex C$. As holomorphic sections must intersect the zero section positively, $H^{0}(D')$ consists only of the zero section. Serre duality implies that $H^{1}(O(D'))=H^{0}(O(-D')\otimes K)^{*}$. Interpret elements of $H^{0}(O(-D')\otimes K)$ as holomorphic one forms which have simple poles at the points of $D'$, and which vanish on external edges of $\ex C$ as required of sections of $K$. There is a natural $I_{f}$-equivariant exact sequence
 \[0\longrightarrow H^{0}( K)\longrightarrow H^{0}(O(-D)\otimes K)\longrightarrow \oplus_{x\in D'} \mathbb C\longrightarrow \mathbb C\longrightarrow 0\] 
where the middle arrow corresponds to taking residues at the points of $D'$, and the last arrow is the sum of the residues (which must be $0$ for a holomorphic $1$ form which vanishes on external edges of $\ex C$.) Because $H^{1}(K)=H^{0}(\mathbb C)^{*}$, which is $1$ dimensional, the only obstruction to constructing a holomorphic $1$ form with specified residues is that the sum of the residues must be $0$. Therefore the above sequence is indeed an exact sequence. The action of $I_{f}$ on $\oplus _{x\in D'}\mathbb C$ is given by the action of $I_{f}$ on $D'$, therefore in $I_{f}$-equivariant $K$ theory, 
\[\oplus_{x\in D'}\mathbb C=\sum_{C_{i}}-a_{i}\mathbb C[I_{f}/I_{i}]\]
We actually need to calculate $H^{0}(O(-D')\otimes K)^{*}$, but the above complex representation is isomorphic to its dual. Therefore, the index of the $\dbar$ operator acting on a holomorphic line bundle with all $a_{i}<0$ is 
\[\mathbb C-H^{0}(K)^{*}+\sum_{C_{i}}a_{i}\mathbb C[I_{f}/I_{i}]\]
 Given a $n$-dimensional $I_{f}$-equivariant complex vector bundle with all $a_{i}\leq -n$, the index of $\dbar$ is
\[n(\mathbb C-H^{0}(K)^{*})+\sum_{C_{i}}a_{i}\mathbb C[I_{f}/I_{i}]\]
 As  a $2$ dimensional vector bundle with all $a_{i}$ sufficiently negative is homotopic to a sum of two holomorphic line bundles, one of which satisfies the above formula, and the other of which may have positive $a_{i}$, the above formula must also hold for all $I_{f}$-equivariant line bundles regardless of the values of $a_{i}\in \mathbb Z$. As all $I_{f}$-equivariant complex vectorbundles on $\ex C$ are homotopic to a sum of holomorphic line bundles, the above formula must be valid for all $I_{f}$-equivariant complex vector bundles. Applying this formula to $D\dbar_{f}^{\mathbb C}$ restricted to $f^{*}T\ex B$ gives the required formula for the index of $D\dbar_{f}^{\mathbb C}$:
 
 \[[D\dbar^{\mathbb C}_{f}]=H^{0}(T^{*}\ex C\otimes K)^{*}-H^{0}(T\ex C)+n(\mathbb C-H^{0}(K)^{*})+\sum_{C_{i}}\lrb{\frac{\abs {I_{i}}}{\abs{I_{f}}}\int_{C_{i}}f^{*}c_{1}}\mathbb C[I_{f}/I_{i}] \]
 
 \stop

\

\begin{remark}\label{multiple cover} In the case that $\ex C$ is a smooth Riemann surface, the sheaf $T^{*}{\ex C}=K$, and the Chevalley-Weil formula from \cite{CW} computes $H^{0}(K)$ and $H^{0}(K\otimes K)$ explicitly. In the simplest case that $\ex C$ is smooth and $I_{f}$ acts freely on $\ex C$, then 
\[H^{0}(K)=\mathbb C\oplus  \mathbb C[I_{f}]^{\oplus (g-1)}\] and 
\[H^{0}(K\otimes K)=\mathbb C[I_{f}]^{\oplus3(g-1)}\] 
where $g$ is the genus of $\ex C/I_{f}$. In this case
\[[D\dbar_{f}^{\mathbb C}]=\lrb{(3-n)(g-1)+c_{1}}\mathbb C[I_{f}]\]
where $c_{1}$ is the integral of the first Chern class over $\ex C/I_{f}$, $g$ is the genus of $\ex C/I_{f}$ and $2n$ is the dimension of the target $\ex B$.
\end{remark}

The right hand side of equation \ref{index} makes sense for any curve in $\Msw_{\co}$, and each term does not depend on which curve $f\in \Msw_{\co}$ is chosen, therefore the class in  $I_{\co}$-equivariant complex $K$ theory which is the right hand side of equation \ref{index} may be associated to $\Msw_{\co}$. (Actually, strictly speaking this class in $I_{\co}$-equivariant $K$ theory depends on a choice of isomorphism of $I_{\co}$ with $I_{f}$, so we can only naturally associate a class in $I_{\co}$-equivariant $K$ theory up to the action of automorphisms of $I_{\co}$.)

 Note that in the case of a family of targets $\hat{\ex B}\longrightarrow \ex B_{0}$, equation \ref{index} gives a formula for the index of $D\dbar_{f}$ relative to $\ex B_{0}$. The index of $D\dbar_{f}:T_{f}\Msw(\hat {\ex B})\longrightarrow \Y(f)$ is given by the right hand side of equation \ref{index} plus the tangent space to $\ex B_{0}$ at the image of $f$.

\

Consider a normally complex perturbation $s$ of $\dbar$ which is transverse to the zero section in the sense of Definition \ref{c def}. Use the notation 
\[\mathcal N:=s^{-1}(0)\subset\Msw\]
 We shall divide $\mathcal N$ into components, 
\[\mathcal N_{\co}:=\mathcal N\cap \Msw_{\co}\] each of which shall be smooth and have codimension $2$ boundary. The  homology classes associated to these different components of $\mathcal N$  can then be used to define some integer-valued holomorphic curve invariants.

The tangent space $T_{f}\Msw_{\co}\subset T_{f}\Msw$ consists of the span  of the  $I_{f}$  invariant vectors of $T_{f}\Msw$ together with the $\mathbb R$-nil vectors in $T_{f}\Msw$. (These $\mathbb R$-nil vectors correspond to the image in $T_{f}\Msw$ of $T_{f}\hat f^{\circ}$, where $\hat f^{\circ}$ is as in Definition \ref{Ifdef}. As $I_{f}$ is only the weak isotropy group, such $\mathbb R$-nil vectors may not be individually fixed by $I_{f}$, but the set of $\mathbb R$-nil vectors is preserved by $I_{f}$.) If $f$ is a holomorphic curve, $D\dbar_{f}:T_{f}\Msw\longrightarrow \Y(f)$ is $I_{f}$ equivariant and vanishes on $\mathbb R$-nil vectors. It follows that $D\dbar_{f}$ sends $T_{f}\Msw_{\co}$ to the $I_{f}$-invariant part of $\Y(f)$. We shall show that index of this operator \[D\dbar_{f}:T_{f}\Msw_{\co}\longrightarrow \Y(f)^{I_{f}}\] is the dimension of $\mathcal N_{\co}$. If $f\in\Msw_{\co}$ is not holomorphic, then a connection on $\Msw$ close to $f$ is needed to define $D\dbar_{f}$, however the index of $D\dbar_{f}:T_{f}\Msw_{\co}\longrightarrow Y(f)^{I_{f}}$ does not depend on the choice of connection or $f\in\Msw_{\co}$.

\begin{lemma}\label{component dimension} 
  $\mathcal N_{\co}\subset \Msw$ is smooth and has  $d(\co)$ dimensional smooth part in the sense that it is locally represented by  the quotient of  smooth  family of curves $\hat g$ by an action of $I_{\co}$ which is trivial on $\totl{\ex F(\hat g)}$, and the dimension of $\totl{\ex F(\hat g)}$ is $d(\co)$.

 The dimension, $d(\co)$ depends only on $\co$. If there is a holomorphic curve $f$ in $\Msw_{\co}$, then $d(\co)$ is equal to the index if the $I_{f}$-invariant part of $D\dbar_{f}$ minus the dimension of the space of $I_{f}$-invariant $\mathbb R$-nil vectors in $T_{f}\Msw$. 
 If there is no holomorphic curve in $\Msw_{\co}$, and $f$ is a curve in a family of targets $\hat{\ex B}\longrightarrow \ex B_{0}$,  then for any curve in $f$ in  $\Msw_{\co}$, $d(\co)$ is equal to the dimension of the invariant vectors in the right hand side of equation \ref{index} plus the dimension of $\ex B_{0}$ minus the dimension of the space of $I_{f}$-invariant $\mathbb R$-nil vectors in $T_{f}\Msw$
  
\end{lemma}

\pf

As this lemma is local, we may work locally around a curve $f\in \mathcal N_{\co}$ in a single Kuranishi chart $(\mathcal U, V,\hat f/G)$. Denote by $\ex F^{\co}\subset\ex F(\hat f)$ the subset consisting of curves in $\Msw_{\co}$. This subset $\ex F^{\co}$ consists of some union of connected components of $\ex F(\hat f)^{\s(f)}$, so it is has a smooth exploded structure. 

 Choose $d$ large enough that the degree $d$ polynomial sections generate $\poly {}{I_{f}}{\ex N_{\ex F(\hat f)^{\s(f)}}}{V}$. Denote by $Z^{\co}$ the subset of $Z\subset\poly {d}{I_{f}}{\ex N_{\ex F(\hat f)^{\s(f)}}}{V}\times_{\totl{\ex F(\hat f)^{\s(f)}}} \ex N_{\ex F(\hat f)^{\s(f)}}$ over $\ex F^{\co}\subset \ex N_{\ex F(\hat f)^{\s(f)}}$.  Restricted to $\ex F^{\co}$, every $I_{f}$-equivariant section of $V$ vanishes on the $I_{f}$-anti-invariant part of $V$, so $Z^{\circ}\subset \poly {d}{I_{f}}{\ex N_{\ex F(\hat f)^{\s(f)}}}{V}\times_{\totl{\ex F(\hat f)^{\s(f)}}}\ex F^{\co}$ is determined by the condition that the section vanishes on the $I_{f}$-invariant part of $V$. The evaluation map at $f$ from $\poly {d}{I_{f}}{\ex N_{\ex F(\hat f)^{\s(f)}}}{V}$ to $V(f)$ surjects onto the $I_{f}$-invariant part of $V(f)$. It follows that $Z^{\co}\subset \poly {d}{I_{f}}{\ex N_{\ex F(\hat f)^{\s(f)}}}{V}\times_{\totl{\ex F(\hat f)^{\s(f)}}}\ex F^{\co}$ is smooth around $f$ and has codimension equal to the dimension of the $I_{f}$-invariant part of $V(f)$. 

Given any nice Whitney stratification of $\totl Z$, $\totl{Z^{\co}}$ is a union of strata because of the condition that every strata of $\totl Z$ projects within a strata of $(\ex F,G)$, and the fact that $\ex F^{\co}$ is a union of connected components of $\ex F^{\s(f)}$. As our normally complex section $s$ is transverse to $0$ in the sense of Definitions \ref{c def}, \ref{nc section def}, and \ref{Z transverse def},  the section $s$ lifts to a section $s^{\s(f)}$ of  $\poly {d}{I_{f}}{\ex N_{\ex F(\hat f)^{\s(f)}}}{V}$ which is transverse to some nice Whitney stratification of $\totl Z$. It follows that $s^{\s(f)}$ is transverse to $Z^{\co}$. As $\mathcal N_{\co}=(s^{\s(f)})^{-1}Z^{\co}/I_{f}$, it follows that $\mathcal N_{\co}$ is locally equal the quotient by $I_{f}$ of  the smooth family of curves $(s^{\s(f)})^{-1}Z^{\co}$. Note that $I_{f}$ acts trivially on the smooth part of $\ex F^{\co}$, which contains $(s^{\s(f)})^{-1}Z^{\co}$.

\

We shall now calculate the dimension of the smooth part of $(s^{\s(f)})^{-1}Z^{\co}$.  The codimension  of  $(s^{\s(f)})^{-1}Z^{\co}$ within   $\ex F^{\co}$ is equal to the dimension of the $I_{f}$-invariant part of $V(f)$.  The tangent space at $f$ of $\ex F^{\co}\subset \ex F(\hat f)$ is equal to the span of the space of $I_{f}$-invariant vectors and the $\mathbb R$-nil vectors. If $f$ is holomorphic, then the index of $D\dbar_{f}$ is equal to $T_{f}\ex F-V(f)$. Note that $T_{f}\ex F\subset \Msw$ contains all the $\mathbb R$-nil vectors in $T_{f}\Msw$. Therefore the dimension of $(s^{\s(f)})^{-1}Z^{\co}$ is equal to the index of the $I_{f}$-invariant part of $D\dbar_{f}$ plus the dimension of the $I_{f}$-anti invariant $\mathbb R$-nil vectors. We arrive at the dimension of the smooth part of $(s^{\s(f)})^{-1}Z^{\co}$ by quotienting out all the $\mathbb R$-nil vectors. Therefore the dimension of the smooth part of $(s^{\s(f)})^{-1}Z^{\co}$ is the index of the $I_{f}$-invariant part of $D\dbar_{f}$ minus the dimension of the $I_{f}$-invariant $\mathbb R$-nil vectors in $T_{f}\Msw$.

In the case that $f$ is not holomorphic, $D\dbar_{f}$ is not well defined, but there is a well defined linearization of $\dbar$ followed by a projection $\pi_{V}:\Y\longrightarrow \Y/V$. The index of $D\pi_{V}\dbar_{f}$ is equal in real $I_{f}$-equivariant $K$-theory to the right hand side of equation \ref{index} plus $V(f)$ plus the the tangent space to $\ex B_{0}$. On the other hand, $D\pi_{V}\dbar_{f}$ is surjective and has kernel equal to $T_{f}\ex F$. It follows that $d(\co)$ satisfies the required formula, and does not depend of the choice of $f\in\Msw_{\co}$ because the right hand side of equation \ref{index} is independent of choice of $f\in \Msw_{\co}$.

\stop

\begin{remark}In the proof above, it was shown that if $s$ is transverse to $0$ in the sense of Definitions \ref{c def}, \ref{nc section def}, and \ref{Z transverse def}, then in each Kuranishi chart $(\mathcal U,V,\hat f/G)$, restricted to $\ex F^{\co}\subset \ex F(\hat f)$, the corresponding section of $V^{I_{\co}}$ is transverse to $0$. The normal bundle to $\ex F^{\co}\subset \ex F(\hat f)$ is given a complex structure when defining the normally complex section $s$, therefore $\ex F^{\co}$ has a canonical orientation relative to $\ex F(\hat f)$, which is itself canonically oriented using the homotopy  of $D\dbar$ to a complex operator in the case of a single target $\ex B$, and is canonically oriented  relative to $\ex B_{0}$ in the case of a family $\hat {\ex B}\longrightarrow \ex B_{0}$ of targets. $V^{I_{\co}}$ is complex, therefore $\mathcal N_{\co}$ inherits a canonical orientation relative to $\ex B_{0}$. Condition \ref{c3} of Definition \ref{c def},  condition \ref{rk3} of Definition \ref{rks def} and the definition of an embedded Kuranishi structure imply that this orientation relative to $\ex B_{0}$ does not depend on the choice of Kuranishi chart used because the derivative of $s$ normal to a small dimension Kuranishi chart within a larger dimensional Kuranishi chart is positively oriented. 
\end{remark}

\begin{lemma}\label{codim2} The closure of $\mathcal N_{\co}\subset \Msw$ consists of the union of $\mathcal N_{\co}$ with a locally finite collection of substacks of $\mathcal N\subset\Msw$ which are smooth and have smooth parts with dimension at most $d(\co)-2$.

\end{lemma}

 Lemma \ref{cs} states that $\mathcal N:=s^{-1}(0)\subset \Msw$ is a closed substack, therefore the closure of of $\mathcal N_{\co}$ is contained within $\mathcal N$. Let $f$ be a curve in the closure of $\mathcal N_{\co}$. As in the proof of Lemma \ref{component dimension}, we may work locally in a Kuranishi chart $(\mathcal U,V,\hat f/G)$ containing $f$. Choosing $d$ large enough, consider $Z\subset\poly {d}{I_{f}}{\ex N_{\ex F(\hat f)^{\s(f)}}}{V}\times_{\totl{\ex F(\hat f)^{\s(f)}}} \ex N_{\ex F(\hat f)^{\s(f)}}$. Denote by $\ex N^{\co}$ the subset of $\ex N_{\ex F(\hat f)^{\s(f)}}$ corresponding to curves in $\Msw_{\co}$, and let $Z^{\co}$ be the subset of $Z$ which projects to $\ex N^{\co}$. Note that both $\totl{\ex N_{\ex F(\hat f)^{\s(f)}}}$ and the smooth part of $\poly {d}{I_{f}}{\ex N_{\ex F(\hat f)^{\s(f)}}}{V}\times_{\totl{\ex F(\hat f)^{\s(f)}}} \ex N_{\ex F(\hat f)^{\s(f)}}$ are bundles over the smooth manifold $\totl{\ex F(\hat f)^{\s(f)}}$ with fibers which are complex algebraic varieties. Each fiber of $\ex N^{\co}$ is a Zariski open subset of a complex subvariety, therefore the same is true of $\totl{Z^{\co}}$. In particular, each fiber of  $\overline{\totl{Z^{\co}}}$ is a complex subvariety, and $\totl Z^{\co}$ is a Zariski open subset. Given any nice Whitney stratification of $\totl Z$, $\totl{Z^{\co}}$ is a union of strata, therefore, $\overline{\totl {Z^{\co}}}$ is also a union of strata. Because on each fiber, $\totl{ Z^{\co}}$ is a Zariski open subset of $\overline{\totl{Z^{\co}}}$, the strata in $\overline{\totl{Z^{\co}}}\setminus\totl{Z^{\co}} $ have dimension at most $\dim \totl {Z^{\co}}-2$. 
 
 As $s$ is transverse to $0$, it lifts to a section $s^{\s(f)}$ of $\poly {d}{I_{f}}{\ex N_{\ex F(\hat f)^{\s(f)}}}{V}$ with graph transverse to every strata of a nice Whitney stratification of $\totl Z$. As proved in Lemma \ref{component dimension}, $s^{\s(f)}$ is transverse to $Z^{\co}$, and the smooth part of  $(s^{\s(f)})^{-1}Z^{\co}$ has dimension $d(\co)$. As every strata in $\overline{\totl{Z^{\co}}}\setminus\totl{Z^{\co}}$ has dimension at most $(\dim \totl{Z^{\co}}-2)$, the smooth part of the inverse image of each of these strata under $s^{\s(f)}$ has dimension at most $d(\co)-2$. Around $f$, the closure of $\mathcal N_{\co}$ is equal to $(s^{\s(f)})^{-1}\overline{Z^{\co}}/I_{f}$, so $\overline{\mathcal N_{\co}}\setminus \mathcal N_{\co}$ locally consists of a finite union of substacks of $\mathcal N$ which are smooth and have smooth parts of dimension at most $d(\co)-2$.
 
 \stop

\

\section{Defining integer curve counts in symplectic manifolds}\label{integer section}

Let $X$ be a compact symplectic manifold.  In Proposition \ref{nc prop}, we proved that there exists a normally complex perturbation $s$ of $\dbar$ which is transverse to the zero section. We also proved  that any two such normally complex perturbations are cobordant through a normally complex perturbation of $\dbar$ which is also transverse to the zero section. (Normally complex sections are said to be transverse to $0$ if they satisfy Definition \ref{Z transverse def}, which differs from the usual definition of transversality.) 

Let $\mathcal N$ indicate $s^{-1}(0)\subset\Msw(X)$. Lemma \ref{component dimension}  allows us to break up $\mathcal N\subset\Msw(X)$ into 
components $\mathcal N_{\co}$ which are smooth, oriented,  have smooth parts of dimension $d(\co)$. Lemma \ref{codim2} tells us that each $\mathcal N_{\co}$  may be compactified to $\overline{\mathcal N_{\co}}\subset\Msw(X)$ by adding a finite number of strata with codimension at least $2$. Recall that in the case that $X$ is smooth, $\Msw_{\co}(X)$ indicates a connected component of the moduli stack of curves with smooth parts that have some fixed number of nodes, and some fixed automorphism group. 

Of most interest is the case when there are no nodes and the automorphism group is trivial--- In this case $d(\co)$ is the usual expected dimension of the moduli stack of holomorphic curves. Given a curve $f$ with no nodes and genus $g$, a choice of oriented diffeomorphism from a fixed genus $g$ surface $\Sigma_{g}$ to the domain, $\ex C(f)$, of $f$ determines a homotopy class in $[\Sigma_{g},X]$. A different choice of oriented diffeomorphism $\Sigma_{g}\longrightarrow \ex C(f)$ changes this homotopy class by the action of an oriented diffeomorphism of $\Sigma_{g}$, so $f$ determines a well defined class 
\[[f]\in [\Sigma_{g},X]/MCG(\Sigma_{g})\]
where $MCG(\Sigma_{g})$ indicates the mapping class group of $\Sigma_{g}$ consisting of the quotient of the group of oriented diffeomorphisms of $\Sigma_{g}$ by the connected component containing the identity. Two curves $f$ and $f'$ with genus $g$ and no nodes or marked points are in the same component $\Msw_{\co}$ if and only if $[f]=[f']\in [\Sigma_{g},X]/MCG(\Sigma_{g})$.

 For the rest of this section, suppose that $\Msw_{\co}(X)$ is the moduli stack of $\C\infty1$ curves in $X$ with no automorphisms,  no nodes or marked points,  with genus $g$, and representing a homotopy class $\beta\in [\Sigma_{g},X]/MCG(\Sigma_{g})$. We shall use $(g,\beta)$ in place of $\co$ in this case.   The corresponding moduli space  $\mathcal N_{g,\beta}$ is a smooth, oriented family of curves in $\Msw_{g,\beta}(X)$ of dimension $d(g,\beta)$  equal to the usual expected dimension of the moduli stack of holomorphic curves.

\

If the expected dimension $d(g,\beta)=0$, then the signed number of curves $n_{g,\beta}$ in $\mathcal N_{g,\beta}$ is a symplectic invariant, because if $s_{0}$ and $s_{1}$ are normally complex perturbations of $\dbar$ transverse to $0$ (for possibly different complex structures tamed by the symplectic form on $X$), then $s_{0}$ and $s_{1}$ are joined by a normally complex perturbation $s_{t}$ of $\dbar$ on $\Msw(X\times[0,1])$ which is also transverse to $0$. Then  $s_{t}^{-1}(0)\cap \Msw_{g,\beta}$ is oriented relative to $[0,1]$, compact, one dimensional,  and has boundary at $0$ equal to $s^{-1}_{0}(0)\cap \Msw_{g,\beta}$, and boundary at $1$ equal to $s^{-1}_{1}(0)\cap \Msw_{g,\beta}$. In other words $\mathcal N_{g,\beta}$ defined using $s_{0}$ is cobordant to $\mathcal N_{g,\beta}'$ defined using $s_{1}$.

\

If the expected dimension $d(g,\beta)$ is greater than $0$,  we can count the number of curves in $\mathcal N_{g,\beta}$ passing through some generic, compact, oriented submanifolds  $\alpha_{1}, \dotsc ,\alpha_{k}$. If the sum of the codimensions of $\alpha_{i}$ is $d(g,\beta)+2k$, then define
\[n_{g,\beta}(\alpha_{1},\dotsc,\alpha_{k})\]
to be the signed number of curves  in $\mathcal N_{g,\beta}$ with points $x_{i}$ sent to $\alpha_{i}$. We shall define $n_{g,\beta}(\alpha_{1},\dotsc,\alpha_{k})$ more precisely, and check that it is a symplectic invariant below in the more general situation that $\alpha_{i}$ are the cohomology classes in $H^{*}(X,\mathbb Z)$  Poincare dual to the manifolds.

 Recall from \cite{evc} that there is a functorial way to add $k$ extra punctures to a family of curves $\hat f$ to obtain a family of curves $\hat f^{+k}$. Evaluation at these marked points defines a natural map
\[\Phi:(\Msw)^{+k}\longrightarrow X^{k}\]
The boundary,  $\overline{\mathcal N_{g,\beta}}\setminus \mathcal N_{g,\beta}$ of $\mathcal N_{g,\beta}$ is compact, and Lemma \ref{codim2} gives that it consists of a finite union of substacks $D_{i}$ with smooth parts that have dimension at most $d(g,\beta)-2$. In particular, this implies that $\Phi(D_{i}^{+k})\subset X^{k}$ will be the image under a smooth map of a smooth manifold with dimension at most $d(g,\beta)+2(k-1)$. Because the image under $\Phi$ of all these $D_{i}^{+k}$ is compact, a generic submanifold of $X^{k}$ with codimension greater than $d(g,\beta)+2(k-1)$ will not intersect $\Phi(D_{i}^{+k})$.

 Let $\alpha_{1},\dotsc,\alpha_{k}\in H^{*}(X,\mathbb Z)$ be pure dimensional singular cohomology classes with sum of dimensions equal to $d(g,\beta)-2k$. Abuse notation a little and think of $\alpha_{j}$ as living in $H^{*}(X^{k},\mathbb Z)$ by pulling it back via the $j$th projection map.  Represent the Poincare dual to $\alpha_{1}\cup\dotsb \cup\alpha_{k}$ by by a cycle $c$ which is smooth and generic in the sense that   it is a formal sum  smooth simplices which are transverse to $\Phi(\mathcal N_{g,\beta}^{+k})$ and do not intersect $\Phi(D^{+k}_{i})$.

Define $n_{g,\beta}(\alpha_{1},\dotsc,\alpha_{k})$ to be the oriented count of the intersection points of $\Phi(\mathcal N_{g,\beta}^{+k})$ with $c$. 

 Any  generic,  smooth, codimension $d(g,\beta)+2k$ cycle  representing $0$ in homology is the boundary of a codimension $d(g,\beta)+2k-1$ chain which does not intersect $\Phi(D^{+k}_{i})$, and intersects $\Phi(\mathcal N_{g,\beta}^{+k})$ in a $1$-chain with boundary the intersection of $\Phi(\mathcal N_{g,\beta}^{+k})$ with the original cycle. Therefore $n_{g,\beta}(\alpha_{1},\dotsc,\alpha_{k})$ does not depend on the choice of a generic smooth cycle $c$ Poincare dual to $\alpha_{1}\cup\dotsb \cup\alpha_{k}$. 

It is similarly easy to show that $n_{g,\beta}(\alpha_{1},\dotsc,\alpha_{k})$ is a symplectic invariant. In particular, suppose that $\mathcal N_{g,\beta}$ and $\mathcal N_{g,\beta}'$ are defined using  normally complex perturbations $s_{0}$ and $s_{1}$ of $\dbar$ which are transverse to $0$ and  defined using possibly different choices of almost complex structure and so on. Then Proposition \ref{nc prop} gives that $s_{0}$ and $s_{1}$ extend to a normally complex perturbation $s_{t}$ of $\dbar$ on $\Msw(X\times \mathbb R)$ which is also transverse to $0$. It follows from \cite{cem} and Lemma \ref{cs} that the projection $s^{-1}_{t}(0)\longrightarrow \mathbb R$ is proper. Use the notation $ \mathcal N_{g,\beta}(t):=s_{t}^{-1}(0)\cap \Msw_{g,\beta}$. As $ \mathcal N_{g,\beta}(t)$ is naturally oriented relative to  $\mathbb R$,  we may orient it by giving $\mathbb R$ the standard orientation.  If $c_{0}$ and $c_{1}$ are generic codimension $d(g,\beta)+2k$ smooth cycles in $X^{k}\times \{0\} $ and $X^{k}\times{1}$ respectively and $c_{0}-c_{1}$ represents $0$ in $H_{*}(X^{k}\times\mathbb R,\mathbb Z)$, then there is a smooth chain $\tilde c$ with boundary $c_{0}-c_{1}$ so that $\tilde  c$ is in general position with respect to $\Phi(\overline{\mathcal N_{g,\beta}(t)^{+k}})$. In particular, for dimension reasons, $\tilde c$ intersects $\Phi(\overline{\mathcal N_{g,\beta}(t)^{+k}})$ only in $\Phi(\mathcal N_{g,\beta}(t)^{+k})$, and intersects $\Phi(\mathcal N_{g,\beta}(t)^{+k})$ in a $1$ chain with boundary equal to $c_{0}\cap\Phi(\mathcal N_{g,\beta}^{+k})-c_{1}\cap \Phi((\mathcal N'_{g,\beta})^{+k})$. In particular, the integers $n_{g,\beta}(\alpha_{1},\dotsc,\alpha_{k})$ are symplectic invariants that do not depend on the choices involved in their construction. 

\begin{remark} An alternative way to define the integers $n_{g,\beta}(\alpha_{1},\dotsc,\alpha_{n})$ is using integration of differential forms. Use the notation \[\begin{tikzcd}\ex C\dar{\pi}\rar {\hat f}&X
\\ \ex F\end{tikzcd}\]  for the family of curves representing $\mathcal N_{g,\beta}$. Represent each $\alpha_{i}$ as a differential form. Then
\[n_{g,\beta}(\alpha_{1},\dotsc, \alpha_{k})=\int_{\ex F}\pi_{!}\hat f^{*}\alpha_{1}\wedge\dotsb\wedge \pi_{!}\hat f^{*}\alpha_{k}\]
\end{remark}

Similar constructions can be made which count smooth, automorphism-free curves with genus $g$ with $k$ marked points, which may be labeled, unlabeled, or only allow morphisms which act on labels by a specified subgroup of $S_{k}$. The resulting invariants will be different from the invariants defined above because adding marked points can kill the automorphisms of a curve. The invariants suggested by Fukaya and Ono in \cite{FOinteger} used $k$ labeled marked points.

\section{Generalizations and suggestions for further research}\label{future section}

\subsection{Gravitational descendants}

The invariants $n_{g,\beta}(\alpha_{1},\dotsc,\alpha_{k})$ may be easily generalized to include a version of the gravitational descendants of Gromov-Witten theory. There are tautological line bundles $\mathcal L_{1},\dotsc,\mathcal L_{k}$ defined over $(\Msw)^{+k}$ as follows: given a family of curves $\hat f$, forgetting the $i$th of the $k$ extra punctures in $\hat f^{+k}$ defines a map $\pi_{i}:\ex F(\hat f^{+n})\longrightarrow \ex F(\hat f^{+n-1})$. This forgetful map has the natural structure of a family of curves equal to $\ex C(\hat f^{+n-1})\longrightarrow \ex F(\hat f^{+n-1})$. On $\ex F(\hat f^{+k})$, we can define the complex line bundle 
\[\mathcal L_{i}:=\ker\pi_{i}\subset T\ex F(\hat f^{+n})\]

Recall that $\mathcal N$ indicates the intersection of a normally complex perturbation of $\dbar$ with $0$, and $\mathcal N^{+k}$ indicates the moduli stack obtained from $\mathcal N$ by adding $k$ extra points to curves in $\mathcal N$. It is easy to define a notion of a normally complex section of $\mathcal L_{i}$ over  $\mathcal N^{+k}$. Given multiplicities $a_{i}\in\mathbb N$, a generic normally complex section $v$ of $\oplus_{i}\mathcal L_{i}^{\oplus a_{i}}$ will intersect the zero section in a set with a nice stratification. Let $\mathcal N_{g,\beta}(X)$ be the subset of $\mathcal N$  consisting of smooth curves with no automorphisms, genus $g$ and representing the homotopy class $\beta$. The subset of $\mathcal N_{g,\beta}^{+k}$ where $v$ is $0$ will have dimension $d(g,\beta)+2\sum_{i}(1-a_{i})$, and will have boundary with smooth part consisting of some union of codimension $2$ manifolds.  Then integer invariants
\[n_{g,\beta}(\tau_{a_{1}}(\alpha_{1}),\dotsc,\tau_{a_{k}}(\alpha_{k}))\]
may be defined by representing the Poincare dual of $a_{1}\cup \dotsb \cup a_{k}$ by a generic smooth cycle $c$ and counting the signed number of points in the intersection of $\Phi^{-1}(c)$ with the subset of $\mathcal N_{g,\beta}^{+k}$ where $v$ is $0$.

\

Similar constructions may be made using any naturally defined complex vector bundles over $\Msw(X)^{+k}$ in place of $\mathcal L_{i}$. For example, a complex vector bundle over $X$ such as $TX$ can be pulled back and tensored  with $\mathcal L_{i}$. Taking normally complex sections of the resulting vector bundle defines invariants.

\subsection{Other integer invariants}

The methods for obtaining integer invariants in this paper may be imitated using other sheaves of sections in place of normally complex sections. Intuitively, using normally complex sections allows us to `throw away' curves with automorphisms from the curve count. If a curve count which `throws away' the contribution of curves mapping to a point was desired, then one could stratify Kuranishi charts using strata which take note of when a component of a curve maps to a point. Then, using normally complex sections  which vanish on any component of a curve which maps to a point, we can define integer invariants which `throw away' contributions from curves with a component mapping to a point.  A comparable approach was taken by Zinger to define  reduced genus $1$ Gromov-Witten invariants in \cite{Zinger2}.

 More ambitiously, one could hope to   define an integer invariant which morally counts embedded holomorphic curves by `throwing away' contributions of holomorphic curves  which factor through simpler holomorphic curves. Such an integer invariant would require embedded Kuranishi structures that have a nice stratification which takes note of when curves factor through simpler curves. To define integer invariants which `throw away' factorizable curves, we also need to be able to choose a complex structure normal to these strata.  If a holomorphic curve $f$ factors through $f_{0}$, then  $D\dbar_{f}$  restricted to the tangent space  of the substack  of curves factoring through $f_{0}$ is complex, so the homotopy of $D\dbar$ to a complex operator determines a homotopy class of  complex structure normal to the curves which factor through $f_{0}$.  Therefore, it is reasonable to expect that some kind of normally complex sections may be defined and used to define integer invariants which `throw away'  factorizable curves.  A partial result in this direction has been achieved by Zinger in \cite{Zinger1}, which compares the contribution to Gromov-Witten invariants of curves which factorize through a given curve $C$ with the Gromov-Witten invariants of $C$. In a recent preprint \cite{IPGV}, Ionel and Parker have given a proof of the Gopakumar-Vafa formula in the setting of $6$ dimensional symplectic Calabi Yau manifolds. In this paper, Ionel and Parker relate Gromov-Witten invariants to a virtual count of elementary curves. It is likely that the integer invariants which throw away factorizable curves using a normally complex perturbation coincide in this case with Ionel and Parker's  virtual count of elementary curves. 

Similarly, a kind of normally complex perturbation may be used to `throw away' contributions from curves with components in a closed $J$-holomorphic submanifold $M$ of $X$. If we do not also throw away the contribution of curves with automorphisms, the resulting invariant will be rational instead of integral. In the case that $M$ is codimension $2$,  this is one interpretation of the construction of relative Gromov-Witten invariants by Ionel and Parker in \cite{IP1}. 

\subsection{Floer homology over the integers}

If sphere bubbles with negative first Chern class can be avoided, then Floer homology is defined over the integers. (This is proved in \cite{Floer}, \cite{HS}, and \cite{Ono}.) In general symplectic manifolds, sphere bubbles with automorphisms make it difficult to define Floer homology over the integers, although, as shown in \cite{FO} and \cite{Liu-Tian},  Floer homology may still be defined over $\mathbb Q$. 

In \cite{FOinteger}, Fukaya and Ono sketched a way of defining Floer homology over the integers, which is essentially the method in this paper  applied to Floer trajectories which allow sphere bubbles. Although \cite{evc} does not apply directly in this situation, very similar analysis gives an embedded Kuranishi structure on the moduli stack of Floer trajectories. Around a Floer trajectory with spherical bubbles, Kuranishi charts may be modeled on some $(\et 1{[0,\infty)})^{n}$-bundle over a fiber product of a moduli space of  Floer trajectories (with no bubbles), and a Kuranishi chart on the moduli stack of spheres. (The $(\et 1{[0,\infty)})^{n}$ bundle may be thought of as gluing data for smoothing $n$ nodes.) In such charts, normally complex sections still make sense, so we can take a generic normally complex perturbation of $\dbar$. Floer homology can then be defined using the resulting perturbed Floer trajectories that do not have any attached bubbles. 

The analysis of whether a Floer trajectory attached to a sphere bubble contributes to the count of these perturbed Floer trajectories depends on analysis of the $\dbar$ equation in the $(\et 1{[0,\infty)})^{n}$ `gluing' directions. Essentially, if $\dbar$ is always complex and nondegenerate in these directions, then such trajectories with bubbles should not contribute.  An optimist might hope that such contributions from bubbled trajectories disappear in generic situations. (Note that the optimistic situation is in contrast to the case of our integer invariants--- even in generic situations, curves with automorphisms may contribute to our count of curves with no automorphisms after $\dbar$ is perturbed to be normally complex. Examples to show that this must be the case are constructed in \cite{taubessw2}.)

\subsection{Invariants counting curves with nodes and automorphisms }
Recall that $\Msw_{\co}(X)$ consists of a connected component of the moduli stack of curves  with smooth parts that have some fixed number of nodes, and some fixed automorphism group. Lemmas \ref{component dimension} and \ref{codim2} give that that the corresponding  component $\mathcal N_{\co}$ of $s^{-1}(0)$  has smooth part of dimension $d(\co)$ and has a compactification which adds in some  boundary with codimension at least $2$. Integer invariants may therefore be defined which count curves in $\mathcal N_{\co}$ in a similar way to the invariants $n_{g,\beta}(\alpha_{1},\dotsc,\alpha_{k})$.

 The curves in $\Msw_{\co}(X)$ can be thought of as created from simple pieces by taking appropriate multiple covers and gluing in a way specified by $\co$. Therefore, it may be expected that all these integer invariants may be deduced from the simple invariants corresponding to curves with no nodes and no automorphisms. Such a deduction is not straight forward for the following reason:  The gluing may be thought of as taking place over some graph $\gamma$ with decorated vertices. The `simple' invariants that are needed are invariants of curves with possibly disconnected domains and marked points labelled appropriately by oriented edges   of $\gamma$, where automorphisms of such labeled curves are allowed to act on labels by an automorphism of $\gamma$. The action of the automorphism group of $\gamma$ affects what is meant by normally complex, so these `simple' invariants may not be easily deducible from invariants which count connected curves with no nodes and no automorphisms. 

\subsection{Relationship with Gromov-Witten invariants}
There is a good possibility that Gromov-Witten invariants of a symplectic manifold $X$ may be reconstructed from the information in the components $\mathcal N_{\co}$ of $\mathcal N$. In particular, on the strata corresponding to $\co$, $s$ is transverse to $0$ in $I_{\co}$-invariant directions, and is `normally complex' in the other directions.  Modification of $s$ to a weighted branched section which is transverse to $0$ and agrees with $s$ in a neighborhood should define a contribution of $\mathcal N_{\co}$ to Gromov-Witten invariants which depends only on bundle data corresponding to these `other directions'. This bundle data is naturally associated with $\mathcal N_{\co}\subset \Msw$, but in general may contain more information than simply $\co$ and the integer invariants described in the previous section.

\subsection{Invariants in the exploded manifold setting}

In the case of curves in an exploded manifold $\ex B$ for which Gromov-compactness holds, integer invariants may be defined corresponding to each component of $\mathcal N_{\co}$. These integer invariants  are invariant under deformations of $\ex B$ parametrized by $\mathbb R$, but are  not invariant under deformations of $\ex B$ in a family of exploded manifolds. It would be desirable to have a version of the invariants $n_{g,\beta}(\alpha_{1},\dotsc,\alpha_{k})$ which do not change in families of exploded manifolds so that the invariants of symplectic manifolds may be calculated using tropical curve counts and the corresponding gluing formulas.

\bibliographystyle{plain}
\bibliography{ref.bib}

\end{document}